\magnification=1200
\def\nop{}
\def\dt{}

\def\title#1{{\titlefont\noindent #1\bigskip}}

\def\author#1{{\largefont\noindent #1}\medskip}

\def\beginlinemode{\endmode
 \begingroup\obeylines\def\endmode{\par\endgroup}}
\let\endmode=\par

\newbox\theaddress
\def\address{\smallskip\beginlinemode\parindent 0in\getaddress}
{\obeylines
\gdef\getaddress #1
 #2
 {#1\gdef\addressee{#2}%
   \global\setbox\theaddress=\vbox\bgroup\raggedright%
    \everypar{\hangindent2em}#2
   \def\endaddress{\egroup\endgroup \copy\theaddress \medskip}}}

\def\thanks#1{\footnote{}{\eightpoint #1}}

\long\def\Abstract#1{{\eightpoint\narrower\vskip\baselineskip\noindent
#1\smallskip}}

\def\skipfirstword#1 {}

\def\ir#1{\csname #1\endcsname}

\def\irrnSection#1#2{\edef\tttempcs{\ir{#2}}
\vskip\baselineskip\penalty-3000
{\largefont\bf\noindent \expandafter\skipfirstword\tttempcs. #1}
\vskip6pt}

\def\irSubsection#1#2{\edef\tttempcs{\ir{#2}}
\vskip\baselineskip\penalty-3000
{\bf\noindent \expandafter\skipfirstword\tttempcs. #1}
\vskip6pt}

\def\irSubsubsection#1#2{\edef\tttempcs{\ir{#2}}
\vskip\baselineskip\penalty-3000
{\noindent \expandafter\skipfirstword\tttempcs. #1}
\vskip6pt}

\def\References{\vskip6pt\noindent{}{\bf References}
\vskip6pt\par}

\def\baselinebreak{\par \ifdim\lastskip<6pt
         \removelastskip\penalty-200\vskip6pt\fi}

\long\def\prclm#1#2#3{\baselinebreak
\noindent{\bf \csname #2\endcsname}:\enspace{\sl #3\par}\baselinebreak}

\def\Prf{\noindent{\bf Proof}: }

\def\rem#1#2{\baselinebreak\noindent{\bf \csname #2\endcsname}:}

\def\qed{{$\diamondsuit$}\vskip6pt}

\nop{}
\def\bibitem#1{\par\indent\llap{\rlap{\bf [#1]}\indent}\hskip3em\hangindent
5em\ignorespaces}

\long\def\eatit#1{}

\def\leftheadlinetext{}
\def\rightheadlinetext{}

\def\leftheadline{{\eightrm\folio\hfil \leftheadlinetext\hfil}}
\def\rightheadline{{\eightrm\hfil\rightheadlinetext\hfil\folio}}

\headline={\ifnum\pageno=1\hfil\else
\ifodd\pageno\rightheadline\else\leftheadline\fi\fi}

\def\tenpoint{\def\rm{\fam0\tenrm}
\textfont0=\tenrm \scriptfont0=\sevenrm \scriptscriptfont0=\fiverm
\textfont1=\teni \scriptfont1=\seveni \scriptscriptfont1=\fivei
\def\mit{\fam1} \def\oldstyle{\fam1\teni}
\textfont2=\tensy \scriptfont2=\sevensy \scriptscriptfont2=\fivesy
\def\cal{\fam2}
\textfont3=\tenex \scriptfont3=\tenex \scriptscriptfont3=\tenex
\def\it{\fam\itfam\tenit} 
\textfont\itfam=\tenit
\def\sl{\fam\slfam\tensl} 
\textfont\slfam=\tensl
\def\bf{\fam\bffam\tenbf} 
\textfont\bffam=\tenbf \scriptfont\bffam=\sevenbf
\scriptscriptfont\bffam=\fivebf
\def\tt{\fam\ttfam\tentt} 
\textfont\ttfam=\tentt
\normalbaselineskip=12pt
\setbox\strutbox=\hbox{\vrule height8.5pt depth3.5pt  width0pt}%
\normalbaselines\rm}

\def\eightpoint{\def\rm{\fam0\eightrm}%
\textfont0=\eightrm \scriptfont0=\sixrm \scriptscriptfont0=\fiverm
\textfont1=\eighti \scriptfont1=\sixi \scriptscriptfont1=\fivei
\def\mit{\fam1} \def\oldstyle{\fam1\eighti}%
\textfont2=\eightsy \scriptfont2=\sixsy \scriptscriptfont2=\fivesy
\def\cal{\fam2}%
\textfont3=\tenex \scriptfont3=\tenex \scriptscriptfont3=\tenex
\def\it{\fam\itfam\eightit} 
\textfont\itfam=\eightit
\def\sl{\fam\slfam\eightsl} 
\textfont\slfam=\eightsl
\def\bf{\fam\bffam\eightbf} 
\textfont\bffam=\eightbf \scriptfont\bffam=\sixbf
\scriptscriptfont\bffam=\fivebf
\def\tt{\fam\ttfam\eighttt} 
\textfont\ttfam=\eighttt
\normalbaselineskip=9pt%
\setbox\strutbox=\hbox{\vrule height7pt depth2pt  width0pt}%
\normalbaselines\rm}

\def\largefont{\def\rm{\fam0\largerm}
\textfont0=\largerm \scriptfont0=\largescriptrm \scriptscriptfont0=\tenrm
\textfont1=\largei \scriptfont1=\largescripti \scriptscriptfont1=\teni
\def\mit{\fam1} \def\oldstyle{\fam1\teni}
\textfont2=\largesy 
\def\cal{\fam2}
\def\it{\fam\itfam\largeit} 
\textfont\itfam=\largeit
\def\sl{\fam\slfam\largesl} 
\textfont\slfam=\largesl
\def\bf{\fam\bffam\largebf} 
\textfont\bffam=\largebf 
\scriptscriptfont\bffam=\fivebf
\def\tt{\fam\ttfam\largett} 
\textfont\ttfam=\largett
\normalbaselineskip=17.28pt
\setbox\strutbox=\hbox{\vrule height12.25pt depth5pt  width0pt}%
\normalbaselines\rm}

\def\titlefont{\def\rm{\fam0\titlerm}
\textfont0=\titlerm \scriptfont0=\largescriptrm \scriptscriptfont0=\tenrm
\textfont1=\titlei \scriptfont1=\largescripti \scriptscriptfont1=\teni
\def\mit{\fam1} \def\oldstyle{\fam1\teni}
\textfont2=\titlesy 
\def\cal{\fam2}
\def\it{\fam\itfam\titleit} 
\textfont\itfam=\titleit
\def\sl{\fam\slfam\titlesl} 
\textfont\slfam=\titlesl
\def\bf{\fam\bffam\titlebf} 
\textfont\bffam=\titlebf 
\scriptscriptfont\bffam=\fivebf
\def\tt{\fam\ttfam\titlett} 
\textfont\ttfam=\titlett
\normalbaselineskip=24.8832pt
\setbox\strutbox=\hbox{\vrule height12.25pt depth5pt  width0pt}%
\normalbaselines\rm}

\nopagenumbers

\font\eightrm=cmr8
\font\eighti=cmmi8
\font\eightsy=cmsy8
\font\eightbf=cmbx8
\font\eighttt=cmtt8
\font\eightit=cmti8
\font\eightsl=cmsl8
\font\sixrm=cmr6
\font\sixi=cmmi6
\font\sixsy=cmsy6
\font\sixbf=cmbx6

\font\largerm=cmr12 at 17.28pt
\font\largei=cmmi12 at 17.28pt
\font\largescriptrm=cmr12 at 14.4pt
\font\largescripti=cmmi12 at 14.4pt
\font\largesy=cmsy10 at 17.28pt
\font\largebf=cmbx12 at 17.28pt
\font\largett=cmtt12 at 17.28pt
\font\largeit=cmti12 at 17.28pt
\font\largesl=cmsl12 at 17.28pt

\font\titlerm=cmr12 at 24.8832pt
\font\titlei=cmmi12 at 24.8832pt
\font\titlesy=cmsy10 at 24.8832pt
\font\titlebf=cmbx12 at 24.8832pt
\font\titlett=cmtt12 at 24.8832pt
\font\titleit=cmti12 at 24.8832pt
\font\titlesl=cmsl12 at 24.8832pt

\tenpoint

\nop{}

\newdimen\myx
\newdimen\myht
\newdimen\myhtB
\newdimen\mydepth
\newdimen\xm
\newdimen\ym
\newdimen\xM
\newdimen\yM
\newdimen\mywidth
\newdimen\mylength
\newdimen\myheight
\newdimen\penwidth
\newdimen\penheight
\newdimen\myskip

\def\myplt #1 #2 {\myskip\baselineskip\baselineskip0pt{
\xm\xmin pt
\ym\ymin pt
\xM\xmax pt
\yM\ymax pt
\myx#1 pt
\myht#2 pt
\ifdim\myx<\xm \else \ifdim\myx>\xM \else \ifdim\myht<\ym \else \ifdim\myht>\yM \else 
\advance\myx by -\xm \advance\myht by -\ym 
\myx\xscale\myx\myht\yscale\myht 
\multiply\myht by-1\advance\myht by\myheight
\hbox to0in{\vbox to0in{\hbox to0in{\hbox to\myx{\hfil}\vbox{\vbox to\myht{\vss}\vrule depth\penheight width\penwidth\vss}\hss}\vss}\hss}
\fi\fi\fi\fi}\baselineskip\myskip}

\def\segplt #1 #2 #3 {\myskip\baselineskip\baselineskip0pt{
\xm\xmin pt
\ym\ymin pt
\xM\xmax pt
\yM\ymax pt
\mylength#1 pt
\myx#2 pt
\myht#3 pt
\penwidth\yscale\penwidth
\myhtB\mylength\advance\myhtB by\myht
\ifdim\myht<\ym \myht\ym\fi
\ifdim\myhtB>\yM \myhtB\yM\fi
\ifdim\myhtB<\myht \mylength0pt\else \mylength\myhtB\advance\mylength by-\myht\fi
\ifdim\myx<\xm \else \ifdim\myx>\xM \else \ifdim\myhtB<\ym \else \ifdim\myht>\yM \else 
\advance\myx by -\xm \advance\myhtB by -\ym 
\myx\xscale\myx\myhtB\yscale\myhtB\mylength\yscale\mylength
\multiply\myhtB by-1\advance\myhtB by\myheight
\hbox to0in{\vbox to0in{\hbox to0in{\hbox to\myx{\hfil}\vbox{\vbox to\myhtB{\vss}\vrule depth\mylength width\penwidth\vss}\hss}\vss}\hss}
\fi\fi\fi\fi}\baselineskip\myskip}


\newcount\mynumber
\def\plt #1 #2 #3 {\mynumber #1\advance\mynumber by1\segplt {\mynumber} #2 #3 }

\def\myput #1 #2 #3 {\myskip\baselineskip\baselineskip0pt{
\xm\xmin pt
\ym\ymin pt
\xM\xmax pt
\yM\ymax pt
\myx#1 pt
\myht#2 pt
\ifdim\myx<\xm \else \ifdim\myx>\xM \else \ifdim\myht<\ym \else \ifdim\myht>\yM \else 
\advance\myx by -\xm \advance\myht by -\ym 
\myx\xscale\myx\myht\yscale\myht 
\multiply\myht by-1\advance\myht by\myheight
\hbox to0in{\vbox to0in{\hbox to0in{\hbox to\myx{\hfil}\vbox{\vbox to\myht{\vss}#3\vss}\hss}\vss}\hss}
\fi\fi\fi\fi}\baselineskip\myskip}

\def\putinaxes #1 #2 {{\xm\xmin pt
\xM\xmax pt
\advance\xM by -\xm 
\xM\xscale\xM
\penwidth\xM
\myplt {\xmin} 0 
\vskip5pt\myput {\xmax} 0 {#1} 
\vskip-5pt}
{\ym\ymin pt\yM\ymax pt
\myput 0 {\ymax} {\hbox to0in{\hss #2\  }}
\advance\yM by -\ym 
\yM\yscale\yM
\penheight\yM
\vskip-\yM
\myplt 0 {\ymin} 
\vskip\yM}}

\def\putinticx #1 #2 #3 {{\penheight#1pt\myplt #2 0
\vskip\penheight\myput #2 0 #3  
\vskip-\penheight}}

\def\putinticy #1 #2 #3 #4 {{\penwidth#1pt\myplt {#2} #3
\myput 0 #3 {\hskip-\penwidth \hbox to0in{\hss #4\  }}
}}

\def\putincaption #1 #2 {{\vskip#1 pt\myput {\xmin} {\ymin} {\hbox{\hsize\mywidth\hfil\vbox{#2}\hfil}} \vskip-#1pt}}



\def\refAH{AH}
\def\refB{B}
\nop{}
\def\refBZ{BZ}
\def\refCat{C1}
\def\refCatb{C2}
\def\refCCMO{CCMO}
\def\refCM{CM1}
\def\refCMb{CM2}
\nop{}
\def\refCMc{CM3}
\def\refE{E1}
\def\refEb{E2}
\def\refFHH{FHH}
\def\refGGR{GGR}
\def\refGi{Gi}
\nop{}
\def\refGIda{GI}
\def\refvanc{H1}
\def\refigp{H2}
\def\refnagconj{H3}
\def\refsurv{H4}
\def\refSC{H5}
\def\refHHF{HHF}
\def\refHi{Hi1}
\def\refHib{Hi2}
\def\refHo{Ho}
\def\refI{I}
\def\refMig{Mg}
\def\refMirb{Mr}
\def\refN{N1}
\def\refNb{N2}
\def\refR{R1}
\def\refRb{R2}
\def\refS{S}
\def\refXb{X}



\def\item#1{\par\indent\indent\llap{\rlap{#1}\indent}\hangindent
2\parindent\ignorespaces}

\def\itemitem#1{\par\indent\indent
\indent\llap{\rlap{#1}\indent}\hangindent
3\parindent\ignorespaces}
\def\binom#1#2{\hbox{$#1 \choose #2$}}

\def\pr#1{\hbox{{\bf P}${}^{#1}$}}
\def\cite#1{[\ir{#1}]}

\def\leftheadlinetext{B. Harbourne and J. Ro\'e}
\def\rightheadlinetext{Linear systems in \pr2}

\title{Linear systems with multiple base}
\vskip-\baselineskip
\title{points in \pr2}
\vskip\baselineskip

\author{Brian Harbourne}

\address
Department of Mathematics and Statistics
University of Nebraska-Lincoln
Lincoln, NE 68588-0323
email: bharbour@math.unl.edu
WEB: http://www.math.unl.edu/$\sim$bharbour/
\smallskip
\endaddress

\author{Joaquim Ro\'e}

\address
Departament d'\`Algebra i Geometria
Universitat de Barcelona
Barcelona 08007, Spain
email: jroevell@mat.ub.es
\smallskip
\nop{}January 30, 2003\endaddress

\Abstract{Abstract: Conjectures for the Hilbert function
$h(n;m)$ and minimal free resolution
of the $m$th symbolic power $I(n;m)$ of the ideal of $n$
general points of \pr2 are verified for a broad range of
values of $m$ and $n$ where both $m$ and $n$ can be large,
including (in the case of the Hilbert
function) for infinitely many $m$ for each square $n>9$ and
(in the case of resolutions) for infinitely many $m$ for each
even square $n>9$. All previous results require either that
$n$ be small or be a square of a special form, or that $m$ be small compared to $n$.
Our results are based on a new approach for bounding
the least degree among curves passing through $n$ general points
of \pr2 with given minimum multiplicities at each point
and for bounding the regularity of the linear system
of all such curves. For simplicity, we work over the complex
numbers.}

\irrnSection{Introduction}{intro}
Consider the ideal $I(n;m)\subset R={\bf C}[\pr2]$ generated by
all forms having multiplicity at least $m$ at $n$ given general points
of \pr2. This is a graded ideal, and thus we can consider the
Hilbert function $h(n;m)$ whose value at each nonnegative integer
$t$ is the dimension $h(n;m)(t) = \hbox{dim} I(n;m)_t$
of the homogeneous component $I(n;m)_t$
of $I(n;m)$ of degree $t$. It is well known that
$h(n;m)(t)\ge \hbox{max}(0,\binom{t+2}{2}-n\binom{m+1}{2})$,
with equality for $t$ sufficiently large.
Denote by $\alpha(n;m)$ the least degree $t$ such that
$h(n;m)(t)>0$ and by $\tau(n;m)$ the least degree $t$
such that $h(n;m)(t)=\binom{t+2}{2}-n\binom{m+1}{2}$;
we refer to $\tau(n;m)$ as the {\it regularity\/} of $I(n;m)$.

For $n\le9$, the Hilbert function \cite{refNb} and minimal
free resolution \cite{refigp} of $I(n;m)$ are known. For
$n>9$, there are in general only conjectures:

\prclm{Conjecture}{introconj}{Let $n\ge10$ and $m\ge0$; then:
\item{(a)} $\alpha(n;m)\ge m\sqrt{n}$;
\item{(b)} $h(n;m)(t)=\hbox{max}(0,\binom{t+2}{2}-n\binom{m+1}{2})$
for each integer $t\ge0$; and
\item{(c)} the minimal free resolution of $I(n;m)$ is
an exact sequence
$$0\to R[-\alpha-2]^d\oplus R[-\alpha-1]^c\to
R[-\alpha-1]^b\oplus R[-\alpha]^a\to I(n;m)\to 0,$$
where $\alpha=\alpha(n;m)$,
$a=h(n;m)(\alpha)$,
$b=\hbox{max}(h(n;m)(\alpha+1)-3h(n;m)(\alpha),0)$,
$c=\hbox{max}(-h(n;m)(\alpha+1)+3h(n;m)(\alpha),0)$,
$d=a+b-c-1$, and $R[i]^j$ is the direct sum of $j$ copies of
the ring $R={\bf C}[\pr2]$, regarded as an $R$-module
with the grading $R[i]_k=R_{k+i}$.}

Note that \ir{introconj}(c) implies \ir{introconj}(b)
which implies \ir{introconj}(a). \ir{introconj}(a)
was posed in \cite{refN} in the form ``$\alpha(n;m)>m\sqrt{n}$
for $n\ge10$ and $m>0$'', together with a proof in case
$n>9$ is a square. \ir{introconj}(c) was posed in \cite{refigp},
together with a determination of the resolution for $n\le9$.
\ir{introconj}(b) is a special case
of more general conjectures posed in different but equivalent
forms by a number of people. In particular, given
general points $p_1,\ldots,p_n\in\pr2$,
let ${\bf m}=(m_1,\ldots,m_n)$ be any sequence of nonnegative integers,
and define $I({\bf m})$ to be the ideal generated by
all forms having multiplicity at least $m_i$ at $p_i$.
We can in the obvious and analogous way define $\alpha({\bf m})$,
$\tau({\bf m})$ and $h({\bf m})$. Equivalent conjectures for $h({\bf m})$ have been
posed in \cite{refvanc}, \cite{refGi} and \cite{refHi}. Ciliberto and
Miranda have recently pointed out that these conjectures are
also equivalent to what seemed to be a weaker conjecture posed
in \cite{refS}; see \cite{refCMc}, or \cite{refsurv}.
However, no general conjecture for the minimal free resolution of $I({\bf m})$
has yet been posed when ${\bf m}$ is arbitrary.

These conjectures for $h({\bf m})$ have been verified in certain
special cases: for $n\le 9$ by
\cite{refNb}; for any $n$ as long as $m_i\le 4$ for
all $i$ by \cite{refMig}; and by \cite{refAH}
for any $m_i$ as long as the maximum of the $m_i$
is sufficiently small compared to the number of
points for which $m_i>0$. In addition,
\ir{introconj}(b) was shown to be true by \cite{refHib}
for $m\le3$, and by \cite{refCMb}
for $m\le12$. (After our paper was submitted for publication,
this was extended
to $m\le20$ by \cite{refCCMO}.) The only result before ours with
$n\ge10$ and arbitrarily large multiplicities was 
that of \cite{refEb}, which verifies \ir{introconj}(b)
for all $m$ as long as $n$ is a power of 4. (This has now
been extended to include all $n$ which are products of powers
of 4 and 9; see \cite{refBZ}.) 
Similarly, \cite{refE} verifies \ir{introconj}(a)
as long as $m$ is no more than about $\sqrt{n}/2$.

Results for resolutions are more limited. A complete solution
for the resolution of $I({\bf m})$ for arbitrary ${\bf m}$
was given in \cite{refCatb} as long as $n\le 5$.
This has been extended to $n\le 8$ by \cite{refFHH}. 
The resolution of $I(n;m)$ was found  
for $n\le 9$ for any $m$ by \cite{refigp}. Also,
\cite{refGGR} shows that \ir{introconj}(c) 
holds for $m=1$. This was extended to $m=2$
by \cite{refI} and to $m=3$ by \cite{refGIda}.
In addition, \cite{refHHF} applies
the result of \cite{refEb} to verify
\ir{introconj}(c) when $n$ is a power of 4,
as long as $m$ is not too small.

In this paper we obtain substantial improvements on
these prior results for all three parts of \ir{introconj}.
For example, we have:

\prclm{Corollary}{introcor}{Let $n\ge10$. Then:
\item{(a)} \ir{introconj}(a) holds as long as $m\le (n-5\sqrt{n})/2$;
\item{(b)} \ir{introconj}(b) holds
for infinitely many $m$ for each square $n\ge 10$; and
\item{(c)} \ir{introconj}(c) holds for infinitely many $m$
for each even square $n\ge10$.}

We also verify \ir{introconj}(a,b,c) for many other values
of $m$ and $n$ (see \ir{vernag}, \ir{verhilb} and \ir{smallm}).
Although the nature of our
approach makes it difficult to give a simple description
of all $m$ and $n$ which we can handle, see Figures 1--4
for graphical representations of some of our results.

Our approach combines and extends the methods of
\cite{refnagconj}, \cite{refHHF}, \cite{refR}, and \cite{refRb},
to obtain improved bounds on
$\alpha({\bf m})$ and $\tau({\bf m})$. Sufficiently good
bounds determine these quantities exactly, which
in many situations is sufficient to also determine
$h({\bf m})$ and the minimal free resolution of $I({\bf m})$.
Our bounds are algorithmic (for an example of the algorithm
applied in a specific case, see \ir{algex}). The input to the algorithm
is the sequence ${\bf m}=(m_1,\ldots,m_n)$ of multiplicities
of the $n$ general points. One also must pick a positive integer $d$
and a positive integer $r\le n$, corresponding to
a specialization of the $n$ points in which the first point
is a general smooth point of an irreducible plane curve of 
degree $d$, each additional point is
infinitely near the previous one, and exactly $r$ of them lie on the curve. 
Using properties of linear series on curves, we then obtain
bounds for such sets of specialized points; by semicontinuity
these bounds also apply in the case of general points. For 
fixed ${\bf m}$, note that different choices of $d$ and
$r$ can lead to different bounds, as is shown particularly
clearly by Figures 3 and 4 at the end of this paper.

An analysis of our
algorithm leads to explicit formulas for
these bounds in certain cases.
To write down these formulas,
given $d>0$, $0<r\le n$ and ${\bf m}=(m_1,\ldots,m_n)$,
define integers $u$ and $\rho$ via $M_n=ur+\rho$,
where $u\ge 0$, $0<\rho\le r$
and $M_i=m_1+\cdots+m_i$ for each $1\le i\le n$.
Note that we can equivalently define $u=\lceil M_n/r\rceil-1$ and
$\rho=M_n-ru$. Also,
given an integer $0<r\le n$, we say that
$(m_1,\ldots,m_n)$ is {\it $r$-semiuniform\/} if
$m_r+1\ge m_1\ge m_2\ge \cdots\ge m_n\ge0$. Note that a nonincreasing
sequence $(m_1,\ldots,m_n)$ of nonnegative integers is
$r$-semiuniform if and only if $m_i$ is either $m_r$ or
$m_r+1$ for every $1\le i\le r$; thus, for example,
$(6,6,6,5,5,5,5,4,4,4)$ is $r$-semiuniform for each $r$ up to
7, but not for $r=8, 9$ or 10.
We now have:

\prclm{Theorem}{introthm}{Given integers
$0<d$ and $0<r\le n$, let ${\bf m}=(m_1,\ldots,m_n)$
be $r$-semiuniform, define $M_i$, $u$ and
$\rho$ as above, denote the genus $(d-1)(d-2)/2$ of a plane
curve of degree $d$ by $g$ and let
$s$ be the largest integer
such that we have both $(s+1)(s+2)\le 2\rho$ and $0\le s <d$.
\item{(a)} If $r\le d^2$, then
\itemitem{(i)} $\alpha({\bf m})\ge 1+
\hbox{min}(\lfloor(M_r+g-1)/d\rfloor,s+ud)$
whenever $d(d+1)/2\le r$, while
\itemitem{(ii)} $\tau({\bf m})\le
\hbox{max}(\lceil(\rho+g-1)/d\rceil+ud,ud+d-2)$.
\item{(b)} If $2r\ge n+d^2$, then
\itemitem{(i)} $\alpha({\bf m})\ge s+ud+1$, and
\itemitem{(ii)} $\tau({\bf m})\le
\hbox{max}(\lceil(M_r+g-1)/d\rceil,ud+d-2)$.
\item{(c)} Say for some $m$ we have $m_i=m$ for all $i$,
and that $rd(d+1)/2\le r^2 \le d^2 n$; then $\alpha{(n;m)}\ge
1+\min (\lfloor (mr+g-1)/d \rfloor, s+ud)$.}

In \ir{algs}, we develop the results needed to state and
analyze our algorithm. In \ir{bnds} we prove \ir{introthm}.
In \ir{NC} we apply \ir{introthm} to prove results less easily
stated but stronger than \ir{introcor}(a), from which \ir{introcor}(a)
is an easy consequence. Similarly, \ir{introcor}(b) is
an immediate consequence of more comprehensive
results that we deduce from \ir{introthm} in \ir{hilb},
and \ir{introcor}(c) is an immediate consequence of more comprehensive
results that we deduce in \ir{res} from \ir{introthm} using \cite{refHHF}.

Since it is hard to easily describe
all of the cases that our method handles, we include
for this purpose some graphs in \ir{comps}, together
with some explicit comparisons of our bounds on $\alpha$ and $\tau$ with
previously known bounds.

\irrnSection{Algorithms}{algs}

In this section we derive algorithms giving bounds
on $\alpha({\bf m})$ and $\tau({\bf m})$.
It is of most interest to give lower bounds for $\alpha$ and
upper bounds for $\tau$, since upper bounds for
$\alpha$ and lower bounds for $\tau$ are known
which are conjectured to be sharp. (See \cite{refsurv}
for a discussion.)

Our method involves a specialization of the $n$ points
as in \cite{refnagconj} (which in turn
was originally inspired by that of \cite{refR}), together
with properties of linear series on curves. Recall that a flex
for a linear series $V$ of dimension $a$ on a curve $C$ is
a point $p\in C$ such that $V-(a+1)p$ is not empty.
In \ir{lemcurve} we use the known result (see p. 235, \cite{refMirb})
that the set of flexes of a linear
series is finite; it is the only place that
we need the characteristic to be zero (although, of course,
everything else refers to this Lemma).
In positive characteristics, complete linear series can
indeed have infinitely many flexes \cite{refHo}.

Before deriving our algorithms, we need two lemmas.

\prclm{Lemma}{lemcurve}{Let $C$ be an irreducible
plane curve of degree $d$,
so $g=(d-1)(d-2)/2$ is the genus of $C$, and
let $q$ be a general point of $C$. Take $D=tL_C-vq$, where
$t\ge 0$ and $v\ge 0$ are integers and $L_C$ is the restriction
to $C$ of a line $L$ in \pr2.
\item{(1)}  If $t \ge d-2$ then $h^0(C, D)=0$ for $t \le (v+g-1)/d$
and $h^1(C, D)=0$ for $t\ge (v+g-1)/d$.
\item{(2)}  If $t < d$ then $h^0(C, D)=0$ for $(t+1)(t+2) \le 2v$.}

\Prf We have $L, C\subset$\pr2.
The linear system $|tL_C|$ is the image of $|tL|$
under restriction to $C$. Note that
$h^1(C, tL_C)=0$ if (and only if) $t\ge d-2$,
in which case $h^0(C,tL_C)=td -g +1$,
whereas, for $t < d$,
$h^0(C,tL_C)=h^0(\pr2,tL)=(t+1)(t+2)/2$.

Since $h^1(C, tL_C)=0$ for $t\ge d-2$, as long as $vq$ imposes
independent conditions on $|tL_C|$
(i.e., $h^0(C,tL_C-vq)=h^0(C,tL_C)-v$), then
$h^1(C, tL_C-vq)=0$ too.
Also, if we show that $h^1(C, tL_C-vq)=0$, then it is easy to
see that $h^1(C, tL_C-v'q)=0$ for all $v'\le v$,
and if we show that $h^0(C, tL_C-vq)=0$, then it is easy to
see that $h^0(C, tL_C-v'q)=0$ for all $v'\ge v$.
So it is enough to show $h^0(C, tL_C-vq)=0$
for $v=h^0(C,tL_C)$. Since
$q$ is general in $C$, we can assume it is not a flex
of $|tL_C|$, and therefore the claim follows.
\qed

Consider $n$ distinct points
$p_1, \ldots, p_n$ of \pr2 and let $X$ be the blow up of the points.
More generally, we can allow the possibility
that some of the points are infinitely near by taking
$p_1\in\pr2=X_0$, $p_2\in X_1$, $\dots$, $p_n\in X_{n-1}$,
where $X_i$, for $0<i\le n$, is the blow up of $X_{i-1}$ at $p_i$,
and we take $X=X_n$. Given integers $t$
and $m_1\ge m_2\ge \cdots\ge m_n\ge 0$,
we denote by $F_t$ the divisor $F_t=tL-m_1E_1-\cdots-m_nE_n$ on $X$,
where $E_i$ is the divisorial inverse image of $p_i$ under the
blow up morphisms $X=X_n\to X_{n-1}\to\cdots\to X_{i-1}$,
and $L$ is the pullback to $X$ of a general line in \pr2.
Note that the divisor classes $[L], [E_1],\ldots,[E_n]$
give a basis for the divisor class group $\hbox{Cl}(X)$ of $X$.

Now, given positive integers $r\le n$ and $d$, we choose our
points $p_1, p_2, \ldots, p_n$ such that $p_1$ is a general smooth point
of an irreducible plane curve $C'$ of degree $d$, and then
choose points $p_2, \ldots, p_n$
so that $p_i$ is infinitely near $p_{i-1}$ for $i\le n$
and so that $p_i$ is a point of the proper
transform of $C'$ on $X_{i-1}$ for $i\le r$
(more precisely, so that $[E_i-E_{i+1}]$ is the class of
an effective, reduced and irreducible divisor for $0<i<n$
and so that the class of the proper transform
of $C'$ to $X$ is $[dL-E_1-\cdots-E_r]$).
Let $C$ denote the proper transform of $C'$ in $X$.

Define divisors $D_j$ and $D_j'$ such that $D_0=F_t$,
$D_j'=D_j-(dL-E_1-\cdots-E_r)$, and such that $D_{j+1}$ is obtained from
$D_j'$ by {\it unloading\/} multiplicities
(i.e., if $D_j'=a_0L-a_1E_1-\cdots-a_nE_n$,
then we permute the $a_1, a_2, \ldots, a_n$ so
that $a_1\ge a_2\ge \cdots\ge a_n$
and set to 0 each which is negative).

Part (1) of the following lemma is used in our algorithms.
Parts (2) and (3) are used in the proof of \ir{introthm}(a,b).
\vskip\baselineskip

\prclm{Lemma}{lemunload}{Let $r\le n$ and $d$ be positive integers.
Given $F_t=tL-m_1E_1-\cdots-m_nE_n$ with
$m_1\ge m_2\ge \cdots\ge m_n\ge 0$,
define $C$, $D_j$ and $D_j'$ for $j\ge 0$ as above.
\item{(1)} For every $j>0$, $h^i(X, D_{j-1}')=h^i(X, D_{j})$, $i=0,1$.
\item{(2)} If $r\le d^2$,
then $D_j\cdot C\le D_{j-1}\cdot C$ for every $j>0$.
\item{(3)} If $2r\ge n+d^2$,
then $D_j\cdot C\ge D_{j-1}\cdot C$ as
long as $D_{j-1}\cdot E_r>0$.}

\noindent{\bf Proof}:
Write $D_{j-1}'=a_0L-a_1E_1-\cdots-a_nE_n$.
Because of the definition of $F_t$ and $D_i$,
we have either $a_1\ge a_2\ge \cdots\ge a_n\ge-1$,
or $a_1\ge a_2\ge \cdots\ge a_i= a_{i+1}-1$ and
$a_{i+1}\ge \cdots\ge a_n\ge -1$ for some $i\le r$.
Therefore, the unloading procedure leading from $D_{j-1}'$
to $D_j$ consists in a number of unloading steps, each of which
either transposes $a_k$ and $a_{k+1}$ (whenever $a_k=a_{k+1}-1$)
or sets $a_n$ to 0 (whenever in the course of the transpositions
we find $a_n=-1$). If we denote the proper transform of the exceptional
divisor of blowing up $p_k$ by $\tilde E_k$ (hence
$[\tilde E_k]=[E_k-E_{k+1}]$ for $k<n$ and $\tilde E_n=E_n$),
this is the same as iteratively subtracting
$\tilde E_k$ from $D_{j-1}'$ whenever
$D_{j-1}' \cdot \tilde E_k=-1$, and it is enough to show
that doing so does not affect cohomology in order to prove (1).

Consider the exact sequence
$0\to {\cal O}_X(D_{j-1}'-\tilde E_k) \to {\cal O}_X(D_{j-1}')
\to {\cal O}_X(D_{j-1}')\otimes{\cal O}_{\tilde E_k} \to 0$.
Since $\tilde E_k \cong \pr1$ and $D_{j-1}' \cdot \tilde E_k=-1$,
we have $h^i({\cal O}_X(D_{j-1}')\otimes{\cal O}_{\tilde E_k})=0$.
Thus, taking cohomology of the sequence,
we see that the cohomology of ${\cal O}_X(D_{j-1}'-\tilde E_k)$
and ${\cal O}_X(D_{j-1}')$ coincide, as claimed.

To prove (2), let $D_{j-1}'=a_0L-a_1E_1-\cdots -a_nE_n$
as above and observe that
$D_{j-1}'\cdot C= D_{j-1}\cdot C-d^2+r$.
Passing from $D_{j-1}'$ to $D_j$ by unloading
might increase some of the coefficients $a_i$ with $0<i\le r$ but
cannot decrease any of these coefficients and hence
cannot increase the intersection with $C$; i.e.,
$D_j\cdot C \le D_{j-1}'\cdot C =D_{j-1}\cdot C -d^2 +r\le D_{j-1}\cdot C$.

To prove (3), let $D_{j-1}'=a_0L-a_1E_1-\cdots -a_nE_n$
as above. If $D_{j-1}\cdot E_r>0$, then $a_i\ge 0$ for all $i>0$
and $a_i+1\ge a_j$ for all $0<i\le r<j\le n$.
Thus passing from $D_{j-1}'$ to $D_j$ may involve swapping
some of the coefficients $a_i$, $0<i\le r$, with some of the
coefficients $a_j$, $j>r$, but there are only $n-r$
coefficients $a_j$ with $j>r$, each of which is at most 1
bigger than the least coefficient $a_i$ with $0<i\le r$, so
passing from $D_{j-1}'$ to $D_j$ can decrease the intersection with
$C$ by at most $n-r$; i.e., $D_j\cdot C \ge D_{j-1}'\cdot C -(n-r)=
D_{j-1}\cdot C -d^2 +r-(n-r)\ge D_{j-1}\cdot C$.
\qed

For \ir{introthm}(c) we need a slightly more general
version of the preceding lemma, restricted however to the case
of uniform multiplicities:
\vskip\baselineskip

\prclm{Lemma}{intlemma}{Let $r\le n$ and $d$ be positive integers.
Let $F_t=tL-mE_1-\cdots-mE_n$, with respect to which we take
$C$, $D_i$ and $D_i'$ for all $i\ge0$ as in \ir{lemunload}, and let
$\omega'$ be the least $i$ such that $D_i\cdot E_i=0$ for all $i>0$.
Then for all $0\le i\le \omega'-1$ we have
$$i\left({{r^2}\over{n}}-d^2\right)-\left(r-{{r^2}\over{n}}\right)
\le D_i \cdot C - D_0\cdot C \le
i\left({{r^2}\over{n}}-d^2\right).$$}

\noindent {\bf Proof}:
Let $A_0=0$, and for $0< k\le n$ let $A_k=-E_1-\cdots-E_k$
and define $t_i$ to be $D_i\cdot L$.
For $0\le i<\omega'$, it is not hard to check that
$D_i=(t_0-id)L-(m-i+q)E_1-\cdots-(m-i+q)E_n+A_\rho$,
where $i(n-r)=qn+\rho$ with $0\le\rho<n$,
and therefore
$$D_i \cdot C - D_0 \cdot C= i(r-d^2) - rq +A_\rho \cdot C.$$
On the other hand, $A_\rho \cdot C=-\min(\rho,r)$,
and it is easy to see that
$$i(n-r){r\over n} \le rq+\min(\rho,r) \le
i(n-r){r\over n}+r-{{r^2}\over n},$$
from which the claim follows.
\qed

We now derive our algorithms. The reader may
find the algorithm easier to follow by looking at the example
we give below.
Let $X$ be obtained by blowing up $n$ general points of
\pr2, let $F_t=tL-m_1E_1-\cdots-m_nE_n$, where
we assume that $m_1\ge\cdots\ge m_n\ge 0$, and choose
any integers $d$ and $r$ such that
$d>0$ and $0<r\le n$. Next, specialize the points
as in \ir{lemunload}. We then have the specialized
surface $X'$. It is convenient to denote the basis of
the divisor class group of $X'$ corresponding to the specialized
points also by $L$, $E_1,\ldots,E_n$, since it will always
be clear whether we are working on $X'$ or on $X$.
By semicontinuity, we know $h^i(X', F_t)\ge h^i(X, F_t)$,
so for any $t$ with $h^i(X', F_t)=0$ we also have $h^i(X, F_t)=0$.

As in \ir{lemunload}, we have on $X'$ the curve $C$
with class $[dL-E_1-\cdots-E_r]$ and genus $g=(d-1)(d-2)/2$
and we have the sequence of divisors $D_j$, where
$D_0=F_t$, $D_j'=D_j-(dL-E_1-\cdots-E_r)$, and
$D_{j+1}$ is obtained from $D_j'$ as above.

Whenever $h^0(X', F_t)=0$, we get the
bound $t+1\le \alpha(m_1,\ldots,m_n)$,
and whenever $h^1(X', F_t)=0$, we get the bound
$t\ge \tau(m_1,\ldots,m_n)$. But, using the exact sequences
$$0\to {\cal O}_{X'}(D_j')\to {\cal O}_{X'}(D_j)
\to {\cal O}_{X'}(D_j)\otimes{\cal O}_C \to 0\eqno\hbox{$(*)$}$$
and \ir{lemunload}(1), we see that $h^0(X', F_t)=0$ if for some $i=I$
we have $h^0(X', D_I)=0$ and for all $0\le i<I$ we have
$h^0(C, {\cal O}_{X'}(D_i)\otimes{\cal O}_C)=0$.
Similarly, $h^1(X', F_t)=0$ if for some $j=J$
we have $h^1(X', D_J)=0$ and for all $0\le j<J$ we have
$h^1(C, {\cal O}_{X'}(D_j)\otimes{\cal O}_C)=0$.
But, ${\cal O}_{X'}(D_j)\otimes{\cal O}_C=
{\cal O}_{X'}((t-jd)L)\otimes{\cal O}_C(-v_jq)$,
where $v_j=(t-jd)d-D_j\cdot C$ and $q$ is the point of $C$
infinitely near to $p_1$. Thus we can apply the criteria
of \ir{lemcurve} to control
$h^0(C, {\cal O}_{X'}(D_i)\otimes{\cal O}_C)$ for $0\le i<I$
and $h^1(C, {\cal O}_{X'}(D_j)\otimes{\cal O}_C)$ for $0\le j<J$.

In particular, for a given value of $t$ we have
$h^0(X', D_I)=0$ if we take $I$ to be the
least $i$ such that $D_i\cdot (L-E_1)<0$. Then,
if $t$ is not too large, for all $0\le i<I$ we also have
either $D_i\cdot C\le g-1$ and $t-id\ge d-2$, or $t-id<0$,
or $t-id<d$ and $(t-id+1)(t-id+2)\le 2v_i$,
which guarantees by \ir{lemcurve} that $h^0(C, D_i)=0$
for all $0\le i<I$. The largest such $t$ then gives
a lower bound for $-1+\alpha(m_1,\ldots,m_n)$;
i.e., $t+1\le \alpha(m_1,\ldots,m_n)$.

Similarly, let $J$ be the least $j$ such that $D_j$ is a multiple of $L$
(i.e., such that if $D_j=a_0L-a_1E_1-\cdots-a_nE_n$,
then $a_1=\cdots=a_n=0$, in which case we have $h^1(X', D_j)=0$).
If we have chosen $t$ sufficiently large, then $D_j\cdot L\ge d-2$ and
$D_j\cdot C\ge g-1$, for all $0\le j<J$;
the least such $t$ then gives
an upper bound for $\tau(m_1,\ldots,m_n)$.

\rem{Example}{algex} For this example, suppose we wish
to get information on $n=18$ general points of multiplicity $m=2$.
We first must pick $d$ and $r$. Any integers $d>0$ and $0<r\le n$ will do, 
but values such that $r/d$ is close
to $\sqrt{n}$ tend to be best; here we will choose $d=4$ and $r=17$.
Now we pick specialized points $p_1,\ldots,p_{18}$ such that
$p_1$ is a general point of a reduced irreducible plane curve
$C'$ of degree $d=4$, and $p_2$ is the
point on $C'$ infinitely near to $p_1$, and so on for the first
$r=17$ points. Thus $p_{17}$ is the
point on $C'$ infinitely near to $p_{16}$, but each point $p_i$
for $i>17$ is taken to be a general point infinitely near to $p_{i-1}$
and so in particular $p_i$ is not on $C'$. For this example, $n=18$ so
there is only one such point, $p_{18}$. 

Now let $C$ be the proper transform of $C'$ to the
surface $X'$ obtained by blowing up each point in turn.
Thus the class of $C$ is $dL-E_1-\cdots-E_{17}$,
and $F_t=tL-2E_1-\cdots-2E_{18}$. 
We want to find the largest $t$ such that we can show
that $h^0(X',F_t)=0$. If we can show that $h^0(X',F_t)=0$,
then we try $t+1$, and we continue increasing $t$ until we are
unable to show that $h^0(X',F_t)=0$. For example,
say $t=8$. Then, using the notation above, $D_0=F_8$ and $D'_0=F_8-C$, and
we use \ir{lemcurve} to determine whether or not
$h^0(C,{\cal O}_{X'}(D_0)\otimes{\cal O}_C)=0$. In this case
$h^0(C,{\cal O}_{X'}(D_0)\otimes{\cal O}_C)=0$ is obvious since
$D_0\cdot C=8\cdot4-17\cdot2=-2$. By the sequence $(*)$ of sheaves 
above, it follows that $h^0(X',D_0)=h^0(X',D'_0)$. 
Next, we {\it unload\/} $D'_0=4L-E_1-\cdots-E_{17}-2E_{18}$
to obtain $D_1$. This just amounts to reordering the coefficients of the $E_i$
to be nondecreasing; thus $D_1=4L-2E_1-E_2-\cdots-E_{17}-E_{18}$.
(Since $p_{18}$ is infinitely near $p_{17}$, $N_{17}=E_{17}-E_{18}$ is
the class of a reduced irreducible curve. But $D'_0\cdot N_{17}<0$,
so $h^0(X',D'_0)=h^0(X',D'_0-N_{17})$, and
$(D'_0-N_{17})\cdot N_{16}<0$, so $h^0(X',D'_0-N_{17})=h^0(X',D'_0-N_{17}-N_{16})$.
Note that 
$D'_0-N_{17}=4L-E_1-\cdots-E_{16}-2E_{17}-E_{18}$, and 
$D'_0-N_{17}-N_{16}=4L-E_1-\cdots-E_{15}-2E_{16}-E_{17}-E_{18}$;
unloading consists of continuing to subtract classes $N_i=E_i-E_{i+1}$
until the coefficients of the $E_i$ are nondecreasing. Fully
unloading $D'_0$ results in the class $D_1$, for which we see
we have $h^0(X',D'_0)=h^0(X',D_1)$.) Now we repeat the process above,
applied to $D_1$ in place of $D_0$.
Again we have $h^0(C,{\cal O}_{X'}(D_1)\otimes{\cal O}_C)=0$,
so $h^0(X',D_1)=h^0(X',D'_1)=h^0(X',D_2)$, but now
$D_2=-E_1-E_2$, so $h^0(X',D_2)=0$, hence $h^0(X',F_8)=0$,
so (by semicontinuity) we know $\alpha(18;2)>8$. At this point we repeat
the whole process with $t=9$, in hopes of improving our bound.
Thus $D_0$ is now $F_9$. However, if for some $i$ it ever happens
that either $h^0(C,{\cal O}_{X'}(D_i)\otimes{\cal O}_C)\ne0$
or $h^0(X',D_i)\ne 0$, then we will be unable to conclude
that $h^0(X',F_9)=0$. In fact, $D_2$ in this case is
$L-E_1-E_2$, so $h^0(X',D_2)\ne 0$, so we are indeed
unable to conclude that $h^0(X',F_9)=0$.
Thus our algorithm cannot
improve on the bound $\alpha(18;2)>8$ so our bound in the end is
$9\le \alpha(18;2)$. Our bound on $\tau$ works similarly,
except we start with $t$ large enough so that $h^1$
vanishes both for $D_j$ for some $j$ and for
the restrictions of $D_i$ to $C$ for all $0\le i <j$. The least such $t$
gives us the bound $t\ge \tau(n;m)$. In the case of our example,
we find $9\ge \tau(18;2)$. 
Thus the 54 conditions imposed by 18 general double points 
on the 55 dimensional space of forms of degree $9$ are independent,
so in fact $h^0(X,F_9)=1$ and $\alpha(18;2)=9$. Moreover, 
since $\alpha(18;2)=9=\tau(18;2)$ and $h^0(X,F_9)=1$, 
we know a minimal set of homogeneous generators
of $I(18;2)$ contains a single generator in degree 9
and $h^0(X,F_{10}) - 3h^0(X,F_9) = 12-3=9$ generators
in degree 10. Generators are never needed
in degrees greater than $\tau+1$, so it follows that
the resolution of $I(18;2)$ is
$0\to R[-11]^9 \to  R[-10]^9\oplus R[-9]   \to I(18;2)\to0$.

\irrnSection{Proof of \ir{introthm}}{bnds}

In this section we prove the explicit bounds on $\alpha$ and
$\tau$ claimed in \ir{introthm}. The analysis of the
algorithm of the previous section, which leads to the proof,
goes differently depending on whether $r$ is relatively small (a) or
big (b) compared to $d$; we also include a separate analysis for the
case of uniform multiplicities (c).

We begin with (a)(i),
so assume $d(d+1)/2\le r\le d^2$ and ${\bf m}=(m_1,\ldots,m_n)$. To show
$1+t\le \alpha({\bf m})$ for
$t=\hbox{min}(\lfloor(M_r+g-1)/d\rfloor,s+ud)$,
it is by our algorithm in \ir{algs} enough (taking $I=u+1$) to show
that $D_I\cdot (L-E_1)<0$, and for all $0\le i<I$ that either
$D_i\cdot C\le g-1$ and $t-id\ge d-2$,
or $t-id<0$, or $0\le t-id<d$ and $(t-id+1)(t-id+2)\le 2v_i$,
where $v_j=(t-jd)d-D_j\cdot C$, as before.

First, by \ir{lemunload}(2), $D_0\cdot C\ge D_1\cdot C\ge\cdots$,
and by hypothesis $t\le\lfloor(M_r+g-1)/d\rfloor$,
so $g-1\ge td-M_r=D_0\cdot C\ge D_1\cdot C\ge\cdots$, as required.
Also by hypothesis, we have $t\le s+ud$.
It follows that $D_I\cdot L=t-Id$
and hence that $D_I\cdot (L-E_1)<0$, as required.
Thus it is now enough to check that $(t-id+1)(t-id+2)\le 2v_i$
for the largest $i$ (call it $i''$) such that $t-id\ge 0$.
If $i''=I-1$ we have $t-i''d=t-ud\le s$ by hypothesis
and hence $(t-i''d+1)(t-i''d+2)\le (s+1)(s+2)\le 2\rho=2v_{I-1}$
by definition of $s$. If $i''<I-1$, we at least have
$t-i''d\le d-1$, so $(t-i''d+1)(t-i''d+2)\le d(d+1)$.
But $v_{I-2}\ge r$ by $r$-semiuniformity,
so $M_r=v_0\ge v_1\ge\cdots\ge v_{I-2}\ge r\ge v_{I-1}=\rho>v_I=0$,
hence $2v_{i''}\ge 2r\ge d(d+1)\ge (t-i''d+1)(t-i''d+2)$,
as we wanted.

We now prove (a)(ii). As always we have $r\le n$;
in addition we assume $r\le d^2$. It follows from
semiuniformity that $D_J$
is a multiple of $L$ for $J=u+1$,
so it suffices to show for
$t=\hbox{max}(\lceil(\rho+g-1)/d\rceil+ud,ud+d-2)$
that $D_j\cdot L\ge d-2$ and $D_j\cdot C\ge g-1$,
for all $0\le j<J$; it then
follows by our algorithm that $t\ge\tau({\bf m})$.
But $t\ge ud+d-2$ ensures that $D_j\cdot L\ge d-2$
for all $0\le j<J$, and, since
$D_0\cdot C\ge D_1\cdot C\ge\cdots\ge D_u\cdot C=(t-ud)d-\rho$
by \ir{lemunload}, $t\ge\lceil(\rho+g-1)/d\rceil+ud$
ensures that $D_j\cdot C\ge g-1$, for all $0\le j<J$.

Next, consider (b). Now we assume that $2r\ge n+d^2$
and that ${\bf m}$ is $r$-semiuniform.
By semiuniformity we have $D_i\cdot E_r>0$
for $i<u$. Now by \ir{lemunload}(3) we have
$D_0\cdot C\le D_1\cdot C\le\cdots\le D_u\cdot C$.

Starting with (b)(i), let $I=u+1$ and $t=s+ud$.
It suffices to check that $D_I\cdot (L-E_1)<0$, and for all $0\le i<I$
that either $D_i\cdot C\le g-1$ and $t-id\ge d-2$,
or $t-id<d$ and $(t-id+1)(t-id+2)\le 2v_i$.
But $D_I\cdot (L-E_1)=t-(u+1)d=s-d<0$, as required, so
now consider $D_u$. Here we have $D_u=sL-(E_1+\cdots+E_\rho)$.
But by hypothesis $s<d$ and $(s+1)(s+2)\le 2\rho=2v_u$,
as required. Finally, consider
$D_i$ for $0\le i<u$. Then $D_i\cdot L=t-id=s+(u-i)d\ge d$,
and $sd-\rho=D_u\cdot C\ge D_i\cdot C$,
so if we prove that $s\le (\rho+g-1)/d$ we will have
$D_i\cdot C\le g-1$ as we want. But it is easy to see that
$sd-g+1=(s+1)(s+2)/2-(s-d+1)(s-d+2)/2 \le (s+1)(s+2)/2$,
and therefore the hypothesis $2 \rho \ge (s+1)(s+2)$
implies $\rho \ge sd-g+1$, so we are done.

Next, we prove (b)(ii). Let $J=u+1$;
then $D_J=(t-Jd)L$ is a multiple of $L$, as by our
algorithm we would want. Now assume in addition that
$t\ge \hbox{max}(\lceil(M_r+g-1)/d\rceil,ud+d-2)$.
We want to verify that $D_j\cdot L\ge d-2$ and
$D_j\cdot C\ge g-1$, for all $0\le j<J$.
First consider $j=0$; since $t\ge(M_r+g-1)/d$ and
$t\ge ud+d-2\ge d-2$, we have $dt-M_r=D_0\cdot C\ge g-1$
and $t=D_0\cdot L\ge d-2$. As for $0<j<J$, we have
$t-jd \ge (ud+d-2)-(ud)=d-2$
and $D_j\cdot C\ge D_0\cdot C\ge g-1$,
which ends the proof of \ir{introthm}(b).

Finally, we prove (c). In the notation of \ir{intlemma}
and its proof, it is easy to check that
$\omega'=\lceil mn/r\rceil=u+1$,
so if $t \le s+ud$, it follows that $t_{\omega'}\le s-d<0$,
and thus $\omega'\ge\omega$, where $\omega$
is the least $i$ such that $t_i<0$. Since $r^2/n-d^2 \le 0$,
it follows from \ir{intlemma} that
$D_i\cdot C\le D_0 \cdot C$ for all $0\le i\le \omega-2$.
If $t\le \lfloor (mr+g-1)/d \rfloor$, then
$D_i\cdot C\le D_0 \cdot C = td-mr \le g-1$.
To conclude that $\alpha(n;m)\ge t+1$, it is now enough to
check that $(t-i'd+1)(t-i'd+2)\le 2v_{i'}$ for $i'=\omega-1$.
If $i'=u$ (i.e., $\omega'=\omega$)
we have $t-i'd=t-ud\le s$ by hypothesis and hence
$(t-i'd+1)(t-i'd+2)\le (s+1)(s+2)\le 2\rho=2v_{i'}$ by
definition of $s$. If $i'<u$ (so  $\omega'>\omega$),
by definition of $i'$ we at least have $t-i'd\le d-1$, so
$(t-i'd+1)(t-i'd+2)\le d(d+1)$.
But $\omega'>\omega$ implies $v_{i'}\ge r$, and
by hypothesis $rd(d+1)/2\le r^2$ (so $d(d+1)\le 2r$);
therefore $2v_{i'}\ge 2r\ge d(d+1)\ge (t-i'd+1)(t-i'd+2)$
as we wanted.
\qed

\rem{Remark}{remexs} Even for uniform multiplicities,
sometimes the best bound determined by \ir{introthm}
comes from parts (a) or (c), and sometimes it comes from
part (b), depending on $n$ and $m$. Sometimes, of course, one can do
better applying our algorithm for values of $r$ and $d$ for
which \ir{introthm} does not apply.

For example, let $\alpha_c(n;m)$ denote the conjectural
value of $\alpha(n;m)$ and let $\tau_c(n;m)$ denote the
conjectural value of $\tau(n;m)$ (i.e., the values of
each assuming \ir{introconj}(b) holds).
Then $\alpha_c(33;29)=168$; the best bound
given by \ir{introthm} is
$\alpha(33;29)\ge 165$, obtained in part (b)
using $r=29$ and $d=5$, or in part (c) using
$r=17$ and $d=3$. Applying our algorithm
with $r=23$ and $d=4$, however, gives
$\alpha(33;29)\ge168$ (and hence
$\alpha(33;29)=\alpha_c(33;29)$).
On the other hand,  $\alpha_c(38;16)=101$
and indeed we obtain $\alpha(38;16)\ge 101$
via \ir{introthm}(b) using $r=37$ and $d=6$,
while the best bound obtainable via
\ir{introthm}(c) is $\alpha(38;16)\ge 98$,
gotten using $r=36$ and $d=6$. In contrast, we obtain
$\alpha(119;13)\ge 146=\alpha_c(119;13)$
via \ir{introthm}(c) using $r=109$ and $d=10$,
while the best bound obtainable via
\ir{introthm}(b) is $\alpha(119;13)\ge 144$,
gotten using $r=100$ and $d=9$.

Similarly, $\tau_c(33;29)=168$; applying our algorithm
with $r=23$ and $d=4$, gives $\tau(33;29)\le169$.
The best bound given by \ir{introthm}(b) is
$\tau(33;29)\le 170$, obtained using $r=29$ and $d=5$, while
the best bound given by \ir{introthm}(a) is
$\tau(33;29)\le 175$, obtained
using $r=33$ and $d=6$.
On the other hand,  $\tau_c(38;16)=101$
and indeed we obtain $\tau(38;16)\le 101$
via \ir{introthm}(b) using $r=37$ and $d=6$,
while the best bound obtainable via
\ir{introthm}(c) is $\tau(38;16)\le 103$,
gotten using $r=36$ and $d=6$. In contrast,
$\tau_c(119;13)=146$, and the best bound obtainable using our algorithm
is $\tau(119;13)\le 147$, obtained using $r=119$ and $d=11$
(and hence \ir{introthm}(a) applies),
while the best bound obtainable via
\ir{introthm}(b) is $\tau(119;13)\le 148$,
gotten using $r=111$ and $d=10$.

\irrnSection{Nagata's Conjecture}{NC}
In this section we prove \ir{introcor}(a) as
an immediate easy-to-state consequence of our following more
involved result. Because \ir{introconj}(a) is known when
$n$ is a square, we need not consider that case.

\prclm{Corollary}{vernag}{Given an integer $n>9$, let
$d = \lfloor \sqrt{n} \rfloor$ (hence $d \ge 3$) and
$\Delta=n-d^2$. Then $\alpha(n;m)\ge m\sqrt{n}$
holds whenever:
\item{(a)} $\Delta$ is odd and $m \le
\hbox{max}(d(d-3),d(d-2)/\Delta)$, or
\item{(b)} $\Delta>0$ is even, and $m \le
\hbox{max}(d(d-3)/2,2d^2/\Delta)$.}

\noindent{\bf Proof}: Both claims follow from \ir{introthm},
with $d = \lfloor \sqrt{n} \rfloor$ and appropriate 
choices of $r$. Consider part (a).
We first prove $\alpha(n;m)\ge m\sqrt{n}$
if $\Delta$ is odd and $m \le d(d-2)/\Delta$.
Apply \ir{introthm}(a)
with $r=d^2$, $u=m$ and $\rho=m \Delta$;
it has to be checked that $md +s +1 \ge m \sqrt{n}$ and
$\lfloor (mr+g-1)/d \rfloor +1 \ge m \sqrt{n}$.
The first inequality is equivalent
to $(s+1)^2 +2(s+1)md \ge m^2 \Delta$.
If $s=0$ then $m \Delta <3$ and the inequality follows
from $d \ge 3$,
whereas if $s=d-1$ then $d^2 > d(d-2) \ge m \Delta$ says
$d^2+2md^2 > m^2 \Delta$.
In all intermediate cases one has $2d(s+1) \ge (s+2)(s+3)$
and $(s+2)(s+3)> 2 \rho= 2m \Delta$
which also imply
$(s+1)^2 +2(s+1)md \ge m^2 \Delta$ easily.
To prove the second inequality it is enough to see
that $md+d/2 -1 \ge m \sqrt{n}$, which is equivalent
to $(d/2-1)^2+md(d-2) \ge m^2 \Delta$, and
this follows from $m \le d(d-2)/\Delta$.

We now prove $\alpha(n;m)\ge m\sqrt{n}$
if $\Delta$ is odd and $m \le d(d-3)$.
We can write $\Delta=2t+1$ for some nonnegative integer $t$,
hence $n=d^2+2t+1$. Apply \ir{introthm}(c) with $r=d^2+t$
(i.e., $r=\lfloor d\sqrt{n}\rfloor$, hence $rd(d+1)/2\le r^2 \le d^2 n$).
Note $\lfloor (mr+g-1)/d\rfloor +1 > (mr+g-1)/d = m(d^2+t)/d+(d-3)/2$,
but $m(d^2+t)/d+(d-3)/2 > m\sqrt{d^2+2t+1} = m\sqrt{n}$
for $m\le d(d-3)$. On the other hand,
since $r^2\le d^2n$, we see that $m\sqrt{n}\le mnd/r$, so
it suffices to show that $mnd/r\le s+ud+1$. If
$s=d-1$, then $s+ud+1=(u+1)d=\lceil mn/r\rceil d\ge mnd/r$
as required, so assume $(s+1)(s+2)\le 2\rho<(s+2)(s+3)$
and $s+2\le d$. Then $r(s+ud+1)=r(s+1)+mnd-d\rho$,
so we need only check that $r(s+1)+mnd-d\rho\ge mnd$,
or even that $r(s+1)\ge d(s+2)(s+3)/2$ (which is
clear if $s=0$ since $d\ge 3$) or that $r\ge d^2(s+3)/(2(s+1))$
(which is also clear since now we may assume $s\ge1$).

Now consider (b). First assume $m \le 2d^2/\Delta$.
Let $\Delta=2t$ and apply \ir{introthm}(b)
with $r=d^2+t$ (so again $r=\lfloor d\sqrt{n}\rfloor$), 
$u=m$ and $\rho=mt$;
it has to be checked that $md +s +1 \ge m \sqrt{n}$
or, equivalently, that $(s+1)^2 +2(s+1)md \ge 2 m^2t$.
If $s=0$ then $mt<3$ and the inequality follows
from $d \ge 3$, whereas if $s=d-1$ then $m\le 2d^2/\Delta$ implies
$d^2+2md^2 > 2 m^2t$. In all intermediate cases one
has $2d(s+1) \ge (s+2)(s+3)$ and $(s+2)(s+3)>2\rho= 2mt$
which also imply $(s+1)^2 +2(s+1)md \ge 2 m^2t$ easily.

Finally, assume $m \le d(d-3)/2$. Again $\Delta=2t$
so $n=d^2+2t$; take $r=d^2+t-1$ and apply \ir{introthm}(c)
in the same manner as previously.
\qed

\noindent{\bf Proof} of \ir{introcor}(a):
It follows from \ir{vernag} that $\alpha(n;m)\ge m\sqrt{n}$
holds for all $n\ge 10$ if $m\le d(d-3)/2$, where
$d=\lfloor\sqrt{n}\rfloor$ (given that
$\alpha(n;m)\ge m\sqrt{n}$ is known
and indeed easy to prove when $n$ is a square).
But $\sqrt{n}\ge d$ and $d^2+2\sqrt{n}\ge n$, so obviously
$d(d-3)/2\ge (n-5\sqrt{n})/2$.
\qed

\irrnSection{Hilbert Functions}{hilb}
We now consider the problem of determining
the Hilbert function of an ideal of the form $I(n;m)$.
Typically \ir{introthm}(b) gives a lower bound $\lambda_\alpha(n;m)$
on $\alpha(n;m)$ which is smaller than the upper bound $\Lambda_\tau(n;m)$
it gives for $\tau(n;m)$, but there are in fact many cases for which
$\lambda_\alpha(n;m)\ge \Lambda_\tau(n;m)$. In any such case,
it follows that $\alpha(n;m)\ge \tau(n;m)$, which clearly
implies \ir{introconj}(b) for the given $n$ and $m$.
This is precisely the method of proof of the next result.

\prclm{Corollary}{verhilb}
{Let $d\ge 3$, $\varepsilon>0$ and $i>0$ be integers,
and consider $n=d^2+2\varepsilon$. Then \ir{introconj}(b)
holds for the given $n$ and $m$
if $m$ falls into one of the following ranges:
\item{(a)} $(d-1)(d-2)/(2\varepsilon)\le m < (d+2)(d+1)/(2\varepsilon)$;
\item{(b)} $(i(d^2+\varepsilon)+(d-1)(d-2)/2)/\varepsilon\le m
\le (id^2+d(d-1)/2)/\varepsilon$;
\item{(c)} $(i(d^2+\varepsilon)+d(d-1)/2)/\varepsilon\le m
\le (id^2+(d+1)d/2)/\varepsilon$; and
\item{(d)} $(i(d^2+\varepsilon)+(d+1)d/2)/\varepsilon\le m
\le (id^2+d(d+3)/2)/\varepsilon$.}

\noindent{\bf Proof}: Case (a) is most easily treated
by considering three subcases:
(a1) $(d-1)(d-2)/(2\varepsilon)\le m < d(d-1)/(2\varepsilon)$;
(a2) $d(d-1)/(2\varepsilon)\le m< (d+1)d/(2\varepsilon)$;
and (a3) $(d+1)d/(2\varepsilon)\le m< (d+2)(d+1)/(2\varepsilon)$.

For the proof, apply \ir{introthm} with
$r=d^2+\varepsilon$, $u=m+i$ and $\rho=m\varepsilon-ir$ (with $i=0$ for
part (a)). The reader will find in cases (a1) and (b)
that $s=d-3$, while $s=d-2$ in cases (a2) and (c),
and $s=d-1$ in cases (a3) and (d).
It follows from \ir{introthm} that
$\lambda_\alpha(n;m)\ge \Lambda_\tau(n;m)$, and hence,
as discussed above, \ir{introconj}(b) holds for the given $n$ and $m$.
\qed

For each $n$, there is a finite set of values of
$d$, $\varepsilon$ and $i$ to which \ir{verhilb} can be
profitably applied. For example, in parts (b), (c) and (d) of
\ir{verhilb} we may assume $i \le (d-1)/\varepsilon$,
$i \le d/\varepsilon$ and $i \le d/\varepsilon$
respectively, as otherwise the
corresponding range of multiplicities is empty.
Thus, \ir{verhilb} determines a finite set $V_n$ of values
of $m$ for which \ir{introconj}(b) must hold.
Between the least and largest $m$ in $V_n$ there
can also be many integers $m$
which are not in $V_n$.
For example, of the 4200 pairs
$(n,m)$ with $10\le n=d^2+2\varepsilon\le 100$
and $1\le m\le 100$, there are 723
with $m\in V_n$. Of these, 308 have $m\le 12$ (and thus for these
\ir{introconj}(b) was verified  by \cite{refCMb}); the other
415 were not known before our results.

It is also noteworthy that in many cases we verify
\ir{introconj}(b) for quite large values of $m$.
In particular, if $n=d^2+2$,
it follows from \ir{verhilb} that $m\in V_n$ for
$m=d(d^2+1)+d(d+1)/2$. Thus we have $243\in V_{38}$, for example,
and $783\in V_{83}$. Apart from special cases when $n$ is a square
(in particular, when $n$ is a power of 4; see \cite{refEb}), no
other method we know can handle such large multiplicities. On the
other hand, as indicated by \ir{introcor}(b), if $n$ is any square larger
than 10, our method also handles arbitrarily large values of $m$,
as we now prove. For the purpose of
stating the result, given any positive integer $i$, let
$l_i$ be the largest integer $j$ such that $j(j+1)\le i$.

\prclm{Corollary}{sqrhilb}{Consider $10\le n=\sigma^2$
general points of \pr2.
Let $k$ be any nonnegative integer, and let $m=x+k(\sigma-1)$,
where $x$ is an integer satisfying $\sigma/2-l_{\sigma}\le x\le \sigma/2$
if $\sigma$ is even, or $(\sigma+1)/2-l_{2\sigma}\le x\le (\sigma+1)/2$
if $\sigma$ is odd. Then \ir{introconj}(b) holds for $I(n;m)$.}

\noindent {\bf Proof}:
We apply \ir{introthm}(c) with $d=\sigma-1$, $r=d\sigma$,
$u=\lceil mn/r\rceil-1=m+k$ and
$\rho=mn-ur=x\sigma$. We claim
that $t_0 \le \min (\lfloor (mr+g-1)/d \rfloor, s+ud)$,
where $t_0=m\sigma+\sigma/2-2$ if $\sigma$ is even
and $t_0=m\sigma+(\sigma-1)/2-2$ if $\sigma$ is odd.
But $t_0\le(mr+g-1)/d$ because
$(mr+g-1)/d=(md\sigma+d(d-3)/2)/d=m\sigma+(\sigma-4)/2$.
To see $t_0 \le s+ud$, note that $t_0 \le s+ud$
simplifies to $x+\sigma/2-2 \le s$ if $\sigma$ is even and to
$x+(\sigma-1)/2-2 \le s$ if $\sigma$ is odd.
Therefore (by definition of $s$) we have to
check that $x+\sigma/2-1 \le d$ and
$(x+\sigma/2-1)(x+\sigma/2)\le 2\sigma x$
if $\sigma$ is even, and that $x+(\sigma-1)/2-1 \le d$ and
$(x+(\sigma-1)/2-1)(x+(\sigma-1)/2)\le 2\sigma x$ if $\sigma$ is odd.
The first inequality follows from $x\le \sigma/2$ and
$x \le (\sigma+1)/2$ respectively. For the second,
substituting $\sigma/2-j$ for $x$ if $\sigma$ is even and
$(\sigma+1)/2-j$ for $x$ if $\sigma$ is odd,
$(x+\sigma/2-1)(x+\sigma/2)\le 2\sigma x$ and
$(x+(\sigma-1)/2-1)(x+(\sigma-1)/2)\le 2\sigma x$ resp.
become $j(j+1)\le \sigma$ if $\sigma$ is even and $j(j+1)\le 2\sigma$ if
$\sigma$ is odd. Thus $(x+\sigma/2-1)(x+\sigma/2)\le 2\sigma x$ and
$(x+(\sigma-1)/2-1)(x+(\sigma-1)/2)\le 2\sigma x$ resp.
hold if $x$ is an integer satisfying $\sigma/2-l_{\sigma}\le x\le \sigma/2$
if $\sigma$ is even, and $(\sigma+1)/2-l_{2\sigma}\le x\le (\sigma+1)/2$
if $\sigma$ is odd.

This shows by \ir{introthm}(c) that $\alpha(n;m)\ge m\sigma+\sigma/2-1$
if $\sigma$ is even and
$\alpha(n;m)\ge m\sigma+(\sigma-1)/2-1$ if $\sigma$ is odd. But
since $n$ points of multiplicity $m$ impose at most
$n\binom{m+1}{2}$ conditions on forms of degree $t$, it
follows that $h(n;m)(t)\ge \binom{t+2}{2}-n\binom{m+1}{2}$,
and it is easy to check  that $\binom{t+2}{2}-n\binom{m+1}{2}>0$
whenever $t\ge m\sigma+\sigma/2-1$ if $\sigma$ is even
and $t\ge m\sigma+(\sigma-1)/2-1$ if $\sigma$ is odd.
Thus in fact we have $\alpha(n;m)=m\sigma+\sigma/2-1$
if $\sigma$ is even and
$\alpha(n;m)=m\sigma+(\sigma-1)/2-1$ if $\sigma$
is odd, whenever $m$ is of the
form $m=x+k(\sigma-1)$, with $x$ as given in the statement
of \ir{sqrhilb}.

Of course, $h(n;m)(t)=0$ for all $t<\alpha(n;m)$, and
by \cite{refHHF}, we know that $h(n;m)(t)=
\binom{t+2}{2}-n\binom{m+1}{2}$
for all $t\ge\alpha(n;m)$ (apply Lemma 5.3 of \cite{refHHF},
keeping in mind our explicit expression for $\alpha(n;m)$).
\qed

Note that \ir{introcor}(b) is an immediate consequence of
the preceding result.

\irrnSection{Resolutions}{res}
We now show how our results verify many cases of
\ir{introconj}(c) also, including cases with $m$ arbitrarily
large. Indeed, in addition to the case that $n$ is an even square
treated in \ir{sqrres} below, we have by the following proposition
the resolution for 121 of the 723 pairs $(n,m)$ with $m\in V_n$
mentioned above, and of these
121, 91 have $m>2$ and hence were not known before.

\prclm{Corollary}{smallm}{Let $d\ge3$ and $\varepsilon \ge 1$
be integers. Then \ir{introconj}(c) holds for each of the
following values of $n$ and $m$, whenever $m$ is an integer:
\item{(a)} $n=d^2 + 2\varepsilon$ and
$m=d(d\pm1)/(2\varepsilon)$, in which case
$\alpha=\alpha(n;m)=md+d-1/2\pm1/2$ and the minimal free
resolution of $I(n;m)$ is
$$0\to R[-\alpha-1]^{\alpha}\to R[-\alpha]^{\alpha+1}\to I(n;m)\to 0;$$
\item{(b)} $n=d^2 + 2\varepsilon$ and
$m=(d(d\pm1)/2-1)/\varepsilon$, in which case
$\alpha=\alpha(n;m)=md+d-3/2\pm1/2$ and the minimal free
resolution of $I(n;m)$ is
$$0\to R[-\alpha-2]^{b+m}\to
R[-\alpha-1]^{b}\oplus R[-\alpha]^{m+1}
\to I(n;m)\to 0,$$
where $b=(m+1)(d-2)+1/2\pm1/2$;
\item{(c)} $n=d^2 + 2$ and $m=d^2+d(d\pm1)/2$, in which case
$\alpha=\alpha(n;m)=(m+1)d+d-3/2\pm1/2$ and the minimal free
resolution of $I(n;m)$ is
$$0\to R[-\alpha-2]^{a+b-1}\to R[-\alpha-1]^{b}\oplus R[-\alpha]^{a}
\to I(n;m)\to 0,$$
where $a=d(d\pm1)/2$ and $b=\alpha+2-d(d\pm1)$.}

\noindent{\bf Proof}: For case (a),
apply \ir{verhilb}(a2) for $m=d(d-1)/(2\varepsilon)$
and \ir{verhilb}(a3) for $m=d(d+1)/(2\varepsilon)$.
It turns out that $\alpha(n;m)>\tau(n;m)$ in these
cases, but it is well known that $I(n;m)$ is generated in degrees
$\tau(n;m)+1$ and less, hence in degree $\alpha(n;m)$, from which
it follows (see the displayed formula following
Definition 2.4 of \cite{refHHF}) that the minimal free
resolution is as claimed.

For cases (b) and (c), it turns out that
$\alpha(n;m)=\tau(n;m)$: for case (b),
apply \ir{verhilb}(a1-2),
resp., while for case (c),
apply \ir{verhilb}(b,c), resp.,
using $i=\varepsilon=1$.
To obtain the resolution, consider
${\bf m}=(m_1,\ldots,m_n)$, where
$m_1=m+1$ and $m_2=\cdots=m_n=m$, and
apply \ir{introthm} to $\alpha({\bf m})$
using $r=d^2+\varepsilon$. It turns out that
$\alpha({\bf m})>\alpha(n;m)$.
Now by Lemma 2.6(b) of \cite{refHHF} it follows that
\ir{introconj}(c) holds and that the minimal free
resolutions are as claimed (again, see the displayed formula following
Definition 2.4 of \cite{refHHF}).
\qed

When $n$ is an even square, \ir{sqrhilb}, together with
Theorem 5.1(a) of \cite{refHHF},
directly implies:

\prclm{Corollary}{sqrres}{Consider $n=\sigma^2$ general points
of \pr2, where $\sigma>3$ is even.
Let $k$ be any nonnegative integer, and let $m=x+k(\sigma-1)$,
where $x$ is an integer satisfying $\sigma/2-l_{\sigma}\le x\le \sigma/2$.
Then \ir{introconj}(c) holds for $I(n;m)$.}

Note that \ir{introcor}(c) is an immediate consequence of
the preceding result.

\irrnSection{Comparisons}{comps}
It is interesting to carry out some comparisons with previously known
bounds. Let, as before, $\alpha_c(n;m)$ denote the conjectural
value of $\alpha(n;m)$ and let $\tau_c(n;m)$ denote the
conjectural value of $\tau(n;m)$ (i.e., the values of
each assuming \ir{introconj}(b) holds).

Suppose that $rd(d+1)/2\le r^2<nd^2$; then for $m$
large enough the bound from
\ir{introthm}(c) is $\alpha(n;m)\ge 1+\lfloor(mr+g-1)/d\rfloor$.
This is better than the bound of
Corollary IV.1.1.2 of \cite{refsurv} (which generalizes the
main theorem of \cite{refnagconj}; see \cite{refSC}
for further generalizations and related results),
which is just $\alpha(n;m)\ge mr/d$. On the other hand,
suppose that $2r\ge n+d^2$. Then $r^2\ge nd^2$ (because
the arithmetic mean is never less than
the geometric mean), so the main theorem of \cite{refnagconj} applies
and gives $\alpha(n;m)\ge mnd/r$. Typically
\ir{introthm}(b) gives a better bound than this, but
if in addition $r$ divides $mn$, the bound from \ir{introthm}(b)
simplifies, also giving $\alpha(n;m)\ge mnd/r$.

In fact, for given $r$ and $d$, the bound in \cite{refnagconj}, and the
generalization given in Theorem IV.1.1.1 of \cite{refsurv}, can be
shown (see \cite{refsurv}) never to give a better bound on $\alpha$
than that given by the algorithms of \ir{algs}.
When $m$ is large enough compared with
$n$, \cite{refnagconj} shows its bound on $\alpha(n;m)$ is better
than those of \cite{refR}, and thus so are the bounds here.

When $m$ is not too large compared with $n$,
the bounds on $\alpha$ given by \ir{introthm}, like the bound
given by the unloading algorithm of \cite{refR}, are among the few
that sometimes give bounds better than the bound
$\alpha(n;m) \ge \lfloor m\sqrt{n}\rfloor+1$
conjectured in \cite{refN}. Consider, for example,
$n=1000$ and $m=13$: \cite{refR} gives $\alpha(n;m)\ge 421$
and \ir{introthm}, using $r=981$ and $d=31$, gives
$\alpha(n;m)\ge 424$, whereas $\lfloor m\sqrt{n}\rfloor+1=412$;
$\alpha_c(n;m)$ is 426 in this case. See \ir{vernag} for
more examples. Moreover, \ir{introthm} is
the only result that we know
which sometimes determines $\alpha(n;m)$ exactly
even for $m$ reasonably large compared to $n$,
when $n$ is not a square.

Here are some comparisons for $\tau$. Bounds on $\tau(n;m)$
given by Hirschowitz \cite{refHi}, Gimigliano \cite{refGi}
and Catalisano \cite{refCat} are on the order of $m\sqrt{2n}$.
Thus, for sufficiently large $m$,
the bound $\tau(n;m)\le m\lceil \sqrt{n}\rceil+
\lceil(\lceil\sqrt{n}\rceil-3)/2\rceil$ given in
\cite{refHHF} for $n>9$ is better. In fact, \cite{refHHF} shows that
$\tau(n;m)\le m\lceil \sqrt{n}\rceil+
\lceil(\lceil\sqrt{n}\rceil-3)/2\rceil$ is an equality
when $n>9$ is a square and $m$ is sufficiently large.
However, when $n$ is a square, \ir{introthm}(a),
using $d^2=r=n$, also gives
this bound (this is to be expected, since the method
we use is based on the method used in \cite{refHHF}),
and when $n$ is not a square, \ir{introthm}(a),
using $d=\lceil \sqrt{n}\rceil$
and $r=n$, gives a bound that is
less than or equal to that of \cite{refHHF} (although never more than
2 smaller). But one can also apply \ir{introthm}(b) using
other values of $r$ and $d$, and often do much better.
In addition, as was pointed out for $\alpha$ above,
\ir{introthm} is the only result that we know
that sometimes determines $\tau(n;m)$ exactly
for values of $m$ and $n$ that can be large,
even when $n$ is not a square.

Other bounds on $\tau$ have also been given.
Bounds given by Xu \cite{refXb} and Ballico \cite{refB}
are on the order of $m\sqrt{n}$, but nonethless
the bound from \cite{refHHF} (and hence \ir{introthm})
is better than Xu's when $n$ is large enough and
better than Ballico's when $m$ is large enough.
For large $m$, the bound given in
\cite{refRb} is also better than those of \cite{refB} and
\cite{refXb}, and by an argument similar to
the one used in \cite{refnagconj}
to compare the bounds on $\alpha$, the bounds here on
$\tau(n;m)$ are better than those of \cite{refRb} when $m$ is
large enough compared with $n$.

For example, for $n=190$ and $m=100$, then $\tau_c(n;m)=1384$,
while \ir{introthm}(b), using $r=180$ and $d=13$, gives
$\tau(n;m)\le 1390$, and we have in addition:
\item{$\bullet$} $\tau(Z)\le 1957$ from \cite{refHi},
\item{$\bullet$} $\tau(Z)\le 1900$ from \cite{refGi},
\item{$\bullet$} $\tau(Z)\le 1899$ from \cite{refCat},
\item{$\bullet$} $\tau(Z)\le 1487$ from \cite{refB},
\item{$\bullet$} $\tau(Z)\le 1465$ from \cite{refXb},
\item{$\bullet$} $\tau(Z)\le 1440$ from \cite{refRb} and
\item{$\bullet$} $\tau(Z)\le 1406$ from \cite{refHHF}.

For a different perspective, we close with
some graphs which show our results in certain ranges.
Figure 1 is a graphical
representation of \ir{verhilb} and \ir{sqrhilb}, while
Figure 2 is a graphical
representation of \ir{smallm} and \ir{sqrres}.
Figure 3 shows all $(n,m)$ for which
\ir{introthm}, using 
$d=\lfloor\sqrt{n}\rfloor$
and $r=\lfloor d\sqrt{n}\rfloor$, implies \ir{introconj}(a).
It was this graph that led us to
the statement of \ir{introcor}(a).
Figure 4 shows for comparison all $(n,m)$ for which
\ir{introthm} implies \ir{introconj}(a), but this time 
using $d=\lfloor\sqrt{n}\rfloor$
and $r=\lfloor(n+d^2)/2\rfloor$. This choice of $r$ and $d$
gives a higher density of cases for which we can conclude that 
\ir{introconj}(a) is true, 
but the graph has a very complicated structure
which does not seem to suggest any simply
stated result.

\vfill\eject

\def\xscale{1}   
\def\yscale{.8}   
\def\xmin{-5}   
\def\ymin{-5}    
\def\xmax{220}     
\def\ymax{220}     
\mywidth220pt     
\myheight220pt   
\penwidth1pt 
\penheight1pt 

\ \vskip-.5in\hskip1in\hbox{\vbox{      
\parindent0in                        
\putincaption 20 {\leftskip-1in\hsize4in\noindent{{\bf Figure 1}:
Graph of all $(n,m)$ with 
$10\le n\le \xmax$ and $0\le m\le \ymax$
such that \ir{verhilb} and \ir{sqrhilb} determine 
the Hilbert function of $I(n;m)$ of $n$ points of multiplicity $m$ 
in ${\bf C}P^2$.}}
\putinaxes  $n$ $m$
\putinticx 5 {100} {100} 
\putinticx 5 {200} {200} 
\putinticy 5 -5 {100} {100} 
\putinticy 5 -5 {200} {200}




\plt 8 11 1 \dt 
\plt 7 11 11 \dt 
\plt 0 11 21 \dt 
\plt 1 11 23 \dt 
\plt 1 11 26 \dt 
\plt 0 11 33 \dt 
\plt 0 11 36 \dt 
\plt 3 13 1 \dt 
\plt 1 13 6 \dt 
\plt 0 13 9 \dt 
\plt 2 15 1 \dt 
\plt 1 15 5 \dt 
\plt 1 17 1 \dt 
\plt 11 18 3 \dt 
\plt 10 18 20 \dt 
\plt 1 18 37 \dt 
\plt 2 18 40 \dt 
\plt 2 18 44 \dt 
\plt 0 18 54 \dt 
\plt 1 18 57 \dt 
\plt 1 18 61 \dt 
\plt 0 18 74 \dt 
\plt 0 18 78 \dt 
\plt 0 19 1 \dt 
\plt 5 20 2 \dt 
\plt 4 20 11 \dt 
\plt 0 20 21 \dt 
\plt 0 20 23 \dt 
\plt 0 21 1 \dt 
\plt 3 22 1 \dt 
\plt 0 22 10 \dt 
\plt 0 23 1 \dt 
\plt 2 24 1 \dt 
\plt 0 25 1 \dt 
\plt 1 26 1 \dt 
\plt 0 27 1 \dt 
\plt 14 27 6 \dt 
\plt 13 27 32 \dt 
\plt 2 27 58 \dt 
\plt 3 27 62 \dt 
\plt 3 27 67 \dt 
\plt 1 27 84 \dt 
\plt 2 27 88 \dt 
\plt 2 27 93 \dt 
\plt 0 27 110 \dt 
\plt 1 27 114 \dt 
\plt 1 27 119 \dt 
\plt 0 27 140 \dt 
\plt 0 27 145 \dt 
\plt 1 28 1 \dt 
\plt 7 29 3 \dt 
\plt 0 29 17 \dt 
\plt 3 29 19 \dt 
\plt 0 29 30 \dt 
\plt 0 29 32 \dt 
\plt 0 29 35 \dt 
\plt 1 30 1 \dt 
\plt 4 31 2 \dt 
\plt 0 31 13 \dt 
\plt 0 31 15 \dt 
\plt 0 32 1 \dt 
\plt 3 33 2 \dt 
\plt 1 33 10 \dt 
\plt 0 34 1 \dt 
\plt 2 35 2 \dt 
\plt 1 35 8 \dt 
\plt 0 36 1 \dt 
\plt 2 37 1 \dt 
\plt 0 38 1 \dt 
\plt 17 38 10 \dt 
\plt 16 38 47 \dt 
\plt 3 38 84 \dt 
\plt 4 38 89 \dt 
\plt 4 38 95 \dt 
\plt 2 38 121 \dt 
\plt 3 38 126 \dt 
\plt 3 38 132 \dt 
\plt 1 38 158 \dt 
\plt 2 38 163 \dt 
\plt 2 38 169 \dt 
\plt 0 38 195 \dt 
\plt 1 38 200 \dt 
\plt 1 38 206 \dt 
\plt 1 39 1 \dt 
\plt 0 40 1 \dt 
\plt 8 40 5 \dt 
\plt 1 40 24 \dt 
\plt 1 40 27 \dt 
\plt 1 40 30 \dt 
\plt 0 40 43 \dt 
\plt 0 40 46 \dt 
\plt 0 40 49 \dt 
\plt 1 41 1 \dt 
\plt 0 42 1 \dt 
\plt 5 42 4 \dt 
\plt 4 42 17 \dt 
\plt 0 42 31 \dt 
\plt 0 42 33 \dt 
\plt 1 43 1 \dt 
\plt 0 44 1 \dt 
\plt 3 44 3 \dt 
\plt 0 44 14 \dt 
\plt 1 45 1 \dt 
\plt 3 46 2 \dt 
\plt 0 47 1 \dt 
\plt 2 48 2 \dt 
\plt 0 49 1 \dt 
\plt 1 50 2 \dt 
\plt 0 51 1 \dt 
\plt 20 51 15 \dt 
\plt 19 51 65 \dt 
\plt 4 51 115 \dt 
\plt 5 51 121 \dt 
\plt 5 51 128 \dt 
\plt 3 51 165 \dt 
\plt 4 51 171 \dt 
\plt 4 51 178 \dt 
\plt 2 51 215 \dt 
\plt 1 52 2 \dt 
\plt 0 53 1 \dt 
\plt 9 53 8 \dt 
\plt 5 53 33 \dt 
\plt 2 53 40 \dt 
\plt 0 53 59 \dt 
\plt 1 53 62 \dt 
\plt 1 53 65 \dt 
\plt 0 53 84 \dt 
\plt 0 53 87 \dt 
\plt 0 53 91 \dt 
\plt 1 54 2 \dt 
\plt 0 55 1 \dt 
\plt 6 55 5 \dt 
\plt 0 55 23 \dt 
\plt 0 55 25 \dt 
\plt 1 55 27 \dt 
\plt 0 55 42 \dt 
\plt 0 55 44 \dt 
\plt 1 56 1 \dt 
\plt 0 57 1 \dt 
\plt 4 57 4 \dt 
\plt 0 57 17 \dt 
\plt 0 57 19 \dt 
\plt 0 57 21 \dt 
\plt 1 58 1 \dt 
\plt 0 59 1 \dt 
\plt 4 59 3 \dt 
\plt 1 59 14 \dt 
\plt 1 60 1 \dt 
\plt 0 61 1 \dt 
\plt 2 61 3 \dt 
\plt 0 61 14 \dt 
\plt 1 62 1 \dt 
\plt 0 63 1 \dt 
\plt 2 63 3 \dt 
\plt 1 63 11 \dt 
\plt 0 64 1 \dt 
\plt 3 65 1 \dt 
\plt 0 66 1 \dt 
\plt 23 66 21 \dt 
\plt 22 66 86 \dt 
\plt 5 66 151 \dt 
\plt 6 66 158 \dt 
\plt 6 66 166 \dt 
\plt 4 66 216 \dt 
\plt 1 67 2 \dt 
\plt 0 68 1 \dt 
\plt 11 68 11 \dt 
\plt 10 68 44 \dt 
\plt 1 68 77 \dt 
\plt 2 68 80 \dt 
\plt 2 68 84 \dt 
\plt 0 68 110 \dt 
\plt 1 68 113 \dt 
\plt 1 68 117 \dt 
\plt 0 68 146 \dt 
\plt 0 68 150 \dt 
\plt 1 69 2 \dt 
\plt 0 70 1 \dt 
\plt 7 70 7 \dt 
\plt 0 70 30 \dt 
\plt 1 70 32 \dt 
\plt 1 70 35 \dt 
\plt 0 70 52 \dt 
\plt 0 70 54 \dt 
\plt 0 70 57 \dt 
\plt 1 71 2 \dt 
\plt 0 72 1 \dt 
\plt 5 72 6 \dt 
\plt 4 72 23 \dt 
\plt 0 72 41 \dt 
\plt 0 72 43 \dt 
\plt 0 73 2 \dt 
\plt 0 74 1 \dt 
\plt 3 74 5 \dt 
\plt 0 74 18 \dt 
\plt 1 74 20 \dt 
\plt 0 75 2 \dt 
\plt 0 76 1 \dt 
\plt 3 76 4 \dt 
\plt 0 76 18 \dt 
\plt 0 77 2 \dt 
\plt 0 78 1 \dt 
\plt 3 78 3 \dt 
\plt 1 79 1 \dt 
\plt 0 80 1 \dt 
\plt 2 80 3 \dt 
\plt 1 81 1 \dt 
\plt 0 82 1 \dt 
\plt 1 82 3 \dt 
\plt 1 83 1 \dt 
\plt 26 83 28 \dt 
\plt 25 83 110 \dt 
\plt 6 83 192 \dt 
\plt 7 83 200 \dt 
\plt 7 83 209 \dt 
\plt 0 84 1 \dt 
\plt 1 84 3 \dt 
\plt 0 85 1 \dt 
\plt 13 85 14 \dt 
\plt 2 85 56 \dt 
\plt 7 85 60 \dt 
\plt 2 85 97 \dt 
\plt 2 85 101 \dt 
\plt 2 85 106 \dt 
\plt 0 85 139 \dt 
\plt 1 85 143 \dt 
\plt 1 85 147 \dt 
\plt 0 85 180 \dt 
\plt 0 85 184 \dt 
\plt 0 85 189 \dt 
\plt 3 86 1 \dt 
\plt 0 87 1 \dt 
\plt 8 87 10 \dt 
\plt 7 87 38 \dt 
\plt 0 87 66 \dt 
\plt 1 87 68 \dt 
\plt 1 87 71 \dt 
\plt 0 87 96 \dt 
\plt 0 87 99 \dt 
\plt 2 88 1 \dt 
\plt 0 89 1 \dt 
\plt 6 89 7 \dt 
\plt 0 89 29 \dt 
\plt 0 89 31 \dt 
\plt 0 89 33 \dt 
\plt 0 89 54 \dt 
\plt 2 90 1 \dt 
\plt 0 91 1 \dt 
\plt 4 91 6 \dt 
\plt 0 91 23 \dt 
\plt 0 91 25 \dt 
\plt 0 91 27 \dt 
\plt 1 92 2 \dt 
\plt 0 93 1 \dt 
\plt 4 93 5 \dt 
\plt 1 93 21 \dt 
\plt 0 94 2 \dt 
\plt 0 95 1 \dt 
\plt 3 95 4 \dt 
\plt 1 95 18 \dt 
\plt 0 96 2 \dt 
\plt 0 97 1 \dt 
\plt 2 97 4 \dt 
\plt 0 98 2 \dt 
\plt 0 99 1 \dt 
\plt 2 99 4 \dt 
\plt 1 99 14 \dt 
\plt 0 100 2 \dt 
\plt 0 101 1 \dt 
\plt 2 101 3 \dt 
\plt 0 102 2 \dt 
\plt 29 102 36 \dt 
\plt 28 102 137 \dt 
\plt 0 103 1 \dt 
\plt 1 103 3 \dt 
\plt 0 104 2 \dt 
\plt 14 104 18 \dt 
\plt 3 104 69 \dt 
\plt 3 104 74 \dt 
\plt 3 104 79 \dt 
\plt 2 104 120 \dt 
\plt 2 104 125 \dt 
\plt 2 104 130 \dt 
\plt 1 104 171 \dt 
\plt 1 104 176 \dt 
\plt 1 104 181 \dt 
\plt 0 105 1 \dt 
\plt 1 105 3 \dt 
\plt 1 106 1 \dt 
\plt 9 106 12 \dt 
\plt 1 106 47 \dt 
\plt 1 106 50 \dt 
\plt 2 106 53 \dt 
\plt 0 106 81 \dt 
\plt 1 106 84 \dt 
\plt 1 106 87 \dt 
\plt 0 106 115 \dt 
\plt 0 106 118 \dt 
\plt 0 107 1 \dt 
\plt 1 107 3 \dt 
\plt 1 108 1 \dt 
\plt 7 108 9 \dt 
\plt 1 108 35 \dt 
\plt 0 108 38 \dt 
\plt 1 108 40 \dt 
\plt 0 108 61 \dt 
\plt 0 108 66 \dt 
\plt 2 109 1 \dt 
\plt 0 110 1 \dt 
\plt 5 110 8 \dt 
\plt 4 110 29 \dt 
\plt 0 110 51 \dt 
\plt 0 110 53 \dt 
\plt 2 111 1 \dt 
\plt 0 112 1 \dt 
\plt 4 112 6 \dt 
\plt 0 112 24 \dt 
\plt 0 112 27 \dt 
\plt 2 113 1 \dt 
\plt 0 114 1 \dt 
\plt 3 114 6 \dt 
\plt 0 114 22 \dt 
\plt 2 115 1 \dt 
\plt 0 116 1 \dt 
\plt 3 116 5 \dt 
\plt 0 116 18 \dt 
\plt 2 117 1 \dt 
\plt 0 118 1 \dt 
\plt 3 118 4 \dt 
\plt 1 119 1 \dt 
\plt 0 120 1 \dt 
\plt 2 120 4 \dt 
\plt 0 121 2 \dt 
\plt 0 122 1 \dt 
\plt 1 122 4 \dt 
\plt 0 123 2 \dt 
\plt 32 123 45 \dt 
\plt 31 123 167 \dt 
\plt 0 124 1 \dt 
\plt 2 124 3 \dt 
\plt 0 125 2 \dt 
\plt 15 125 23 \dt 
\plt 9 125 84 \dt 
\plt 4 125 95 \dt 
\plt 2 125 146 \dt 
\plt 3 125 151 \dt 
\plt 3 125 156 \dt 
\plt 2 125 207 \dt 
\plt 2 125 212 \dt 
\plt 2 125 218 \dt 
\plt 0 126 1 \dt 
\plt 2 126 3 \dt 
\plt 0 127 2 \dt 
\plt 10 127 15 \dt 
\plt 1 127 57 \dt 
\plt 2 127 60 \dt 
\plt 2 127 64 \dt 
\plt 1 127 98 \dt 
\plt 1 127 101 \dt 
\plt 1 127 105 \dt 
\plt 0 127 139 \dt 
\plt 0 127 143 \dt 
\plt 0 127 146 \dt 
\plt 0 128 1 \dt 
\plt 1 128 3 \dt 
\plt 0 129 2 \dt 
\plt 7 129 12 \dt 
\plt 3 129 43 \dt 
\plt 1 129 48 \dt 
\plt 0 129 74 \dt 
\plt 0 129 77 \dt 
\plt 0 129 79 \dt 
\plt 0 130 1 \dt 
\plt 1 130 3 \dt 
\plt 0 131 2 \dt 
\plt 6 131 9 \dt 
\plt 0 131 35 \dt 
\plt 0 131 37 \dt 
\plt 0 131 39 \dt 
\plt 0 132 1 \dt 
\plt 1 132 3 \dt 
\plt 0 133 2 \dt 
\plt 4 133 8 \dt 
\plt 0 133 29 \dt 
\plt 0 133 31 \dt 
\plt 0 133 33 \dt 
\plt 0 134 1 \dt 
\plt 0 134 3 \dt 
\plt 0 135 2 \dt 
\plt 4 135 7 \dt 
\plt 0 135 25 \dt 
\plt 0 135 28 \dt 
\plt 2 136 1 \dt 
\plt 0 137 1 \dt 
\plt 3 137 6 \dt 
\plt 1 137 22 \dt 
\plt 2 138 1 \dt 
\plt 0 139 1 \dt 
\plt 3 139 5 \dt 
\plt 0 139 22 \dt 
\plt 2 140 1 \dt 
\plt 0 141 1 \dt 
\plt 2 141 5 \dt 
\plt 2 142 1 \dt 
\plt 0 143 1 \dt 
\plt 2 143 5 \dt 
\plt 1 143 17 \dt 
\plt 1 144 1 \dt 
\plt 0 145 1 \dt 
\plt 2 145 4 \dt 
\plt 1 146 1 \dt 
\plt 35 146 55 \dt 
\plt 20 146 200 \dt 
\plt 0 147 1 \dt 
\plt 1 147 4 \dt 
\plt 1 148 1 \dt 
\plt 17 148 28 \dt 
\plt 16 148 101 \dt 
\plt 3 148 174 \dt 
\plt 4 148 179 \dt 
\plt 4 148 185 \dt 
\plt 0 149 1 \dt 
\plt 1 149 4 \dt 
\plt 1 150 1 \dt 
\plt 11 150 19 \dt 
\plt 10 150 68 \dt 
\plt 1 150 117 \dt 
\plt 2 150 120 \dt 
\plt 2 150 124 \dt 
\plt 0 150 166 \dt 
\plt 1 150 169 \dt 
\plt 1 150 173 \dt 
\plt 0 150 218 \dt 
\plt 0 151 1 \dt 
\plt 2 151 3 \dt 
\plt 1 152 1 \dt 
\plt 8 152 14 \dt 
\plt 1 152 51 \dt 
\plt 1 152 54 \dt 
\plt 1 152 57 \dt 
\plt 0 152 88 \dt 
\plt 0 152 91 \dt 
\plt 0 152 94 \dt 
\plt 0 153 1 \dt 
\plt 1 153 3 \dt 
\plt 0 154 2 \dt 
\plt 7 154 11 \dt 
\plt 3 154 41 \dt 
\plt 0 154 46 \dt 
\plt 0 154 73 \dt 
\plt 0 155 1 \dt 
\plt 1 155 3 \dt 
\plt 0 156 2 \dt 
\plt 5 156 10 \dt 
\plt 4 156 35 \dt 
\plt 0 156 61 \dt 
\plt 0 156 63 \dt 
\plt 0 157 1 \dt 
\plt 1 157 3 \dt 
\plt 0 158 2 \dt 
\plt 4 158 8 \dt 
\plt 1 158 30 \dt 
\plt 0 158 33 \dt 
\plt 0 159 1 \dt 
\plt 1 159 3 \dt 
\plt 0 160 2 \dt 
\plt 4 160 7 \dt 
\plt 0 160 26 \dt 
\plt 0 160 29 \dt 
\plt 0 161 1 \dt 
\plt 0 161 3 \dt 
\plt 0 162 2 \dt 
\plt 3 162 7 \dt 
\plt 0 162 26 \dt 
\plt 0 163 1 \dt 
\plt 0 163 3 \dt 
\plt 0 164 2 \dt 
\plt 3 164 6 \dt 
\plt 1 164 21 \dt 
\plt 0 165 1 \dt 
\plt 0 165 3 \dt 
\plt 3 166 5 \dt 
\plt 2 167 1 \dt 
\plt 2 168 5 \dt 
\plt 2 169 1 \dt 
\plt 1 170 5 \dt 
\plt 2 171 1 \dt 
\plt 38 171 66 \dt 
\plt 0 172 1 \dt 
\plt 2 172 4 \dt 
\plt 1 173 1 \dt 
\plt 19 173 33 \dt 
\plt 4 173 119 \dt 
\plt 11 173 125 \dt 
\plt 4 173 204 \dt 
\plt 4 173 210 \dt 
\plt 3 173 217 \dt 
\plt 0 174 1 \dt 
\plt 2 174 4 \dt 
\plt 1 175 1 \dt 
\plt 12 175 22 \dt 
\plt 2 175 80 \dt 
\plt 2 175 84 \dt 
\plt 3 175 88 \dt 
\plt 1 175 137 \dt 
\plt 2 175 141 \dt 
\plt 2 175 145 \dt 
\plt 1 175 194 \dt 
\plt 1 175 198 \dt 
\plt 0 175 203 \dt 
\plt 0 176 1 \dt 
\plt 1 176 4 \dt 
\plt 1 177 1 \dt 
\plt 9 177 17 \dt 
\plt 1 177 60 \dt 
\plt 5 177 63 \dt 
\plt 1 177 103 \dt 
\plt 1 177 106 \dt 
\plt 0 177 110 \dt 
\plt 0 178 1 \dt 
\plt 1 178 4 \dt 
\plt 1 179 1 \dt 
\plt 6 179 14 \dt 
\plt 1 179 48 \dt 
\plt 3 179 51 \dt 
\plt 0 179 83 \dt 
\plt 0 179 88 \dt 
\plt 0 180 1 \dt 
\plt 1 180 4 \dt 
\plt 1 181 1 \dt 
\plt 6 181 11 \dt 
\plt 0 181 41 \dt 
\plt 0 181 43 \dt 
\plt 0 181 45 \dt 
\plt 0 182 1 \dt 
\plt 1 182 3 \dt 
\plt 1 183 1 \dt 
\plt 4 183 10 \dt 
\plt 0 183 35 \dt 
\plt 0 183 37 \dt 
\plt 0 183 39 \dt 
\plt 0 184 1 \dt 
\plt 1 184 3 \dt 
\plt 1 185 1 \dt 
\plt 4 185 9 \dt 
\plt 0 185 32 \dt 
\plt 0 185 34 \dt 
\plt 0 186 1 \dt 
\plt 1 186 3 \dt 
\plt 1 187 1 \dt 
\plt 3 187 8 \dt 
\plt 0 187 30 \dt 
\plt 0 188 1 \dt 
\plt 1 188 3 \dt 
\plt 1 189 1 \dt 
\plt 3 189 7 \dt 
\plt 1 189 26 \dt 
\plt 0 190 1 \dt 
\plt 0 190 3 \dt 
\plt 0 191 2 \dt 
\plt 3 191 6 \dt 
\plt 0 192 1 \dt 
\plt 0 192 3 \dt 
\plt 0 193 2 \dt 
\plt 2 193 6 \dt 
\plt 0 194 1 \dt 
\plt 0 194 3 \dt 
\plt 0 195 2 \dt 
\plt 2 195 6 \dt 
\plt 1 195 20 \dt 
\plt 0 196 1 \dt 
\plt 0 196 3 \dt 
\plt 0 197 2 \dt 
\plt 2 197 5 \dt 
\plt 0 198 1 \dt 
\plt 0 198 3 \dt 
\plt 41 198 78 \dt 
\plt 1 199 5 \dt 
\plt 2 200 1 \dt 
\plt 20 200 39 \dt 
\plt 5 200 138 \dt 
\plt 5 200 145 \dt 
\plt 5 200 152 \dt 
\plt 1 201 5 \dt 
\plt 2 202 1 \dt 
\plt 13 202 26 \dt 
\plt 2 202 93 \dt 
\plt 3 202 97 \dt 
\plt 3 202 102 \dt 
\plt 2 202 159 \dt 
\plt 2 202 163 \dt 
\plt 2 202 168 \dt 
\plt 2 203 4 \dt 
\plt 2 204 1 \dt 
\plt 9 204 20 \dt 
\plt 1 204 70 \dt 
\plt 2 204 73 \dt 
\plt 1 204 77 \dt 
\plt 0 204 120 \dt 
\plt 1 204 123 \dt 
\plt 0 204 127 \dt 
\plt 0 204 173 \dt 
\plt 1 205 4 \dt 
\plt 1 206 1 \dt 
\plt 7 206 16 \dt 
\plt 1 206 56 \dt 
\plt 1 206 59 \dt 
\plt 1 206 62 \dt 
\plt 0 206 96 \dt 
\plt 0 206 99 \dt 
\plt 0 206 102 \dt 
\plt 1 207 4 \dt 
\plt 1 208 1 \dt 
\plt 6 208 13 \dt 
\plt 0 208 47 \dt 
\plt 1 208 49 \dt 
\plt 0 208 52 \dt 
\plt 0 208 85 \dt 
\plt 1 209 4 \dt 
\plt 1 210 1 \dt 
\plt 5 210 12 \dt 
\plt 4 210 41 \dt 
\plt 0 210 71 \dt 
\plt 0 210 73 \dt 
\plt 0 211 1 \dt 
\plt 0 211 4 \dt 
\plt 1 212 1 \dt 
\plt 4 212 10 \dt 
\plt 0 212 37 \dt 
\plt 0 212 39 \dt 
\plt 0 213 1 \dt 
\plt 1 213 3 \dt 
\plt 1 214 1 \dt 
\plt 4 214 9 \dt 
\plt 0 214 33 \dt 
\plt 0 214 35 \dt 
\plt 0 215 1 \dt 
\plt 1 215 3 \dt 
\plt 1 216 1 \dt 
\plt 3 216 8 \dt 
\plt 0 216 30 \dt 
\plt 0 217 1 \dt 
\plt 1 217 3 \dt 
\plt 1 218 1 \dt 
\plt 2 218 8 \dt 
\plt 0 218 26 \dt 
\plt 0 219 1 \dt 
\plt 1 219 3 \dt 
\plt 1 220 1 \dt 
\plt 2 220 7 \dt 
\plt 0 220 25 \dt 
\plt 1 16 1 \dt
\plt 1 16 4 \dt
\plt 1 16 7 \dt
\plt 1 16 10 \dt
\plt 1 16 13 \dt
\plt 1 16 16 \dt
\plt 1 16 19 \dt
\plt 1 16 22 \dt
\plt 1 16 25 \dt
\plt 1 16 28 \dt
\plt 1 16 31 \dt
\plt 1 16 34 \dt
\plt 1 16 37 \dt
\plt 1 16 40 \dt
\plt 1 16 43 \dt
\plt 1 16 46 \dt
\plt 1 16 49 \dt
\plt 1 16 52 \dt
\plt 1 16 55 \dt
\plt 1 16 58 \dt
\plt 1 16 61 \dt
\plt 1 16 64 \dt
\plt 1 16 67 \dt
\plt 1 16 70 \dt
\plt 1 16 73 \dt
\plt 1 16 76 \dt
\plt 1 16 79 \dt
\plt 1 16 82 \dt
\plt 1 16 85 \dt
\plt 1 16 88 \dt
\plt 1 16 91 \dt
\plt 1 16 94 \dt
\plt 1 16 97 \dt
\plt 1 16 100 \dt
\plt 1 16 103 \dt
\plt 1 16 106 \dt
\plt 1 16 109 \dt
\plt 1 16 112 \dt
\plt 1 16 115 \dt
\plt 1 16 118 \dt
\plt 1 16 121 \dt
\plt 1 16 124 \dt
\plt 1 16 127 \dt
\plt 1 16 130 \dt
\plt 1 16 133 \dt
\plt 1 16 136 \dt
\plt 1 16 139 \dt
\plt 1 16 142 \dt
\plt 1 16 145 \dt
\plt 1 16 148 \dt
\plt 1 16 151 \dt
\plt 1 16 154 \dt
\plt 1 16 157 \dt
\plt 1 16 160 \dt
\plt 1 16 163 \dt
\plt 1 16 166 \dt
\plt 1 16 169 \dt
\plt 1 16 172 \dt
\plt 1 16 175 \dt
\plt 1 16 178 \dt
\plt 1 16 181 \dt
\plt 1 16 184 \dt
\plt 1 16 187 \dt
\plt 1 16 190 \dt
\plt 1 16 193 \dt
\plt 1 16 196 \dt
\plt 1 16 199 \dt
\plt 1 16 202 \dt
\plt 1 16 205 \dt
\plt 1 16 208 \dt
\plt 1 16 211 \dt
\plt 1 16 214 \dt
\plt 1 16 217 \dt
\plt 2 25 1 \dt
\plt 2 25 5 \dt
\plt 2 25 9 \dt
\plt 2 25 13 \dt
\plt 2 25 17 \dt
\plt 2 25 21 \dt
\plt 2 25 25 \dt
\plt 2 25 29 \dt
\plt 2 25 33 \dt
\plt 2 25 37 \dt
\plt 2 25 41 \dt
\plt 2 25 45 \dt
\plt 2 25 49 \dt
\plt 2 25 53 \dt
\plt 2 25 57 \dt
\plt 2 25 61 \dt
\plt 2 25 65 \dt
\plt 2 25 69 \dt
\plt 2 25 73 \dt
\plt 2 25 77 \dt
\plt 2 25 81 \dt
\plt 2 25 85 \dt
\plt 2 25 89 \dt
\plt 2 25 93 \dt
\plt 2 25 97 \dt
\plt 2 25 101 \dt
\plt 2 25 105 \dt
\plt 2 25 109 \dt
\plt 2 25 113 \dt
\plt 2 25 117 \dt
\plt 2 25 121 \dt
\plt 2 25 125 \dt
\plt 2 25 129 \dt
\plt 2 25 133 \dt
\plt 2 25 137 \dt
\plt 2 25 141 \dt
\plt 2 25 145 \dt
\plt 2 25 149 \dt
\plt 2 25 153 \dt
\plt 2 25 157 \dt
\plt 2 25 161 \dt
\plt 2 25 165 \dt
\plt 2 25 169 \dt
\plt 2 25 173 \dt
\plt 2 25 177 \dt
\plt 2 25 181 \dt
\plt 2 25 185 \dt
\plt 2 25 189 \dt
\plt 2 25 193 \dt
\plt 2 25 197 \dt
\plt 2 25 201 \dt
\plt 2 25 205 \dt
\plt 2 25 209 \dt
\plt 2 25 213 \dt
\plt 2 25 217 \dt
\plt 2 25 221 \dt
\plt 2 36 1 \dt
\plt 2 36 6 \dt
\plt 2 36 11 \dt
\plt 2 36 16 \dt
\plt 2 36 21 \dt
\plt 2 36 26 \dt
\plt 2 36 31 \dt
\plt 2 36 36 \dt
\plt 2 36 41 \dt
\plt 2 36 46 \dt
\plt 2 36 51 \dt
\plt 2 36 56 \dt
\plt 2 36 61 \dt
\plt 2 36 66 \dt
\plt 2 36 71 \dt
\plt 2 36 76 \dt
\plt 2 36 81 \dt
\plt 2 36 86 \dt
\plt 2 36 91 \dt
\plt 2 36 96 \dt
\plt 2 36 101 \dt
\plt 2 36 106 \dt
\plt 2 36 111 \dt
\plt 2 36 116 \dt
\plt 2 36 121 \dt
\plt 2 36 126 \dt
\plt 2 36 131 \dt
\plt 2 36 136 \dt
\plt 2 36 141 \dt
\plt 2 36 146 \dt
\plt 2 36 151 \dt
\plt 2 36 156 \dt
\plt 2 36 161 \dt
\plt 2 36 166 \dt
\plt 2 36 171 \dt
\plt 2 36 176 \dt
\plt 2 36 181 \dt
\plt 2 36 186 \dt
\plt 2 36 191 \dt
\plt 2 36 196 \dt
\plt 2 36 201 \dt
\plt 2 36 206 \dt
\plt 2 36 211 \dt
\plt 2 36 216 \dt
\plt 2 36 221 \dt
\plt 3 49 1 \dt
\plt 3 49 7 \dt
\plt 3 49 13 \dt
\plt 3 49 19 \dt
\plt 3 49 25 \dt
\plt 3 49 31 \dt
\plt 3 49 37 \dt
\plt 3 49 43 \dt
\plt 3 49 49 \dt
\plt 3 49 55 \dt
\plt 3 49 61 \dt
\plt 3 49 67 \dt
\plt 3 49 73 \dt
\plt 3 49 79 \dt
\plt 3 49 85 \dt
\plt 3 49 91 \dt
\plt 3 49 97 \dt
\plt 3 49 103 \dt
\plt 3 49 109 \dt
\plt 3 49 115 \dt
\plt 3 49 121 \dt
\plt 3 49 127 \dt
\plt 3 49 133 \dt
\plt 3 49 139 \dt
\plt 3 49 145 \dt
\plt 3 49 151 \dt
\plt 3 49 157 \dt
\plt 3 49 163 \dt
\plt 3 49 169 \dt
\plt 3 49 175 \dt
\plt 3 49 181 \dt
\plt 3 49 187 \dt
\plt 3 49 193 \dt
\plt 3 49 199 \dt
\plt 3 49 205 \dt
\plt 3 49 211 \dt
\plt 3 49 217 \dt
\plt 2 64 2 \dt
\plt 2 64 9 \dt
\plt 2 64 16 \dt
\plt 2 64 23 \dt
\plt 2 64 30 \dt
\plt 2 64 37 \dt
\plt 2 64 44 \dt
\plt 2 64 51 \dt
\plt 2 64 58 \dt
\plt 2 64 65 \dt
\plt 2 64 72 \dt
\plt 2 64 79 \dt
\plt 2 64 86 \dt
\plt 2 64 93 \dt
\plt 2 64 100 \dt
\plt 2 64 107 \dt
\plt 2 64 114 \dt
\plt 2 64 121 \dt
\plt 2 64 128 \dt
\plt 2 64 135 \dt
\plt 2 64 142 \dt
\plt 2 64 149 \dt
\plt 2 64 156 \dt
\plt 2 64 163 \dt
\plt 2 64 170 \dt
\plt 2 64 177 \dt
\plt 2 64 184 \dt
\plt 2 64 191 \dt
\plt 2 64 198 \dt
\plt 2 64 205 \dt
\plt 2 64 212 \dt
\plt 2 64 219 \dt
\plt 3 81 2 \dt
\plt 3 81 10 \dt
\plt 3 81 18 \dt
\plt 3 81 26 \dt
\plt 3 81 34 \dt
\plt 3 81 42 \dt
\plt 3 81 50 \dt
\plt 3 81 58 \dt
\plt 3 81 66 \dt
\plt 3 81 74 \dt
\plt 3 81 82 \dt
\plt 3 81 90 \dt
\plt 3 81 98 \dt
\plt 3 81 106 \dt
\plt 3 81 114 \dt
\plt 3 81 122 \dt
\plt 3 81 130 \dt
\plt 3 81 138 \dt
\plt 3 81 146 \dt
\plt 3 81 154 \dt
\plt 3 81 162 \dt
\plt 3 81 170 \dt
\plt 3 81 178 \dt
\plt 3 81 186 \dt
\plt 3 81 194 \dt
\plt 3 81 202 \dt
\plt 3 81 210 \dt
\plt 3 81 218 \dt
\plt 2 100 3 \dt
\plt 2 100 12 \dt
\plt 2 100 21 \dt
\plt 2 100 30 \dt
\plt 2 100 39 \dt
\plt 2 100 48 \dt
\plt 2 100 57 \dt
\plt 2 100 66 \dt
\plt 2 100 75 \dt
\plt 2 100 84 \dt
\plt 2 100 93 \dt
\plt 2 100 102 \dt
\plt 2 100 111 \dt
\plt 2 100 120 \dt
\plt 2 100 129 \dt
\plt 2 100 138 \dt
\plt 2 100 147 \dt
\plt 2 100 156 \dt
\plt 2 100 165 \dt
\plt 2 100 174 \dt
\plt 2 100 183 \dt
\plt 2 100 192 \dt
\plt 2 100 201 \dt
\plt 2 100 210 \dt
\plt 2 100 219 \dt
\plt 4 121 2 \dt
\plt 4 121 12 \dt
\plt 4 121 22 \dt
\plt 4 121 32 \dt
\plt 4 121 42 \dt
\plt 4 121 52 \dt
\plt 4 121 62 \dt
\plt 4 121 72 \dt
\plt 4 121 82 \dt
\plt 4 121 92 \dt
\plt 4 121 102 \dt
\plt 4 121 112 \dt
\plt 4 121 122 \dt
\plt 4 121 132 \dt
\plt 4 121 142 \dt
\plt 4 121 152 \dt
\plt 4 121 162 \dt
\plt 4 121 172 \dt
\plt 4 121 182 \dt
\plt 4 121 192 \dt
\plt 4 121 202 \dt
\plt 4 121 212 \dt
\plt 3 144 3 \dt
\plt 3 144 14 \dt
\plt 3 144 25 \dt
\plt 3 144 36 \dt
\plt 3 144 47 \dt
\plt 3 144 58 \dt
\plt 3 144 69 \dt
\plt 3 144 80 \dt
\plt 3 144 91 \dt
\plt 3 144 102 \dt
\plt 3 144 113 \dt
\plt 3 144 124 \dt
\plt 3 144 135 \dt
\plt 3 144 146 \dt
\plt 3 144 157 \dt
\plt 3 144 168 \dt
\plt 3 144 179 \dt
\plt 3 144 190 \dt
\plt 3 144 201 \dt
\plt 3 144 212 \dt
\plt 4 169 3 \dt
\plt 4 169 15 \dt
\plt 4 169 27 \dt
\plt 4 169 39 \dt
\plt 4 169 51 \dt
\plt 4 169 63 \dt
\plt 4 169 75 \dt
\plt 4 169 87 \dt
\plt 4 169 99 \dt
\plt 4 169 111 \dt
\plt 4 169 123 \dt
\plt 4 169 135 \dt
\plt 4 169 147 \dt
\plt 4 169 159 \dt
\plt 4 169 171 \dt
\plt 4 169 183 \dt
\plt 4 169 195 \dt
\plt 4 169 207 \dt
\plt 4 169 219 \dt
\plt 3 196 4 \dt
\plt 3 196 17 \dt
\plt 3 196 30 \dt
\plt 3 196 43 \dt
\plt 3 196 56 \dt
\plt 3 196 69 \dt
\plt 3 196 82 \dt
\plt 3 196 95 \dt
\plt 3 196 108 \dt
\plt 3 196 121 \dt
\plt 3 196 134 \dt
\plt 3 196 147 \dt
\plt 3 196 160 \dt
\plt 3 196 173 \dt
\plt 3 196 186 \dt
\plt 3 196 199 \dt
\plt 3 196 212 \dt
}}

\def\xscale{1}   
\def\yscale{.8}   
\def\xmin{-5}   
\def\ymin{-5}    
\def\xmax{220}     
\def\ymax{220}     
\mywidth220pt     
\myheight220pt   
\penwidth1pt 
\penheight1pt 

\vskip3.5in\hskip1in\hbox{\vbox{      
\parindent0in                        
\putincaption 20 {\leftskip-1in\hsize4in\noindent{{\bf Figure 2}:
Graph of all $(n,m)$ with 
$10\le n\le \xmax$ and $0\le m\le \ymax$
such that \ir{smallm} and \ir{sqrres} determine 
the minimal free resolution of $I(n;m)$ of $n$ points of multiplicity $m$ 
in ${\bf C}P^2$.}}
\putinaxes  $n$ $m$
\putinticx 5 {100} {100} 
\putinticx 5 {200} {200} 
\putinticy 5 -5 {100} {100} 
\putinticy 5 -5 {200} {200}




\plt 1 11 2 \dt 
\plt 1 11 5 \dt 
\plt 0 11 12 \dt 
\plt 0 11 15 \dt 
\plt 0 13 1 \dt 
\plt 0 13 3 \dt 
\plt 1 15 1 \dt 
\plt 1 18 5 \dt 
\plt 1 18 9 \dt 
\plt 0 18 22 \dt 
\plt 0 18 26 \dt 
\plt 0 19 1 \dt 
\plt 0 20 3 \dt 
\plt 0 20 5 \dt 
\plt 0 21 1 \dt 
\plt 1 22 2 \dt 
\plt 1 26 1 \dt 
\plt 1 27 9 \dt 
\plt 1 27 14 \dt 
\plt 0 27 35 \dt 
\plt 0 27 40 \dt 
\plt 0 28 1 \dt 
\plt 0 29 5 \dt 
\plt 0 29 7 \dt 
\plt 0 31 3 \dt 
\plt 0 31 5 \dt 
\plt 0 34 1 \dt 
\plt 1 35 2 \dt 
\plt 0 36 1 \dt 
\plt 1 38 14 \dt 
\plt 1 38 20 \dt 
\plt 0 38 51 \dt 
\plt 0 38 57 \dt 
\plt 0 39 2 \dt 
\plt 0 40 7 \dt 
\plt 0 40 10 \dt 
\plt 0 42 5 \dt 
\plt 0 42 7 \dt 
\plt 0 43 1 \dt 
\plt 0 44 5 \dt 
\plt 0 45 1 \dt 
\plt 1 46 3 \dt 
\plt 1 50 2 \dt 
\plt 1 51 20 \dt 
\plt 1 51 27 \dt 
\plt 0 51 70 \dt 
\plt 0 51 77 \dt 
\plt 0 53 1 \dt 
\plt 0 53 10 \dt 
\plt 0 53 14 \dt 
\plt 0 55 1 \dt 
\plt 0 55 7 \dt 
\plt 0 55 9 \dt 
\plt 0 56 2 \dt 
\plt 0 57 5 \dt 
\plt 0 57 7 \dt 
\plt 0 59 4 \dt 
\plt 1 63 3 \dt 
\plt 0 64 1 \dt 
\plt 0 66 1 \dt 
\plt 1 66 27 \dt 
\plt 1 66 35 \dt 
\plt 0 66 92 \dt 
\plt 0 66 100 \dt 
\plt 0 67 3 \dt 
\plt 0 68 14 \dt 
\plt 0 68 18 \dt 
\plt 0 69 2 \dt 
\plt 0 70 9 \dt 
\plt 0 70 12 \dt 
\plt 0 72 7 \dt 
\plt 0 72 9 \dt 
\plt 0 74 7 \dt 
\plt 0 76 1 \dt 
\plt 0 76 6 \dt 
\plt 0 77 2 \dt 
\plt 0 78 1 \dt 
\plt 1 78 4 \dt 
\plt 1 82 3 \dt 
\plt 1 83 35 \dt 
\plt 1 83 44 \dt 
\plt 0 83 117 \dt 
\plt 0 83 126 \dt 
\plt 0 85 18 \dt 
\plt 0 85 22 \dt 
\plt 0 87 12 \dt 
\plt 0 87 15 \dt 
\plt 0 88 3 \dt 
\plt 0 89 1 \dt 
\plt 0 89 9 \dt 
\plt 0 89 11 \dt 
\plt 0 91 1 \dt 
\plt 0 91 7 \dt 
\plt 0 91 9 \dt 
\plt 0 92 2 \dt 
\plt 0 93 6 \dt 
\plt 0 95 5 \dt 
\plt 1 99 4 \dt 
\plt 0 100 2 \dt 
\plt 1 102 44 \dt 
\plt 1 102 54 \dt 
\plt 0 102 145 \dt 
\plt 0 102 155 \dt 
\plt 0 103 1 \dt 
\plt 0 103 4 \dt 
\plt 0 104 22 \dt 
\plt 0 104 27 \dt 
\plt 0 105 1 \dt 
\plt 0 105 3 \dt 
\plt 0 106 15 \dt 
\plt 0 106 18 \dt 
\plt 0 108 11 \dt 
\plt 0 110 9 \dt 
\plt 0 110 11 \dt 
\plt 0 111 3 \dt 
\plt 0 112 9 \dt 
\plt 0 117 2 \dt 
\plt 0 118 1 \dt 
\plt 1 118 5 \dt 
\plt 0 120 1 \dt 
\plt 1 122 4 \dt 
\plt 1 123 54 \dt 
\plt 1 123 65 \dt 
\plt 0 123 176 \dt 
\plt 0 123 187 \dt 
\plt 0 125 2 \dt 
\plt 0 125 27 \dt 
\plt 0 125 33 \dt 
\plt 0 127 18 \dt 
\plt 0 127 22 \dt 
\plt 0 130 3 \dt 
\plt 0 131 11 \dt 
\plt 0 131 13 \dt 
\plt 0 133 9 \dt 
\plt 0 133 11 \dt 
\plt 0 134 1 \dt 
\plt 0 136 1 \dt 
\plt 0 136 3 \dt 
\plt 0 139 6 \dt 
\plt 1 143 5 \dt 
\plt 0 144 2 \dt 
\plt 1 146 65 \dt 
\plt 1 146 77 \dt 
\plt 0 146 210 \dt 
\plt 0 147 5 \dt 
\plt 0 148 33 \dt 
\plt 0 148 39 \dt 
\plt 0 150 22 \dt 
\plt 0 150 26 \dt 
\plt 0 151 1 \dt 
\plt 0 153 1 \dt 
\plt 0 154 2 \dt 
\plt 0 154 13 \dt 
\plt 0 156 11 \dt 
\plt 0 156 13 \dt 
\plt 0 157 3 \dt 
\plt 0 158 11 \dt 
\plt 0 165 3 \dt 
\plt 1 166 6 \dt 
\plt 0 169 1 \dt 
\plt 1 170 5 \dt 
\plt 0 171 1 \dt 
\plt 1 171 77 \dt 
\plt 1 171 90 \dt 
\plt 0 173 39 \dt 
\plt 0 173 45 \dt 
\plt 0 175 2 \dt 
\plt 0 175 26 \dt 
\plt 0 175 30 \dt 
\plt 0 179 18 \dt 
\plt 0 181 13 \dt 
\plt 0 181 15 \dt 
\plt 0 183 11 \dt 
\plt 0 183 13 \dt 
\plt 0 187 2 \dt 
\plt 0 187 10 \dt 
\plt 0 188 1 \dt 
\plt 0 188 3 \dt 
\plt 0 189 9 \dt 
\plt 0 190 1 \dt 
\plt 0 191 7 \dt 
\plt 1 195 6 \dt 
\plt 0 196 3 \dt 
\plt 1 198 90 \dt 
\plt 1 198 104 \dt 
\plt 0 199 6 \dt 
\plt 0 200 45 \dt 
\plt 0 200 52 \dt 
\plt 0 202 30 \dt 
\plt 0 202 35 \dt 
\plt 0 204 26 \dt 
\plt 0 205 5 \dt 
\plt 0 206 18 \dt 
\plt 0 206 21 \dt 
\plt 0 208 1 \dt 
\plt 0 208 15 \dt 
\plt 1 210 1 \dt 
\plt 0 210 13 \dt 
\plt 0 210 15 \dt 
\plt 0 212 13 \dt 
\plt 0 214 10 \dt 
\plt 0 216 9 \dt 
\plt 1 16 1 \dt
\plt 1 16 4 \dt
\plt 1 16 7 \dt
\plt 1 16 10 \dt
\plt 1 16 13 \dt
\plt 1 16 16 \dt
\plt 1 16 19 \dt
\plt 1 16 22 \dt
\plt 1 16 25 \dt
\plt 1 16 28 \dt
\plt 1 16 31 \dt
\plt 1 16 34 \dt
\plt 1 16 37 \dt
\plt 1 16 40 \dt
\plt 1 16 43 \dt
\plt 1 16 46 \dt
\plt 1 16 49 \dt
\plt 1 16 52 \dt
\plt 1 16 55 \dt
\plt 1 16 58 \dt
\plt 1 16 61 \dt
\plt 1 16 64 \dt
\plt 1 16 67 \dt
\plt 1 16 70 \dt
\plt 1 16 73 \dt
\plt 1 16 76 \dt
\plt 1 16 79 \dt
\plt 1 16 82 \dt
\plt 1 16 85 \dt
\plt 1 16 88 \dt
\plt 1 16 91 \dt
\plt 1 16 94 \dt
\plt 1 16 97 \dt
\plt 1 16 100 \dt
\plt 1 16 103 \dt
\plt 1 16 106 \dt
\plt 1 16 109 \dt
\plt 1 16 112 \dt
\plt 1 16 115 \dt
\plt 1 16 118 \dt
\plt 1 16 121 \dt
\plt 1 16 124 \dt
\plt 1 16 127 \dt
\plt 1 16 130 \dt
\plt 1 16 133 \dt
\plt 1 16 136 \dt
\plt 1 16 139 \dt
\plt 1 16 142 \dt
\plt 1 16 145 \dt
\plt 1 16 148 \dt
\plt 1 16 151 \dt
\plt 1 16 154 \dt
\plt 1 16 157 \dt
\plt 1 16 160 \dt
\plt 1 16 163 \dt
\plt 1 16 166 \dt
\plt 1 16 169 \dt
\plt 1 16 172 \dt
\plt 1 16 175 \dt
\plt 1 16 178 \dt
\plt 1 16 181 \dt
\plt 1 16 184 \dt
\plt 1 16 187 \dt
\plt 1 16 190 \dt
\plt 1 16 193 \dt
\plt 1 16 196 \dt
\plt 1 16 199 \dt
\plt 1 16 202 \dt
\plt 1 16 205 \dt
\plt 1 16 208 \dt
\plt 1 16 211 \dt
\plt 1 16 214 \dt
\plt 1 16 217 \dt
\plt 1 16 220 \dt
\plt 2 36 1 \dt
\plt 2 36 6 \dt
\plt 2 36 11 \dt
\plt 2 36 16 \dt
\plt 2 36 21 \dt
\plt 2 36 26 \dt
\plt 2 36 31 \dt
\plt 2 36 36 \dt
\plt 2 36 41 \dt
\plt 2 36 46 \dt
\plt 2 36 51 \dt
\plt 2 36 56 \dt
\plt 2 36 61 \dt
\plt 2 36 66 \dt
\plt 2 36 71 \dt
\plt 2 36 76 \dt
\plt 2 36 81 \dt
\plt 2 36 86 \dt
\plt 2 36 91 \dt
\plt 2 36 96 \dt
\plt 2 36 101 \dt
\plt 2 36 106 \dt
\plt 2 36 111 \dt
\plt 2 36 116 \dt
\plt 2 36 121 \dt
\plt 2 36 126 \dt
\plt 2 36 131 \dt
\plt 2 36 136 \dt
\plt 2 36 141 \dt
\plt 2 36 146 \dt
\plt 2 36 151 \dt
\plt 2 36 156 \dt
\plt 2 36 161 \dt
\plt 2 36 166 \dt
\plt 2 36 171 \dt
\plt 2 36 176 \dt
\plt 2 36 181 \dt
\plt 2 36 186 \dt
\plt 2 36 191 \dt
\plt 2 36 196 \dt
\plt 2 36 201 \dt
\plt 2 36 206 \dt
\plt 2 36 211 \dt
\plt 2 36 216 \dt
\plt 2 36 221 \dt
\plt 2 64 2 \dt
\plt 2 64 9 \dt
\plt 2 64 16 \dt
\plt 2 64 23 \dt
\plt 2 64 30 \dt
\plt 2 64 37 \dt
\plt 2 64 44 \dt
\plt 2 64 51 \dt
\plt 2 64 58 \dt
\plt 2 64 65 \dt
\plt 2 64 72 \dt
\plt 2 64 79 \dt
\plt 2 64 86 \dt
\plt 2 64 93 \dt
\plt 2 64 100 \dt
\plt 2 64 107 \dt
\plt 2 64 114 \dt
\plt 2 64 121 \dt
\plt 2 64 128 \dt
\plt 2 64 135 \dt
\plt 2 64 142 \dt
\plt 2 64 149 \dt
\plt 2 64 156 \dt
\plt 2 64 163 \dt
\plt 2 64 170 \dt
\plt 2 64 177 \dt
\plt 2 64 184 \dt
\plt 2 64 191 \dt
\plt 2 64 198 \dt
\plt 2 64 205 \dt
\plt 2 64 212 \dt
\plt 2 64 219 \dt
\plt 2 100 3 \dt
\plt 2 100 12 \dt
\plt 2 100 21 \dt
\plt 2 100 30 \dt
\plt 2 100 39 \dt
\plt 2 100 48 \dt
\plt 2 100 57 \dt
\plt 2 100 66 \dt
\plt 2 100 75 \dt
\plt 2 100 84 \dt
\plt 2 100 93 \dt
\plt 2 100 102 \dt
\plt 2 100 111 \dt
\plt 2 100 120 \dt
\plt 2 100 129 \dt
\plt 2 100 138 \dt
\plt 2 100 147 \dt
\plt 2 100 156 \dt
\plt 2 100 165 \dt
\plt 2 100 174 \dt
\plt 2 100 183 \dt
\plt 2 100 192 \dt
\plt 2 100 201 \dt
\plt 2 100 210 \dt
\plt 2 100 219 \dt
\plt 3 144 3 \dt
\plt 3 144 14 \dt
\plt 3 144 25 \dt
\plt 3 144 36 \dt
\plt 3 144 47 \dt
\plt 3 144 58 \dt
\plt 3 144 69 \dt
\plt 3 144 80 \dt
\plt 3 144 91 \dt
\plt 3 144 102 \dt
\plt 3 144 113 \dt
\plt 3 144 124 \dt
\plt 3 144 135 \dt
\plt 3 144 146 \dt
\plt 3 144 157 \dt
\plt 3 144 168 \dt
\plt 3 144 179 \dt
\plt 3 144 190 \dt
\plt 3 144 201 \dt
\plt 3 144 212 \dt
\plt 3 196 4 \dt
\plt 3 196 17 \dt
\plt 3 196 30 \dt
\plt 3 196 43 \dt
\plt 3 196 56 \dt
\plt 3 196 69 \dt
\plt 3 196 82 \dt
\plt 3 196 95 \dt
\plt 3 196 108 \dt
\plt 3 196 121 \dt
\plt 3 196 134 \dt
\plt 3 196 147 \dt
\plt 3 196 160 \dt
\plt 3 196 173 \dt
\plt 3 196 186 \dt
\plt 3 196 199 \dt
\plt 3 196 212 \dt
}}

\vfil\eject

\def\xscale{1}   
\def\yscale{.8}   
\def\xmin{-5}   
\def\ymin{-5}    
\def\xmax{220}     
\def\ymax{220}     
\mywidth220pt     
\myheight220pt   
\penwidth1pt 
\penheight1pt 

\ \vskip-.5in\hskip1in\hbox{\vbox{      
\parindent0in                        
\putincaption 20 {\leftskip-1in\hsize4in\noindent{{\bf Figure 3}:
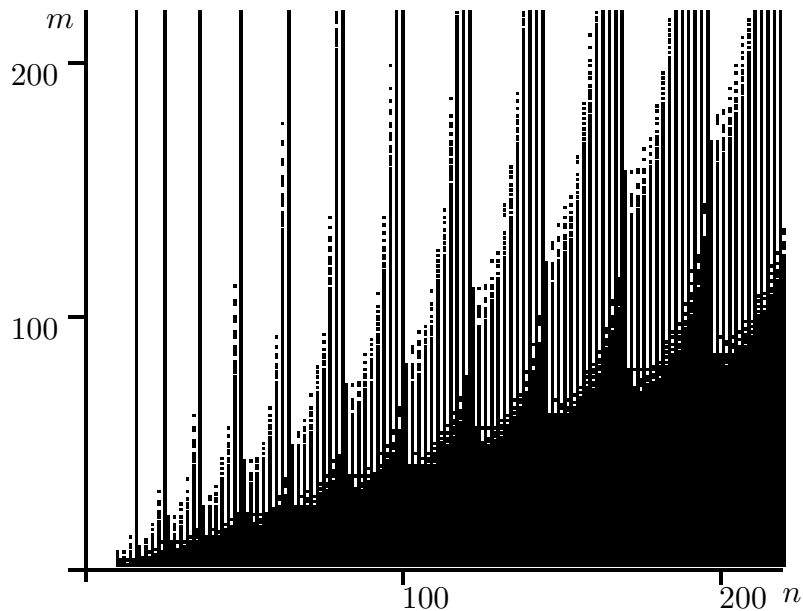
Graph of all $(n,m)$ for $n$ points of multiplicity $m$
in ${\bf C}P^2$ with $10\le n\le \xmax$ and $0\le m\le \ymax$, 
such that \ir{introthm},
using $d=\lfloor \sqrt{n}\rfloor$ and $r= \lfloor d\sqrt{n}\rfloor$,
implies $\alpha(n;m) \ge m\sqrt{n}$.}}
\putinaxes  $n$ $m$
\putinticx 5 {100} {100} 
\putinticx 5 {200} {200} 
\putinticy 5 -5 {100} {100} 
\putinticy 5 -5 {200} {200}




\plt 5 10 1 \dt 
\plt 2 11 1 \dt 
\plt 3 12 1 \dt 
\plt 0 12 6 \dt 
\plt 0 13 1 \dt 
\plt 0 13 3 \dt 
\plt 5 14 1 \dt 
\plt 1 14 8 \dt 
\plt 0 14 12 \dt 
\plt 2 15 1 \dt 
\plt 219 16 1 \dt 
\plt 7 17 1 \dt 
\plt 3 18 1 \dt 
\plt 7 19 1 \dt 
\plt 1 19 10 \dt 
\plt 3 20 1 \dt 
\plt 0 20 6 \dt 
\plt 12 21 1 \dt 
\plt 0 21 15 \dt 
\plt 0 21 17 \dt 
\plt 4 22 1 \dt 
\plt 0 22 7 \dt 
\plt 17 23 1 \dt 
\plt 2 23 20 \dt 
\plt 1 23 25 \dt 
\plt 0 23 30 \dt 
\plt 5 24 1 \dt 
\plt 0 24 10 \dt 
\plt 219 25 1 \dt 
\plt 19 26 1 \dt 
\plt 9 27 1 \dt 
\plt 12 28 1 \dt 
\plt 2 28 15 \dt 
\plt 0 28 20 \dt 
\plt 6 29 1 \dt 
\plt 0 29 10 \dt 
\plt 15 30 1 \dt 
\plt 0 30 18 \dt 
\plt 0 30 20 \dt 
\plt 0 30 23 \dt 
\plt 0 30 25 \dt 
\plt 7 31 1 \dt 
\plt 0 31 10 \dt 
\plt 21 32 1 \dt 
\plt 1 32 24 \dt 
\plt 0 32 27 \dt 
\plt 0 32 30 \dt 
\plt 0 32 35 \dt 
\plt 9 33 1 \dt 
\plt 0 33 12 \dt 
\plt 39 34 1 \dt 
\plt 3 34 42 \dt 
\plt 2 34 48 \dt 
\plt 1 34 54 \dt 
\plt 0 34 60 \dt 
\plt 9 35 1 \dt 
\plt 3 35 12 \dt 
\plt 219 36 1 \dt 
\plt 23 37 1 \dt 
\plt 11 38 1 \dt 
\plt 23 39 1 \dt 
\plt 1 39 27 \dt 
\plt 11 40 1 \dt 
\plt 0 40 15 \dt 
\plt 23 41 1 \dt 
\plt 1 41 26 \dt 
\plt 0 41 29 \dt 
\plt 0 41 32 \dt 
\plt 11 42 1 \dt 
\plt 0 42 14 \dt 
\plt 34 43 1 \dt 
\plt 0 43 37 \dt 
\plt 0 43 39 \dt 
\plt 0 43 41 \dt 
\plt 0 43 43 \dt 
\plt 14 44 1 \dt 
\plt 0 44 17 \dt 
\plt 39 45 1 \dt 
\plt 1 45 42 \dt 
\plt 1 45 45 \dt 
\plt 1 45 48 \dt 
\plt 0 45 52 \dt 
\plt 0 45 55 \dt 
\plt 15 46 1 \dt 
\plt 0 46 19 \dt 
\plt 74 47 1 \dt 
\plt 4 47 77 \dt 
\plt 3 47 84 \dt 
\plt 2 47 91 \dt 
\plt 2 47 97 \dt 
\plt 1 47 104 \dt 
\plt 0 47 111 \dt 
\plt 20 48 1 \dt 
\plt 219 49 1 \dt 
\plt 41 50 1 \dt 
\plt 20 51 1 \dt 
\plt 32 52 1 \dt 
\plt 2 52 35 \dt 
\plt 0 52 42 \dt 
\plt 16 53 1 \dt 
\plt 0 53 21 \dt 
\plt 36 54 1 \dt 
\plt 1 54 39 \dt 
\plt 1 54 42 \dt 
\plt 15 55 1 \dt 
\plt 1 55 18 \dt 
\plt 0 55 21 \dt 
\plt 42 56 1 \dt 
\plt 0 56 45 \dt 
\plt 0 56 47 \dt 
\plt 0 56 49 \dt 
\plt 20 57 1 \dt 
\plt 50 58 1 \dt 
\plt 0 58 53 \dt 
\plt 1 58 55 \dt 
\plt 0 58 58 \dt 
\plt 0 58 60 \dt 
\plt 0 58 63 \dt 
\plt 22 59 1 \dt 
\plt 72 60 1 \dt 
\plt 2 60 75 \dt 
\plt 1 60 79 \dt 
\plt 1 60 83 \dt 
\plt 0 60 87 \dt 
\plt 0 60 91 \dt 
\plt 23 61 1 \dt 
\plt 1 61 27 \dt 
\plt 132 62 1 \dt 
\plt 5 62 135 \dt 
\plt 4 62 143 \dt 
\plt 3 62 151 \dt 
\plt 2 62 159 \dt 
\plt 1 62 167 \dt 
\plt 0 62 175 \dt 
\plt 27 63 1 \dt 
\plt 3 63 32 \dt 
\plt 219 64 1 \dt 
\plt 47 65 1 \dt 
\plt 23 66 1 \dt 
\plt 47 67 1 \dt 
\plt 1 67 52 \dt 
\plt 23 68 1 \dt 
\plt 0 68 28 \dt 
\plt 47 69 1 \dt 
\plt 2 69 50 \dt 
\plt 1 69 54 \dt 
\plt 0 69 58 \dt 
\plt 23 70 1 \dt 
\plt 1 70 26 \dt 
\plt 0 70 30 \dt 
\plt 52 71 1 \dt 
\plt 1 71 55 \dt 
\plt 0 71 58 \dt 
\plt 1 71 60 \dt 
\plt 0 71 63 \dt 
\plt 0 71 68 \dt 
\plt 23 72 1 \dt 
\plt 0 72 26 \dt 
\plt 0 72 28 \dt 
\plt 68 73 1 \dt 
\plt 0 73 71 \dt 
\plt 0 73 73 \dt 
\plt 0 73 75 \dt 
\plt 0 73 77 \dt 
\plt 0 73 79 \dt 
\plt 28 74 1 \dt 
\plt 0 74 31 \dt 
\plt 0 74 33 \dt 
\plt 77 75 1 \dt 
\plt 1 75 80 \dt 
\plt 1 75 83 \dt 
\plt 0 75 86 \dt 
\plt 0 75 89 \dt 
\plt 0 75 92 \dt 
\plt 29 76 1 \dt 
\plt 1 76 32 \dt 
\plt 0 76 36 \dt 
\plt 109 77 1 \dt 
\plt 2 77 112 \dt 
\plt 2 77 116 \dt 
\plt 2 77 120 \dt 
\plt 1 77 125 \dt 
\plt 1 77 129 \dt 
\plt 0 77 134 \dt 
\plt 0 77 138 \dt 
\plt 33 78 1 \dt 
\plt 2 78 36 \dt 
\plt 0 78 42 \dt 
\plt 203 79 1 \dt 
\plt 6 79 206 \dt 
\plt 5 79 215 \dt 
\plt 43 80 1 \dt 
\plt 219 81 1 \dt 
\plt 71 82 1 \dt 
\plt 35 83 1 \dt 
\plt 59 84 1 \dt 
\plt 3 84 63 \dt 
\plt 0 84 72 \dt 
\plt 30 85 1 \dt 
\plt 0 85 36 \dt 
\plt 64 86 1 \dt 
\plt 1 86 68 \dt 
\plt 1 86 72 \dt 
\plt 29 87 1 \dt 
\plt 1 87 32 \dt 
\plt 0 87 36 \dt 
\plt 69 88 1 \dt 
\plt 1 88 72 \dt 
\plt 1 88 75 \dt 
\plt 0 88 78 \dt 
\plt 0 88 81 \dt 
\plt 0 88 84 \dt 
\plt 33 89 1 \dt 
\plt 0 89 36 \dt 
\plt 0 89 39 \dt 
\plt 72 90 1 \dt 
\plt 7 90 75 \dt 
\plt 0 90 84 \dt 
\plt 0 90 86 \dt 
\plt 0 90 88 \dt 
\plt 0 90 90 \dt 
\plt 36 91 1 \dt 
\plt 86 92 1 \dt 
\plt 3 92 89 \dt 
\plt 0 92 94 \dt 
\plt 0 92 96 \dt 
\plt 1 92 98 \dt 
\plt 0 92 101 \dt 
\plt 0 92 103 \dt 
\plt 0 92 108 \dt 
\plt 37 93 1 \dt 
\plt 0 93 40 \dt 
\plt 0 93 45 \dt 
\plt 116 94 1 \dt 
\plt 1 94 119 \dt 
\plt 1 94 122 \dt 
\plt 1 94 125 \dt 
\plt 0 94 129 \dt 
\plt 0 94 132 \dt 
\plt 0 94 135 \dt 
\plt 0 94 138 \dt 
\plt 41 95 1 \dt 
\plt 1 95 44 \dt 
\plt 0 95 48 \dt 
\plt 156 96 1 \dt 
\plt 3 96 159 \dt 
\plt 2 96 164 \dt 
\plt 2 96 169 \dt 
\plt 1 96 174 \dt 
\plt 1 96 179 \dt 
\plt 0 96 184 \dt 
\plt 0 96 189 \dt 
\plt 0 96 198 \dt 
\plt 44 97 1 \dt 
\plt 2 97 47 \dt 
\plt 1 97 53 \dt 
\plt 219 98 1 \dt 
\plt 53 99 1 \dt 
\plt 3 99 60 \dt 
\plt 219 100 1 \dt 
\plt 79 101 1 \dt 
\plt 39 102 1 \dt 
\plt 72 103 1 \dt 
\plt 5 103 75 \dt 
\plt 2 103 85 \dt 
\plt 39 104 1 \dt 
\plt 0 104 45 \dt 
\plt 75 105 1 \dt 
\plt 2 105 78 \dt 
\plt 2 105 83 \dt 
\plt 1 105 88 \dt 
\plt 0 105 93 \dt 
\plt 39 106 1 \dt 
\plt 0 106 43 \dt 
\plt 86 107 1 \dt 
\plt 1 107 89 \dt 
\plt 1 107 92 \dt 
\plt 0 107 95 \dt 
\plt 39 108 1 \dt 
\plt 1 108 42 \dt 
\plt 0 108 45 \dt 
\plt 89 109 1 \dt 
\plt 3 109 92 \dt 
\plt 0 109 97 \dt 
\plt 0 109 99 \dt 
\plt 0 109 102 \dt 
\plt 0 109 104 \dt 
\plt 0 109 109 \dt 
\plt 39 110 1 \dt 
\plt 3 110 42 \dt 
\plt 0 110 47 \dt 
\plt 0 110 49 \dt 
\plt 112 111 1 \dt 
\plt 0 111 115 \dt 
\plt 0 111 117 \dt 
\plt 0 111 119 \dt 
\plt 0 111 121 \dt 
\plt 0 111 123 \dt 
\plt 0 111 125 \dt 
\plt 48 112 1 \dt 
\plt 0 112 51 \dt 
\plt 0 112 53 \dt 
\plt 122 113 1 \dt 
\plt 1 113 125 \dt 
\plt 0 113 128 \dt 
\plt 1 113 130 \dt 
\plt 0 113 133 \dt 
\plt 0 113 136 \dt 
\plt 0 113 138 \dt 
\plt 0 113 141 \dt 
\plt 47 114 1 \dt 
\plt 1 114 50 \dt 
\plt 0 114 53 \dt 
\plt 0 114 56 \dt 
\plt 150 115 1 \dt 
\plt 5 115 153 \dt 
\plt 1 115 160 \dt 
\plt 2 115 163 \dt 
\plt 1 115 167 \dt 
\plt 0 115 171 \dt 
\plt 1 115 174 \dt 
\plt 0 115 178 \dt 
\plt 0 115 185 \dt 
\plt 50 116 1 \dt 
\plt 2 116 53 \dt 
\plt 1 116 57 \dt 
\plt 0 116 61 \dt 
\plt 216 117 1 \dt 
\plt 1 117 219 \dt 
\plt 56 118 1 \dt 
\plt 3 118 59 \dt 
\plt 1 118 66 \dt 
\plt 219 119 1 \dt 
\plt 74 120 1 \dt 
\plt 219 121 1 \dt 
\plt 109 122 1 \dt 
\plt 54 123 1 \dt 
\plt 94 124 1 \dt 
\plt 4 124 99 \dt 
\plt 0 124 110 \dt 
\plt 48 125 1 \dt 
\plt 0 125 55 \dt 
\plt 96 126 1 \dt 
\plt 3 126 99 \dt 
\plt 1 126 105 \dt 
\plt 1 126 110 \dt 
\plt 47 127 1 \dt 
\plt 1 127 50 \dt 
\plt 0 127 55 \dt 
\plt 104 128 1 \dt 
\plt 1 128 107 \dt 
\plt 1 128 110 \dt 
\plt 0 128 114 \dt 
\plt 0 128 121 \dt 
\plt 49 129 1 \dt 
\plt 1 129 52 \dt 
\plt 0 129 55 \dt 
\plt 113 130 1 \dt 
\plt 0 130 116 \dt 
\plt 0 130 119 \dt 
\plt 0 130 121 \dt 
\plt 0 130 124 \dt 
\plt 50 131 1 \dt 
\plt 0 131 53 \dt 
\plt 1 131 55 \dt 
\plt 0 131 58 \dt 
\plt 121 132 1 \dt 
\plt 0 132 124 \dt 
\plt 0 132 126 \dt 
\plt 0 132 128 \dt 
\plt 0 132 130 \dt 
\plt 0 132 132 \dt 
\plt 0 132 139 \dt 
\plt 0 132 141 \dt 
\plt 0 132 143 \dt 
\plt 55 133 1 \dt 
\plt 0 133 58 \dt 
\plt 0 133 60 \dt 
\plt 137 134 1 \dt 
\plt 0 134 140 \dt 
\plt 3 134 142 \dt 
\plt 0 134 147 \dt 
\plt 0 134 149 \dt 
\plt 0 134 151 \dt 
\plt 0 134 154 \dt 
\plt 0 134 156 \dt 
\plt 0 134 158 \dt 
\plt 58 135 1 \dt 
\plt 0 135 61 \dt 
\plt 0 135 63 \dt 
\plt 0 135 66 \dt 
\plt 166 136 1 \dt 
\plt 1 136 169 \dt 
\plt 1 136 172 \dt 
\plt 1 136 175 \dt 
\plt 0 136 178 \dt 
\plt 0 136 181 \dt 
\plt 0 136 184 \dt 
\plt 0 136 187 \dt 
\plt 62 137 1 \dt 
\plt 1 137 65 \dt 
\plt 0 137 68 \dt 
\plt 211 138 1 \dt 
\plt 2 138 214 \dt 
\plt 2 138 218 \dt 
\plt 68 139 1 \dt 
\plt 1 139 72 \dt 
\plt 0 139 77 \dt 
\plt 219 140 1 \dt 
\plt 76 141 1 \dt 
\plt 2 141 80 \dt 
\plt 0 141 88 \dt 
\plt 219 142 1 \dt 
\plt 87 143 1 \dt 
\plt 3 143 96 \dt 
\plt 219 144 1 \dt 
\plt 119 145 1 \dt 
\plt 59 146 1 \dt 
\plt 111 147 1 \dt 
\plt 6 147 114 \dt 
\plt 2 147 126 \dt 
\plt 59 148 1 \dt 
\plt 0 148 66 \dt 
\plt 120 149 1 \dt 
\plt 2 149 123 \dt 
\plt 1 149 129 \dt 
\plt 0 149 135 \dt 
\plt 59 150 1 \dt 
\plt 1 150 63 \dt 
\plt 123 151 1 \dt 
\plt 2 151 126 \dt 
\plt 1 151 130 \dt 
\plt 1 151 134 \dt 
\plt 0 151 138 \dt 
\plt 0 151 142 \dt 
\plt 59 152 1 \dt 
\plt 1 152 62 \dt 
\plt 0 152 66 \dt 
\plt 128 153 1 \dt 
\plt 1 153 131 \dt 
\plt 1 153 134 \dt 
\plt 1 153 137 \dt 
\plt 0 153 140 \dt 
\plt 0 153 143 \dt 
\plt 0 153 146 \dt 
\plt 62 154 1 \dt 
\plt 0 154 65 \dt 
\plt 0 154 68 \dt 
\plt 143 155 1 \dt 
\plt 0 155 146 \dt 
\plt 0 155 148 \dt 
\plt 1 155 150 \dt 
\plt 0 155 153 \dt 
\plt 0 155 155 \dt 
\plt 0 155 160 \dt 
\plt 0 155 162 \dt 
\plt 66 156 1 \dt 
\plt 0 156 69 \dt 
\plt 0 156 71 \dt 
\plt 166 157 1 \dt 
\plt 0 157 169 \dt 
\plt 0 157 171 \dt 
\plt 0 157 173 \dt 
\plt 0 157 175 \dt 
\plt 0 157 177 \dt 
\plt 0 157 179 \dt 
\plt 0 157 181 \dt 
\plt 0 157 183 \dt 
\plt 70 158 1 \dt 
\plt 0 158 73 \dt 
\plt 0 158 75 \dt 
\plt 0 158 77 \dt 
\plt 177 159 1 \dt 
\plt 3 159 180 \dt 
\plt 3 159 185 \dt 
\plt 0 159 190 \dt 
\plt 0 159 193 \dt 
\plt 0 159 195 \dt 
\plt 0 159 198 \dt 
\plt 0 159 200 \dt 
\plt 0 159 205 \dt 
\plt 0 159 210 \dt 
\plt 69 160 1 \dt 
\plt 1 160 72 \dt 
\plt 0 160 75 \dt 
\plt 1 160 77 \dt 
\plt 0 160 80 \dt 
\plt 213 161 1 \dt 
\plt 1 161 216 \dt 
\plt 1 161 219 \dt 
\plt 75 162 1 \dt 
\plt 1 162 78 \dt 
\plt 1 162 81 \dt 
\plt 0 162 85 \dt 
\plt 219 163 1 \dt 
\plt 85 164 1 \dt 
\plt 2 164 88 \dt 
\plt 1 164 93 \dt 
\plt 219 165 1 \dt 
\plt 92 166 1 \dt 
\plt 4 166 95 \dt 
\plt 1 166 104 \dt 
\plt 219 167 1 \dt 
\plt 101 168 1 \dt 
\plt 10 168 104 \dt 
\plt 219 169 1 \dt 
\plt 155 170 1 \dt 
\plt 77 171 1 \dt 
\plt 138 172 1 \dt 
\plt 4 172 143 \dt 
\plt 0 172 156 \dt 
\plt 70 173 1 \dt 
\plt 0 173 78 \dt 
\plt 140 174 1 \dt 
\plt 3 174 143 \dt 
\plt 1 174 150 \dt 
\plt 1 174 156 \dt 
\plt 68 175 1 \dt 
\plt 2 175 72 \dt 
\plt 0 175 78 \dt 
\plt 145 176 1 \dt 
\plt 2 176 148 \dt 
\plt 1 176 152 \dt 
\plt 1 176 156 \dt 
\plt 0 176 161 \dt 
\plt 0 176 165 \dt 
\plt 71 177 1 \dt 
\plt 1 177 74 \dt 
\plt 0 177 78 \dt 
\plt 157 178 1 \dt 
\plt 0 178 160 \dt 
\plt 0 178 163 \dt 
\plt 0 178 166 \dt 
\plt 0 178 169 \dt 
\plt 72 179 1 \dt 
\plt 1 179 75 \dt 
\plt 1 179 78 \dt 
\plt 159 180 1 \dt 
\plt 3 180 162 \dt 
\plt 3 180 167 \dt 
\plt 0 180 172 \dt 
\plt 0 180 175 \dt 
\plt 0 180 177 \dt 
\plt 0 180 180 \dt 
\plt 0 180 182 \dt 
\plt 73 181 1 \dt 
\plt 3 181 76 \dt 
\plt 0 181 81 \dt 
\plt 182 182 1 \dt 
\plt 0 182 185 \dt 
\plt 0 182 187 \dt 
\plt 0 182 189 \dt 
\plt 0 182 191 \dt 
\plt 0 182 193 \dt 
\plt 0 182 195 \dt 
\plt 78 183 1 \dt 
\plt 0 183 81 \dt 
\plt 0 183 83 \dt 
\plt 0 183 85 \dt 
\plt 0 183 87 \dt 
\plt 0 183 89 \dt 
\plt 200 184 1 \dt 
\plt 0 184 203 \dt 
\plt 0 184 205 \dt 
\plt 1 184 207 \dt 
\plt 0 184 210 \dt 
\plt 0 184 212 \dt 
\plt 0 184 214 \dt 
\plt 0 184 216 \dt 
\plt 83 185 1 \dt 
\plt 0 185 86 \dt 
\plt 0 185 88 \dt 
\plt 0 185 91 \dt 
\plt 0 185 93 \dt 
\plt 219 186 1 \dt 
\plt 84 187 1 \dt 
\plt 1 187 87 \dt 
\plt 1 187 90 \dt 
\plt 0 187 93 \dt 
\plt 0 187 96 \dt 
\plt 219 188 1 \dt 
\plt 93 189 1 \dt 
\plt 1 189 96 \dt 
\plt 0 189 100 \dt 
\plt 0 189 104 \dt 
\plt 219 190 1 \dt 
\plt 98 191 1 \dt 
\plt 3 191 101 \dt 
\plt 2 191 106 \dt 
\plt 0 191 112 \dt 
\plt 0 191 117 \dt 
\plt 219 192 1 \dt 
\plt 109 193 1 \dt 
\plt 5 193 112 \dt 
\plt 2 193 121 \dt 
\plt 219 194 1 \dt 
\plt 129 195 1 \dt 
\plt 3 195 140 \dt 
\plt 219 196 1 \dt 
\plt 167 197 1 \dt 
\plt 83 198 1 \dt 
\plt 158 199 1 \dt 
\plt 7 199 161 \dt 
\plt 3 199 175 \dt 
\plt 83 200 1 \dt 
\plt 0 200 91 \dt 
\plt 162 201 1 \dt 
\plt 4 201 165 \dt 
\plt 2 201 172 \dt 
\plt 1 201 179 \dt 
\plt 78 202 1 \dt 
\plt 3 202 81 \dt 
\plt 1 202 88 \dt 
\plt 168 203 1 \dt 
\plt 2 203 171 \dt 
\plt 2 203 175 \dt 
\plt 1 203 180 \dt 
\plt 0 203 185 \dt 
\plt 0 203 189 \dt 
\plt 83 204 1 \dt 
\plt 1 204 87 \dt 
\plt 0 204 91 \dt 
\plt 178 205 1 \dt 
\plt 1 205 181 \dt 
\plt 1 205 184 \dt 
\plt 0 205 188 \dt 
\plt 0 205 191 \dt 
\plt 0 205 195 \dt 
\plt 0 205 198 \dt 
\plt 86 206 1 \dt 
\plt 0 206 90 \dt 
\plt 0 206 93 \dt 
\plt 184 207 1 \dt 
\plt 3 207 187 \dt 
\plt 1 207 192 \dt 
\plt 1 207 195 \dt 
\plt 0 207 198 \dt 
\plt 0 207 201 \dt 
\plt 0 207 203 \dt 
\plt 0 207 206 \dt 
\plt 0 207 209 \dt 
\plt 86 208 1 \dt 
\plt 3 208 89 \dt 
\plt 0 208 94 \dt 
\plt 0 208 97 \dt 
\plt 202 209 1 \dt 
\plt 0 209 205 \dt 
\plt 0 209 207 \dt 
\plt 1 209 209 \dt 
\plt 0 209 212 \dt 
\plt 0 209 214 \dt 
\plt 0 209 216 \dt 
\plt 90 210 1 \dt 
\plt 0 210 93 \dt 
\plt 0 210 95 \dt 
\plt 0 210 97 \dt 
\plt 219 211 1 \dt 
\plt 98 212 1 \dt 
\plt 0 212 101 \dt 
\plt 0 212 103 \dt 
\plt 0 212 105 \dt 
\plt 0 212 107 \dt 
\plt 219 213 1 \dt 
\plt 100 214 1 \dt 
\plt 0 214 103 \dt 
\plt 0 214 106 \dt 
\plt 0 214 108 \dt 
\plt 219 215 1 \dt 
\plt 106 216 1 \dt 
\plt 1 216 109 \dt 
\plt 0 216 113 \dt 
\plt 0 216 116 \dt 
\plt 0 216 119 \dt 
\plt 219 217 1 \dt 
\plt 112 218 1 \dt 
\plt 2 218 115 \dt 
\plt 0 218 120 \dt 
\plt 0 218 124 \dt 
\plt 219 219 1 \dt 
\plt 122 220 1 \dt 
\plt 2 220 126 \dt 
\plt 1 220 132 \dt 
}}

\def\xscale{1}   
\def\yscale{.8}   
\def\xmin{-5}   
\def\ymin{-5}    
\def\xmax{220}     
\def\ymax{220}     
\mywidth220pt     
\myheight220pt   
\penwidth1pt 
\penheight1pt 

\vskip3.5in\hskip1in\hbox{\vbox{      
\parindent0in                        
\putincaption 20 {\leftskip-1in\hsize4in\noindent{{\bf Figure 4}:
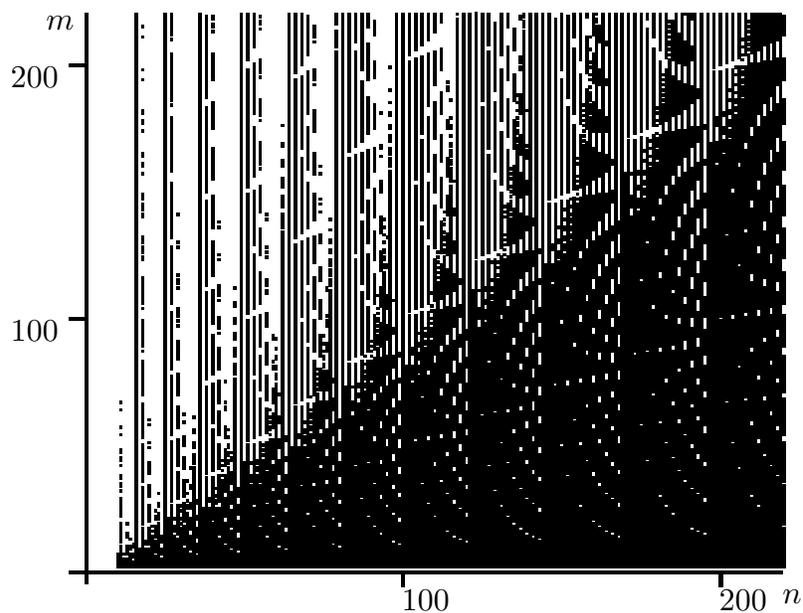
Graph of all $(n,m)$ for $n$ points of multiplicity $m$
in ${\bf C}P^2$ with $10\le n\le \xmax$ and $0\le m\le \ymax$, 
such that \ir{introthm},
using $d=\lfloor \sqrt{n}\rfloor$ and $r= \lfloor (n+d^2)/2\rfloor$,
implies $\alpha(n;m) \ge m\sqrt{n}$.}}
\putinaxes  $n$ $m$
\putinticx 5 {100} {100} 
\putinticx 5 {200} {200} 
\putinticy 5 -5 {100} {100} 
\putinticy 5 -5 {200} {200}




\plt 5 10 1 \dt 
\plt 8 11 1 \dt 
\plt 7 11 11 \dt 
\plt 7 11 21 \dt 
\plt 0 11 31 \dt 
\plt 1 11 33 \dt 
\plt 1 11 36 \dt 
\plt 0 11 41 \dt 
\plt 1 11 43 \dt 
\plt 1 11 46 \dt 
\plt 0 11 53 \dt 
\plt 0 11 56 \dt 
\plt 0 11 63 \dt 
\plt 0 11 66 \dt 
\plt 3 12 1 \dt 
\plt 0 12 6 \dt 
\plt 3 13 1 \dt 
\plt 3 13 6 \dt 
\plt 1 13 13 \dt 
\plt 0 13 18 \dt 
\plt 5 14 1 \dt 
\plt 1 14 8 \dt 
\plt 0 14 12 \dt 
\plt 2 15 1 \dt 
\plt 1 15 5 \dt 
\plt 1 15 9 \dt 
\plt 219 16 1 \dt 
\plt 7 17 1 \dt 
\plt 15 18 1 \dt 
\plt 14 18 18 \dt 
\plt 14 18 35 \dt 
\plt 13 18 52 \dt 
\plt 13 18 69 \dt 
\plt 0 18 86 \dt 
\plt 10 18 88 \dt 
\plt 0 18 103 \dt 
\plt 10 18 105 \dt 
\plt 1 18 122 \dt 
\plt 2 18 125 \dt 
\plt 2 18 129 \dt 
\plt 1 18 139 \dt 
\plt 2 18 142 \dt 
\plt 2 18 146 \dt 
\plt 0 18 156 \dt 
\plt 1 18 159 \dt 
\plt 1 18 163 \dt 
\plt 0 18 173 \dt 
\plt 1 18 176 \dt 
\plt 1 18 180 \dt 
\plt 0 18 193 \dt 
\plt 0 18 197 \dt 
\plt 0 18 210 \dt 
\plt 0 18 214 \dt 
\plt 7 19 1 \dt 
\plt 1 19 10 \dt 
\plt 7 20 1 \dt 
\plt 6 20 10 \dt 
\plt 6 20 19 \dt 
\plt 4 20 29 \dt 
\plt 4 20 38 \dt 
\plt 0 20 48 \dt 
\plt 0 20 50 \dt 
\plt 0 20 57 \dt 
\plt 0 20 59 \dt 
\plt 12 21 1 \dt 
\plt 0 21 15 \dt 
\plt 0 21 17 \dt 
\plt 4 22 1 \dt 
\plt 4 22 7 \dt 
\plt 4 22 13 \dt 
\plt 1 22 20 \dt 
\plt 0 22 23 \dt 
\plt 17 23 1 \dt 
\plt 2 23 20 \dt 
\plt 1 23 25 \dt 
\plt 0 23 30 \dt 
\plt 3 24 1 \dt 
\plt 2 24 6 \dt 
\plt 2 24 11 \dt 
\plt 219 25 1 \dt 
\plt 19 26 1 \dt 
\plt 24 27 1 \dt 
\plt 23 27 27 \dt 
\plt 23 27 53 \dt 
\plt 22 27 79 \dt 
\plt 22 27 105 \dt 
\plt 21 27 131 \dt 
\plt 21 27 157 \dt 
\plt 0 27 183 \dt 
\plt 18 27 185 \dt 
\plt 0 27 209 \dt 
\plt 9 27 211 \dt 
\plt 12 28 1 \dt 
\plt 2 28 15 \dt 
\plt 0 28 20 \dt 
\plt 11 29 1 \dt 
\plt 11 29 14 \dt 
\plt 10 29 28 \dt 
\plt 10 29 41 \dt 
\plt 8 29 56 \dt 
\plt 8 29 69 \dt 
\plt 7 29 83 \dt 
\plt 0 29 96 \dt 
\plt 0 29 98 \dt 
\plt 3 29 100 \dt 
\plt 0 29 111 \dt 
\plt 1 29 113 \dt 
\plt 0 29 116 \dt 
\plt 0 29 127 \dt 
\plt 0 29 129 \dt 
\plt 0 29 140 \dt 
\plt 15 30 1 \dt 
\plt 0 30 18 \dt 
\plt 0 30 20 \dt 
\plt 0 30 23 \dt 
\plt 0 30 25 \dt 
\plt 7 31 1 \dt 
\plt 7 31 10 \dt 
\plt 7 31 19 \dt 
\plt 6 31 29 \dt 
\plt 4 31 40 \dt 
\plt 3 31 49 \dt 
\plt 0 31 58 \dt 
\plt 0 31 61 \dt 
\plt 21 32 1 \dt 
\plt 1 32 24 \dt 
\plt 0 32 27 \dt 
\plt 0 32 30 \dt 
\plt 0 32 35 \dt 
\plt 5 33 1 \dt 
\plt 5 33 8 \dt 
\plt 4 33 16 \dt 
\plt 2 33 24 \dt 
\plt 1 33 32 \dt 
\plt 0 33 40 \dt 
\plt 39 34 1 \dt 
\plt 3 34 42 \dt 
\plt 2 34 48 \dt 
\plt 1 34 54 \dt 
\plt 0 34 60 \dt 
\plt 4 35 1 \dt 
\plt 3 35 7 \dt 
\plt 3 35 13 \dt 
\plt 1 35 20 \dt 
\plt 1 35 26 \dt 
\plt 219 36 1 \dt 
\plt 23 37 1 \dt 
\plt 35 38 1 \dt 
\plt 34 38 38 \dt 
\plt 34 38 75 \dt 
\plt 33 38 112 \dt 
\plt 33 38 149 \dt 
\plt 32 38 186 \dt 
\plt 23 39 1 \dt 
\plt 1 39 27 \dt 
\plt 17 40 1 \dt 
\plt 16 40 20 \dt 
\plt 16 40 39 \dt 
\plt 15 40 58 \dt 
\plt 15 40 77 \dt 
\plt 13 40 97 \dt 
\plt 13 40 116 \dt 
\plt 12 40 135 \dt 
\plt 12 40 154 \dt 
\plt 10 40 174 \dt 
\plt 10 40 193 \dt 
\plt 0 40 212 \dt 
\plt 1 40 214 \dt 
\plt 1 40 217 \dt 
\plt 0 40 220 \dt 
\plt 23 41 1 \dt 
\plt 1 41 26 \dt 
\plt 0 41 29 \dt 
\plt 0 41 32 \dt 
\plt 11 42 1 \dt 
\plt 10 42 14 \dt 
\plt 10 42 27 \dt 
\plt 9 42 40 \dt 
\plt 9 42 53 \dt 
\plt 8 42 66 \dt 
\plt 8 42 79 \dt 
\plt 0 42 93 \dt 
\plt 4 42 95 \dt 
\plt 0 42 106 \dt 
\plt 4 42 108 \dt 
\plt 0 42 122 \dt 
\plt 0 42 124 \dt 
\plt 0 42 135 \dt 
\plt 0 42 137 \dt 
\plt 34 43 1 \dt 
\plt 0 43 37 \dt 
\plt 0 43 39 \dt 
\plt 0 43 41 \dt 
\plt 0 43 43 \dt 
\plt 8 44 1 \dt 
\plt 7 44 11 \dt 
\plt 7 44 21 \dt 
\plt 6 44 31 \dt 
\plt 6 44 41 \dt 
\plt 4 44 52 \dt 
\plt 1 44 63 \dt 
\plt 0 44 66 \dt 
\plt 0 44 74 \dt 
\plt 39 45 1 \dt 
\plt 1 45 42 \dt 
\plt 1 45 45 \dt 
\plt 1 45 48 \dt 
\plt 0 45 52 \dt 
\plt 0 45 55 \dt 
\plt 6 46 1 \dt 
\plt 6 46 9 \dt 
\plt 6 46 17 \dt 
\plt 4 46 26 \dt 
\plt 4 46 34 \dt 
\plt 1 46 43 \dt 
\plt 0 46 46 \dt 
\plt 74 47 1 \dt 
\plt 4 47 77 \dt 
\plt 3 47 84 \dt 
\plt 2 47 91 \dt 
\plt 2 47 97 \dt 
\plt 1 47 104 \dt 
\plt 0 47 111 \dt 
\plt 5 48 1 \dt 
\plt 4 48 8 \dt 
\plt 4 48 15 \dt 
\plt 3 48 22 \dt 
\plt 3 48 29 \dt 
\plt 219 49 1 \dt 
\plt 41 50 1 \dt 
\plt 48 51 1 \dt 
\plt 47 51 51 \dt 
\plt 47 51 101 \dt 
\plt 46 51 151 \dt 
\plt 19 51 201 \dt 
\plt 32 52 1 \dt 
\plt 2 52 35 \dt 
\plt 0 52 42 \dt 
\plt 23 53 1 \dt 
\plt 23 53 26 \dt 
\plt 22 53 52 \dt 
\plt 22 53 77 \dt 
\plt 21 53 103 \dt 
\plt 21 53 128 \dt 
\plt 19 53 155 \dt 
\plt 19 53 180 \dt 
\plt 14 53 206 \dt 
\plt 36 54 1 \dt 
\plt 1 54 39 \dt 
\plt 1 54 42 \dt 
\plt 15 55 1 \dt 
\plt 15 55 18 \dt 
\plt 15 55 35 \dt 
\plt 14 55 53 \dt 
\plt 13 55 71 \dt 
\plt 12 55 88 \dt 
\plt 12 55 105 \dt 
\plt 10 55 124 \dt 
\plt 10 55 141 \dt 
\plt 0 55 158 \dt 
\plt 8 55 160 \dt 
\plt 0 55 177 \dt 
\plt 6 55 179 \dt 
\plt 0 55 194 \dt 
\plt 5 55 196 \dt 
\plt 0 55 213 \dt 
\plt 1 55 215 \dt 
\plt 0 55 218 \dt 
\plt 42 56 1 \dt 
\plt 0 56 45 \dt 
\plt 0 56 47 \dt 
\plt 0 56 49 \dt 
\plt 11 57 1 \dt 
\plt 11 57 14 \dt 
\plt 11 57 27 \dt 
\plt 10 57 40 \dt 
\plt 9 57 54 \dt 
\plt 9 57 67 \dt 
\plt 8 57 81 \dt 
\plt 5 57 96 \dt 
\plt 1 57 109 \dt 
\plt 2 57 112 \dt 
\plt 0 57 123 \dt 
\plt 0 57 125 \dt 
\plt 0 57 127 \dt 
\plt 0 57 138 \dt 
\plt 0 57 140 \dt 
\plt 50 58 1 \dt 
\plt 0 58 53 \dt 
\plt 1 58 55 \dt 
\plt 0 58 58 \dt 
\plt 0 58 60 \dt 
\plt 0 58 63 \dt 
\plt 19 59 1 \dt 
\plt 8 59 22 \dt 
\plt 8 59 33 \dt 
\plt 7 59 44 \dt 
\plt 5 59 56 \dt 
\plt 5 59 66 \dt 
\plt 1 59 79 \dt 
\plt 0 59 82 \dt 
\plt 1 59 91 \dt 
\plt 72 60 1 \dt 
\plt 2 60 75 \dt 
\plt 1 60 79 \dt 
\plt 1 60 83 \dt 
\plt 0 60 87 \dt 
\plt 0 60 91 \dt 
\plt 7 61 1 \dt 
\plt 7 61 10 \dt 
\plt 6 61 19 \dt 
\plt 6 61 28 \dt 
\plt 5 61 38 \dt 
\plt 3 61 48 \dt 
\plt 2 61 58 \dt 
\plt 0 61 69 \dt 
\plt 132 62 1 \dt 
\plt 5 62 135 \dt 
\plt 4 62 143 \dt 
\plt 3 62 151 \dt 
\plt 2 62 159 \dt 
\plt 1 62 167 \dt 
\plt 0 62 175 \dt 
\plt 6 63 1 \dt 
\plt 5 63 9 \dt 
\plt 5 63 17 \dt 
\plt 4 63 25 \dt 
\plt 4 63 33 \dt 
\plt 1 63 43 \dt 
\plt 1 63 51 \dt 
\plt 219 64 1 \dt 
\plt 47 65 1 \dt 
\plt 63 66 1 \dt 
\plt 62 66 66 \dt 
\plt 62 66 131 \dt 
\plt 24 66 196 \dt 
\plt 47 67 1 \dt 
\plt 1 67 52 \dt 
\plt 31 68 1 \dt 
\plt 30 68 34 \dt 
\plt 30 68 67 \dt 
\plt 29 68 100 \dt 
\plt 29 68 133 \dt 
\plt 28 68 166 \dt 
\plt 21 68 199 \dt 
\plt 47 69 1 \dt 
\plt 2 69 50 \dt 
\plt 1 69 54 \dt 
\plt 0 69 58 \dt 
\plt 20 70 1 \dt 
\plt 20 70 23 \dt 
\plt 20 70 45 \dt 
\plt 19 70 68 \dt 
\plt 19 70 90 \dt 
\plt 17 70 113 \dt 
\plt 17 70 135 \dt 
\plt 16 70 158 \dt 
\plt 16 70 180 \dt 
\plt 15 70 203 \dt 
\plt 52 71 1 \dt 
\plt 1 71 55 \dt 
\plt 0 71 58 \dt 
\plt 1 71 60 \dt 
\plt 0 71 63 \dt 
\plt 0 71 68 \dt 
\plt 15 72 1 \dt 
\plt 14 72 18 \dt 
\plt 14 72 35 \dt 
\plt 13 72 52 \dt 
\plt 13 72 69 \dt 
\plt 12 72 86 \dt 
\plt 12 72 103 \dt 
\plt 10 72 121 \dt 
\plt 10 72 138 \dt 
\plt 8 72 156 \dt 
\plt 8 72 173 \dt 
\plt 0 72 191 \dt 
\plt 4 72 193 \dt 
\plt 0 72 208 \dt 
\plt 4 72 210 \dt 
\plt 68 73 1 \dt 
\plt 0 73 71 \dt 
\plt 0 73 73 \dt 
\plt 0 73 75 \dt 
\plt 0 73 77 \dt 
\plt 0 73 79 \dt 
\plt 25 74 1 \dt 
\plt 11 74 28 \dt 
\plt 11 74 42 \dt 
\plt 10 74 56 \dt 
\plt 8 74 71 \dt 
\plt 8 74 84 \dt 
\plt 7 74 99 \dt 
\plt 5 74 114 \dt 
\plt 3 74 129 \dt 
\plt 0 74 141 \dt 
\plt 0 74 144 \dt 
\plt 0 74 146 \dt 
\plt 0 74 156 \dt 
\plt 0 74 159 \dt 
\plt 77 75 1 \dt 
\plt 1 75 80 \dt 
\plt 1 75 83 \dt 
\plt 0 75 86 \dt 
\plt 0 75 89 \dt 
\plt 0 75 92 \dt 
\plt 21 76 1 \dt 
\plt 9 76 24 \dt 
\plt 8 76 36 \dt 
\plt 7 76 48 \dt 
\plt 6 76 60 \dt 
\plt 7 76 71 \dt 
\plt 3 76 86 \dt 
\plt 2 76 96 \dt 
\plt 0 76 100 \dt 
\plt 1 76 110 \dt 
\plt 109 77 1 \dt 
\plt 2 77 112 \dt 
\plt 2 77 116 \dt 
\plt 2 77 120 \dt 
\plt 1 77 125 \dt 
\plt 1 77 129 \dt 
\plt 0 77 134 \dt 
\plt 0 77 138 \dt 
\plt 8 78 1 \dt 
\plt 8 78 11 \dt 
\plt 7 78 21 \dt 
\plt 7 78 31 \dt 
\plt 6 78 42 \dt 
\plt 4 78 53 \dt 
\plt 4 78 63 \dt 
\plt 1 78 74 \dt 
\plt 203 79 1 \dt 
\plt 6 79 206 \dt 
\plt 5 79 215 \dt 
\plt 7 80 1 \dt 
\plt 6 80 10 \dt 
\plt 6 80 19 \dt 
\plt 5 80 28 \dt 
\plt 5 80 37 \dt 
\plt 3 80 47 \dt 
\plt 3 80 56 \dt 
\plt 219 81 1 \dt 
\plt 71 82 1 \dt 
\plt 80 83 1 \dt 
\plt 79 83 83 \dt 
\plt 55 83 165 \dt 
\plt 59 84 1 \dt 
\plt 3 84 63 \dt 
\plt 0 84 72 \dt 
\plt 39 85 1 \dt 
\plt 39 85 42 \dt 
\plt 38 85 84 \dt 
\plt 38 85 125 \dt 
\plt 37 85 167 \dt 
\plt 12 85 208 \dt 
\plt 64 86 1 \dt 
\plt 1 86 68 \dt 
\plt 1 86 72 \dt 
\plt 26 87 1 \dt 
\plt 25 87 29 \dt 
\plt 25 87 57 \dt 
\plt 24 87 85 \dt 
\plt 24 87 113 \dt 
\plt 23 87 141 \dt 
\plt 23 87 169 \dt 
\plt 22 87 197 \dt 
\plt 69 88 1 \dt 
\plt 1 88 72 \dt 
\plt 1 88 75 \dt 
\plt 0 88 78 \dt 
\plt 0 88 81 \dt 
\plt 0 88 84 \dt 
\plt 19 89 1 \dt 
\plt 19 89 22 \dt 
\plt 19 89 43 \dt 
\plt 18 89 64 \dt 
\plt 17 89 86 \dt 
\plt 17 89 107 \dt 
\plt 16 89 129 \dt 
\plt 14 89 151 \dt 
\plt 14 89 172 \dt 
\plt 14 89 193 \dt 
\plt 6 89 214 \dt 
\plt 72 90 1 \dt 
\plt 7 90 75 \dt 
\plt 0 90 84 \dt 
\plt 0 90 86 \dt 
\plt 0 90 88 \dt 
\plt 0 90 90 \dt 
\plt 15 91 1 \dt 
\plt 15 91 18 \dt 
\plt 15 91 35 \dt 
\plt 13 91 53 \dt 
\plt 13 91 70 \dt 
\plt 13 91 87 \dt 
\plt 11 91 105 \dt 
\plt 11 91 122 \dt 
\plt 11 91 139 \dt 
\plt 9 91 157 \dt 
\plt 9 91 174 \dt 
\plt 6 91 194 \dt 
\plt 5 91 211 \dt 
\plt 86 92 1 \dt 
\plt 3 92 89 \dt 
\plt 0 92 94 \dt 
\plt 0 92 96 \dt 
\plt 1 92 98 \dt 
\plt 0 92 101 \dt 
\plt 0 92 103 \dt 
\plt 0 92 108 \dt 
\plt 12 93 1 \dt 
\plt 12 93 15 \dt 
\plt 11 93 30 \dt 
\plt 11 93 44 \dt 
\plt 10 93 59 \dt 
\plt 10 93 73 \dt 
\plt 9 93 88 \dt 
\plt 7 93 104 \dt 
\plt 7 93 118 \dt 
\plt 6 93 133 \dt 
\plt 4 93 149 \dt 
\plt 0 93 163 \dt 
\plt 1 93 166 \dt 
\plt 0 93 180 \dt 
\plt 116 94 1 \dt 
\plt 1 94 119 \dt 
\plt 1 94 122 \dt 
\plt 1 94 125 \dt 
\plt 0 94 129 \dt 
\plt 0 94 132 \dt 
\plt 0 94 135 \dt 
\plt 0 94 138 \dt 
\plt 23 95 1 \dt 
\plt 10 95 26 \dt 
\plt 9 95 39 \dt 
\plt 8 95 52 \dt 
\plt 8 95 64 \dt 
\plt 7 95 77 \dt 
\plt 5 95 91 \dt 
\plt 5 95 103 \dt 
\plt 1 95 119 \dt 
\plt 0 95 131 \dt 
\plt 156 96 1 \dt 
\plt 3 96 159 \dt 
\plt 2 96 164 \dt 
\plt 2 96 169 \dt 
\plt 1 96 174 \dt 
\plt 1 96 179 \dt 
\plt 0 96 184 \dt 
\plt 0 96 189 \dt 
\plt 0 96 198 \dt 
\plt 9 97 1 \dt 
\plt 9 97 12 \dt 
\plt 8 97 23 \dt 
\plt 8 97 34 \dt 
\plt 7 97 46 \dt 
\plt 6 97 57 \dt 
\plt 6 97 68 \dt 
\plt 4 97 80 \dt 
\plt 2 97 93 \dt 
\plt 0 97 106 \dt 
\plt 219 98 1 \dt 
\plt 8 99 1 \dt 
\plt 7 99 11 \dt 
\plt 7 99 21 \dt 
\plt 6 99 31 \dt 
\plt 6 99 41 \dt 
\plt 4 99 52 \dt 
\plt 4 99 62 \dt 
\plt 1 99 74 \dt 
\plt 1 99 84 \dt 
\plt 219 100 1 \dt 
\plt 79 101 1 \dt 
\plt 99 102 1 \dt 
\plt 98 102 102 \dt 
\plt 17 102 203 \dt 
\plt 72 103 1 \dt 
\plt 5 103 75 \dt 
\plt 2 103 85 \dt 
\plt 49 104 1 \dt 
\plt 48 104 52 \dt 
\plt 48 104 103 \dt 
\plt 47 104 154 \dt 
\plt 15 104 205 \dt 
\plt 75 105 1 \dt 
\plt 2 105 78 \dt 
\plt 2 105 83 \dt 
\plt 1 105 88 \dt 
\plt 0 105 93 \dt 
\plt 32 106 1 \dt 
\plt 32 106 35 \dt 
\plt 32 106 69 \dt 
\plt 31 106 104 \dt 
\plt 31 106 138 \dt 
\plt 30 106 172 \dt 
\plt 13 106 207 \dt 
\plt 86 107 1 \dt 
\plt 1 107 89 \dt 
\plt 1 107 92 \dt 
\plt 0 107 95 \dt 
\plt 24 108 1 \dt 
\plt 23 108 27 \dt 
\plt 23 108 53 \dt 
\plt 22 108 79 \dt 
\plt 22 108 105 \dt 
\plt 21 108 131 \dt 
\plt 21 108 157 \dt 
\plt 20 108 183 \dt 
\plt 11 108 209 \dt 
\plt 89 109 1 \dt 
\plt 3 109 92 \dt 
\plt 0 109 97 \dt 
\plt 0 109 99 \dt 
\plt 0 109 102 \dt 
\plt 0 109 104 \dt 
\plt 0 109 109 \dt 
\plt 19 110 1 \dt 
\plt 18 110 22 \dt 
\plt 18 110 43 \dt 
\plt 17 110 64 \dt 
\plt 17 110 85 \dt 
\plt 16 110 106 \dt 
\plt 16 110 127 \dt 
\plt 14 110 149 \dt 
\plt 14 110 170 \dt 
\plt 13 110 191 \dt 
\plt 8 110 212 \dt 
\plt 112 111 1 \dt 
\plt 0 111 115 \dt 
\plt 0 111 117 \dt 
\plt 0 111 119 \dt 
\plt 0 111 121 \dt 
\plt 0 111 123 \dt 
\plt 0 111 125 \dt 
\plt 33 112 1 \dt 
\plt 15 112 36 \dt 
\plt 14 112 54 \dt 
\plt 13 112 72 \dt 
\plt 13 112 89 \dt 
\plt 13 112 107 \dt 
\plt 12 112 125 \dt 
\plt 11 112 143 \dt 
\plt 10 112 161 \dt 
\plt 8 112 180 \dt 
\plt 8 112 197 \dt 
\plt 4 112 216 \dt 
\plt 122 113 1 \dt 
\plt 1 113 125 \dt 
\plt 0 113 128 \dt 
\plt 1 113 130 \dt 
\plt 0 113 133 \dt 
\plt 0 113 136 \dt 
\plt 0 113 138 \dt 
\plt 0 113 141 \dt 
\plt 13 114 1 \dt 
\plt 13 114 16 \dt 
\plt 13 114 31 \dt 
\plt 12 114 47 \dt 
\plt 11 114 62 \dt 
\plt 10 114 78 \dt 
\plt 10 114 93 \dt 
\plt 9 114 109 \dt 
\plt 8 114 124 \dt 
\plt 7 114 140 \dt 
\plt 7 114 155 \dt 
\plt 3 114 174 \dt 
\plt 1 114 189 \dt 
\plt 0 114 192 \dt 
\plt 150 115 1 \dt 
\plt 5 115 153 \dt 
\plt 1 115 160 \dt 
\plt 2 115 163 \dt 
\plt 1 115 167 \dt 
\plt 0 115 171 \dt 
\plt 1 115 174 \dt 
\plt 0 115 178 \dt 
\plt 0 115 185 \dt 
\plt 11 116 1 \dt 
\plt 11 116 14 \dt 
\plt 10 116 28 \dt 
\plt 10 116 41 \dt 
\plt 9 116 55 \dt 
\plt 8 116 69 \dt 
\plt 7 116 83 \dt 
\plt 7 116 96 \dt 
\plt 6 116 110 \dt 
\plt 4 116 125 \dt 
\plt 2 116 140 \dt 
\plt 0 116 153 \dt 
\plt 216 117 1 \dt 
\plt 1 117 219 \dt 
\plt 10 118 1 \dt 
\plt 10 118 13 \dt 
\plt 9 118 25 \dt 
\plt 9 118 37 \dt 
\plt 7 118 50 \dt 
\plt 7 118 62 \dt 
\plt 7 118 74 \dt 
\plt 4 118 88 \dt 
\plt 4 118 100 \dt 
\plt 1 118 113 \dt 
\plt 219 119 1 \dt 
\plt 9 120 1 \dt 
\plt 8 120 12 \dt 
\plt 8 120 23 \dt 
\plt 7 120 34 \dt 
\plt 7 120 45 \dt 
\plt 6 120 56 \dt 
\plt 6 120 67 \dt 
\plt 3 120 80 \dt 
\plt 3 120 91 \dt 
\plt 219 121 1 \dt 
\plt 109 122 1 \dt 
\plt 120 123 1 \dt 
\plt 97 123 123 \dt 
\plt 94 124 1 \dt 
\plt 4 124 99 \dt 
\plt 0 124 110 \dt 
\plt 59 125 1 \dt 
\plt 59 125 62 \dt 
\plt 58 125 124 \dt 
\plt 35 125 185 \dt 
\plt 96 126 1 \dt 
\plt 3 126 99 \dt 
\plt 1 126 105 \dt 
\plt 1 126 110 \dt 
\plt 39 127 1 \dt 
\plt 39 127 42 \dt 
\plt 39 127 83 \dt 
\plt 38 127 125 \dt 
\plt 38 127 166 \dt 
\plt 13 127 207 \dt 
\plt 104 128 1 \dt 
\plt 1 128 107 \dt 
\plt 1 128 110 \dt 
\plt 0 128 114 \dt 
\plt 0 128 121 \dt 
\plt 29 129 1 \dt 
\plt 29 129 32 \dt 
\plt 29 129 63 \dt 
\plt 28 129 94 \dt 
\plt 27 129 126 \dt 
\plt 27 129 157 \dt 
\plt 26 129 189 \dt 
\plt 0 129 220 \dt 
\plt 113 130 1 \dt 
\plt 0 130 116 \dt 
\plt 0 130 119 \dt 
\plt 0 130 121 \dt 
\plt 0 130 124 \dt 
\plt 23 131 1 \dt 
\plt 23 131 26 \dt 
\plt 23 131 51 \dt 
\plt 22 131 76 \dt 
\plt 22 131 101 \dt 
\plt 21 131 127 \dt 
\plt 20 131 152 \dt 
\plt 20 131 177 \dt 
\plt 17 131 203 \dt 
\plt 121 132 1 \dt 
\plt 0 132 124 \dt 
\plt 0 132 126 \dt 
\plt 0 132 128 \dt 
\plt 0 132 130 \dt 
\plt 0 132 132 \dt 
\plt 0 132 139 \dt 
\plt 0 132 141 \dt 
\plt 0 132 143 \dt 
\plt 19 133 1 \dt 
\plt 19 133 22 \dt 
\plt 18 133 43 \dt 
\plt 18 133 64 \dt 
\plt 17 133 86 \dt 
\plt 16 133 107 \dt 
\plt 16 133 128 \dt 
\plt 15 133 150 \dt 
\plt 14 133 171 \dt 
\plt 13 133 193 \dt 
\plt 6 133 214 \dt 
\plt 137 134 1 \dt 
\plt 0 134 140 \dt 
\plt 3 134 142 \dt 
\plt 0 134 147 \dt 
\plt 0 134 149 \dt 
\plt 0 134 151 \dt 
\plt 0 134 154 \dt 
\plt 0 134 156 \dt 
\plt 0 134 158 \dt 
\plt 16 135 1 \dt 
\plt 16 135 19 \dt 
\plt 16 135 37 \dt 
\plt 15 135 56 \dt 
\plt 14 135 74 \dt 
\plt 14 135 92 \dt 
\plt 13 135 111 \dt 
\plt 13 135 129 \dt 
\plt 11 135 148 \dt 
\plt 11 135 166 \dt 
\plt 10 135 185 \dt 
\plt 8 135 205 \dt 
\plt 166 136 1 \dt 
\plt 1 136 169 \dt 
\plt 1 136 172 \dt 
\plt 1 136 175 \dt 
\plt 0 136 178 \dt 
\plt 0 136 181 \dt 
\plt 0 136 184 \dt 
\plt 0 136 187 \dt 
\plt 14 137 1 \dt 
\plt 14 137 17 \dt 
\plt 13 137 33 \dt 
\plt 13 137 49 \dt 
\plt 13 137 65 \dt 
\plt 11 137 82 \dt 
\plt 11 137 98 \dt 
\plt 9 137 115 \dt 
\plt 9 137 131 \dt 
\plt 9 137 147 \dt 
\plt 7 137 164 \dt 
\plt 7 137 180 \dt 
\plt 1 137 197 \dt 
\plt 2 137 200 \dt 
\plt 1 137 217 \dt 
\plt 211 138 1 \dt 
\plt 2 138 214 \dt 
\plt 2 138 218 \dt 
\plt 12 139 1 \dt 
\plt 12 139 15 \dt 
\plt 12 139 29 \dt 
\plt 11 139 44 \dt 
\plt 10 139 59 \dt 
\plt 10 139 73 \dt 
\plt 9 139 88 \dt 
\plt 8 139 103 \dt 
\plt 7 139 118 \dt 
\plt 5 139 134 \dt 
\plt 5 139 148 \dt 
\plt 3 139 162 \dt 
\plt 0 139 167 \dt 
\plt 0 139 181 \dt 
\plt 219 140 1 \dt 
\plt 11 141 1 \dt 
\plt 11 141 14 \dt 
\plt 10 141 27 \dt 
\plt 10 141 40 \dt 
\plt 9 141 53 \dt 
\plt 8 141 67 \dt 
\plt 8 141 80 \dt 
\plt 6 141 94 \dt 
\plt 6 141 107 \dt 
\plt 5 141 120 \dt 
\plt 2 141 136 \dt 
\plt 219 142 1 \dt 
\plt 10 143 1 \dt 
\plt 9 143 13 \dt 
\plt 9 143 25 \dt 
\plt 8 143 37 \dt 
\plt 8 143 49 \dt 
\plt 7 143 61 \dt 
\plt 7 143 73 \dt 
\plt 5 143 86 \dt 
\plt 5 143 98 \dt 
\plt 1 143 113 \dt 
\plt 1 143 125 \dt 
\plt 219 144 1 \dt 
\plt 119 145 1 \dt 
\plt 143 146 1 \dt 
\plt 74 146 146 \dt 
\plt 111 147 1 \dt 
\plt 6 147 114 \dt 
\plt 2 147 126 \dt 
\plt 71 148 1 \dt 
\plt 70 148 74 \dt 
\plt 70 148 147 \dt 
\plt 0 148 220 \dt 
\plt 120 149 1 \dt 
\plt 2 149 123 \dt 
\plt 1 149 129 \dt 
\plt 0 149 135 \dt 
\plt 47 150 1 \dt 
\plt 46 150 50 \dt 
\plt 46 150 99 \dt 
\plt 45 150 148 \dt 
\plt 23 150 197 \dt 
\plt 123 151 1 \dt 
\plt 2 151 126 \dt 
\plt 1 151 130 \dt 
\plt 1 151 134 \dt 
\plt 0 151 138 \dt 
\plt 0 151 142 \dt 
\plt 35 152 1 \dt 
\plt 34 152 38 \dt 
\plt 34 152 75 \dt 
\plt 33 152 112 \dt 
\plt 33 152 149 \dt 
\plt 32 152 186 \dt 
\plt 128 153 1 \dt 
\plt 1 153 131 \dt 
\plt 1 153 134 \dt 
\plt 1 153 137 \dt 
\plt 0 153 140 \dt 
\plt 0 153 143 \dt 
\plt 0 153 146 \dt 
\plt 57 154 1 \dt 
\plt 27 154 60 \dt 
\plt 27 154 90 \dt 
\plt 26 154 120 \dt 
\plt 25 154 150 \dt 
\plt 25 154 180 \dt 
\plt 10 154 210 \dt 
\plt 143 155 1 \dt 
\plt 0 155 146 \dt 
\plt 0 155 148 \dt 
\plt 1 155 150 \dt 
\plt 0 155 153 \dt 
\plt 0 155 155 \dt 
\plt 0 155 160 \dt 
\plt 0 155 162 \dt 
\plt 23 156 1 \dt 
\plt 22 156 26 \dt 
\plt 22 156 51 \dt 
\plt 21 156 76 \dt 
\plt 21 156 101 \dt 
\plt 20 156 126 \dt 
\plt 20 156 151 \dt 
\plt 19 156 176 \dt 
\plt 19 156 201 \dt 
\plt 166 157 1 \dt 
\plt 0 157 169 \dt 
\plt 0 157 171 \dt 
\plt 0 157 173 \dt 
\plt 0 157 175 \dt 
\plt 0 157 177 \dt 
\plt 0 157 179 \dt 
\plt 0 157 181 \dt 
\plt 0 157 183 \dt 
\plt 41 158 1 \dt 
\plt 19 158 44 \dt 
\plt 18 158 66 \dt 
\plt 18 158 87 \dt 
\plt 17 158 109 \dt 
\plt 16 158 131 \dt 
\plt 16 158 152 \dt 
\plt 15 158 174 \dt 
\plt 14 158 196 \dt 
\plt 2 158 218 \dt 
\plt 177 159 1 \dt 
\plt 3 159 180 \dt 
\plt 3 159 185 \dt 
\plt 0 159 190 \dt 
\plt 0 159 193 \dt 
\plt 0 159 195 \dt 
\plt 0 159 198 \dt 
\plt 0 159 200 \dt 
\plt 0 159 205 \dt 
\plt 0 159 210 \dt 
\plt 17 160 1 \dt 
\plt 16 160 20 \dt 
\plt 16 160 39 \dt 
\plt 15 160 58 \dt 
\plt 15 160 77 \dt 
\plt 14 160 96 \dt 
\plt 14 160 115 \dt 
\plt 13 160 134 \dt 
\plt 12 160 154 \dt 
\plt 11 160 173 \dt 
\plt 11 160 192 \dt 
\plt 9 160 211 \dt 
\plt 213 161 1 \dt 
\plt 1 161 216 \dt 
\plt 1 161 219 \dt 
\plt 15 162 1 \dt 
\plt 14 162 18 \dt 
\plt 14 162 35 \dt 
\plt 13 162 52 \dt 
\plt 13 162 69 \dt 
\plt 12 162 86 \dt 
\plt 12 162 103 \dt 
\plt 10 162 121 \dt 
\plt 10 162 138 \dt 
\plt 9 162 155 \dt 
\plt 8 162 173 \dt 
\plt 6 162 191 \dt 
\plt 6 162 208 \dt 
\plt 219 163 1 \dt 
\plt 13 164 1 \dt 
\plt 13 164 16 \dt 
\plt 13 164 31 \dt 
\plt 12 164 47 \dt 
\plt 12 164 62 \dt 
\plt 11 164 78 \dt 
\plt 11 164 93 \dt 
\plt 10 164 109 \dt 
\plt 8 164 125 \dt 
\plt 7 164 141 \dt 
\plt 6 164 157 \dt 
\plt 5 164 173 \dt 
\plt 2 164 191 \dt 
\plt 1 164 207 \dt 
\plt 219 165 1 \dt 
\plt 12 166 1 \dt 
\plt 12 166 15 \dt 
\plt 11 166 29 \dt 
\plt 11 166 43 \dt 
\plt 10 166 57 \dt 
\plt 10 166 71 \dt 
\plt 9 166 86 \dt 
\plt 8 166 100 \dt 
\plt 7 166 115 \dt 
\plt 5 166 130 \dt 
\plt 4 166 145 \dt 
\plt 1 166 160 \dt 
\plt 219 167 1 \dt 
\plt 11 168 1 \dt 
\plt 10 168 14 \dt 
\plt 10 168 27 \dt 
\plt 9 168 40 \dt 
\plt 9 168 53 \dt 
\plt 8 168 66 \dt 
\plt 8 168 79 \dt 
\plt 6 168 93 \dt 
\plt 6 168 106 \dt 
\plt 4 168 120 \dt 
\plt 4 168 133 \dt 
\plt 219 169 1 \dt 
\plt 155 170 1 \dt 
\plt 168 171 1 \dt 
\plt 49 171 171 \dt 
\plt 138 172 1 \dt 
\plt 4 172 143 \dt 
\plt 0 172 156 \dt 
\plt 83 173 1 \dt 
\plt 83 173 86 \dt 
\plt 48 173 172 \dt 
\plt 140 174 1 \dt 
\plt 3 174 143 \dt 
\plt 1 174 150 \dt 
\plt 1 174 156 \dt 
\plt 55 175 1 \dt 
\plt 55 175 58 \dt 
\plt 55 175 115 \dt 
\plt 47 175 173 \dt 
\plt 145 176 1 \dt 
\plt 2 176 148 \dt 
\plt 1 176 152 \dt 
\plt 1 176 156 \dt 
\plt 0 176 161 \dt 
\plt 0 176 165 \dt 
\plt 41 177 1 \dt 
\plt 41 177 44 \dt 
\plt 41 177 87 \dt 
\plt 40 177 130 \dt 
\plt 39 177 174 \dt 
\plt 3 177 217 \dt 
\plt 157 178 1 \dt 
\plt 0 178 160 \dt 
\plt 0 178 163 \dt 
\plt 0 178 166 \dt 
\plt 0 178 169 \dt 
\plt 67 179 1 \dt 
\plt 32 179 70 \dt 
\plt 32 179 105 \dt 
\plt 31 179 140 \dt 
\plt 30 179 175 \dt 
\plt 10 179 210 \dt 
\plt 159 180 1 \dt 
\plt 3 180 162 \dt 
\plt 3 180 167 \dt 
\plt 0 180 172 \dt 
\plt 0 180 175 \dt 
\plt 0 180 177 \dt 
\plt 0 180 180 \dt 
\plt 0 180 182 \dt 
\plt 27 181 1 \dt 
\plt 27 181 30 \dt 
\plt 26 181 59 \dt 
\plt 26 181 88 \dt 
\plt 25 181 118 \dt 
\plt 24 181 147 \dt 
\plt 24 181 176 \dt 
\plt 15 181 205 \dt 
\plt 182 182 1 \dt 
\plt 0 182 185 \dt 
\plt 0 182 187 \dt 
\plt 0 182 189 \dt 
\plt 0 182 191 \dt 
\plt 0 182 193 \dt 
\plt 0 182 195 \dt 
\plt 23 183 1 \dt 
\plt 23 183 26 \dt 
\plt 22 183 51 \dt 
\plt 22 183 76 \dt 
\plt 22 183 101 \dt 
\plt 20 183 127 \dt 
\plt 20 183 152 \dt 
\plt 20 183 177 \dt 
\plt 18 183 202 \dt 
\plt 200 184 1 \dt 
\plt 0 184 203 \dt 
\plt 0 184 205 \dt 
\plt 1 184 207 \dt 
\plt 0 184 210 \dt 
\plt 0 184 212 \dt 
\plt 0 184 214 \dt 
\plt 0 184 216 \dt 
\plt 20 185 1 \dt 
\plt 20 185 23 \dt 
\plt 19 185 45 \dt 
\plt 19 185 67 \dt 
\plt 19 185 89 \dt 
\plt 18 185 111 \dt 
\plt 17 185 134 \dt 
\plt 16 185 156 \dt 
\plt 15 185 179 \dt 
\plt 15 185 201 \dt 
\plt 219 186 1 \dt 
\plt 37 187 1 \dt 
\plt 17 187 40 \dt 
\plt 17 187 60 \dt 
\plt 16 187 80 \dt 
\plt 15 187 100 \dt 
\plt 14 187 120 \dt 
\plt 14 187 140 \dt 
\plt 13 187 160 \dt 
\plt 12 187 180 \dt 
\plt 11 187 200 \dt 
\plt 0 187 220 \dt 
\plt 219 188 1 \dt 
\plt 33 189 1 \dt 
\plt 16 189 36 \dt 
\plt 15 189 54 \dt 
\plt 14 189 72 \dt 
\plt 13 189 91 \dt 
\plt 13 189 108 \dt 
\plt 12 189 127 \dt 
\plt 11 189 145 \dt 
\plt 10 189 163 \dt 
\plt 9 189 182 \dt 
\plt 9 189 199 \dt 
\plt 1 189 219 \dt 
\plt 219 190 1 \dt 
\plt 14 191 1 \dt 
\plt 14 191 17 \dt 
\plt 14 191 33 \dt 
\plt 13 191 50 \dt 
\plt 13 191 66 \dt 
\plt 12 191 83 \dt 
\plt 10 191 100 \dt 
\plt 10 191 116 \dt 
\plt 9 191 133 \dt 
\plt 8 191 150 \dt 
\plt 7 191 167 \dt 
\plt 6 191 184 \dt 
\plt 5 191 201 \dt 
\plt 3 191 217 \dt 
\plt 219 192 1 \dt 
\plt 13 193 1 \dt 
\plt 13 193 16 \dt 
\plt 12 193 31 \dt 
\plt 12 193 46 \dt 
\plt 11 193 61 \dt 
\plt 11 193 76 \dt 
\plt 9 193 92 \dt 
\plt 9 193 107 \dt 
\plt 9 193 122 \dt 
\plt 6 193 139 \dt 
\plt 6 193 154 \dt 
\plt 5 193 169 \dt 
\plt 3 193 186 \dt 
\plt 219 194 1 \dt 
\plt 12 195 1 \dt 
\plt 11 195 15 \dt 
\plt 11 195 29 \dt 
\plt 10 195 43 \dt 
\plt 10 195 57 \dt 
\plt 9 195 71 \dt 
\plt 9 195 85 \dt 
\plt 7 195 100 \dt 
\plt 7 195 114 \dt 
\plt 5 195 129 \dt 
\plt 5 195 143 \dt 
\plt 1 195 160 \dt 
\plt 1 195 174 \dt 
\plt 219 196 1 \dt 
\plt 167 197 1 \dt 
\plt 195 198 1 \dt 
\plt 22 198 198 \dt 
\plt 158 199 1 \dt 
\plt 7 199 161 \dt 
\plt 3 199 175 \dt 
\plt 97 200 1 \dt 
\plt 96 200 100 \dt 
\plt 21 200 199 \dt 
\plt 162 201 1 \dt 
\plt 4 201 165 \dt 
\plt 2 201 172 \dt 
\plt 1 201 179 \dt 
\plt 64 202 1 \dt 
\plt 64 202 67 \dt 
\plt 64 202 133 \dt 
\plt 20 202 200 \dt 
\plt 168 203 1 \dt 
\plt 2 203 171 \dt 
\plt 2 203 175 \dt 
\plt 1 203 180 \dt 
\plt 0 203 185 \dt 
\plt 0 203 189 \dt 
\plt 48 204 1 \dt 
\plt 47 204 51 \dt 
\plt 47 204 101 \dt 
\plt 46 204 151 \dt 
\plt 19 204 201 \dt 
\plt 178 205 1 \dt 
\plt 1 205 181 \dt 
\plt 1 205 184 \dt 
\plt 0 205 188 \dt 
\plt 0 205 191 \dt 
\plt 0 205 195 \dt 
\plt 0 205 198 \dt 
\plt 38 206 1 \dt 
\plt 38 206 41 \dt 
\plt 38 206 81 \dt 
\plt 37 206 121 \dt 
\plt 37 206 161 \dt 
\plt 18 206 202 \dt 
\plt 184 207 1 \dt 
\plt 3 207 187 \dt 
\plt 1 207 192 \dt 
\plt 1 207 195 \dt 
\plt 0 207 198 \dt 
\plt 0 207 201 \dt 
\plt 0 207 203 \dt 
\plt 0 207 206 \dt 
\plt 0 207 209 \dt 
\plt 65 208 1 \dt 
\plt 31 208 68 \dt 
\plt 30 208 102 \dt 
\plt 30 208 135 \dt 
\plt 29 208 169 \dt 
\plt 17 208 203 \dt 
\plt 202 209 1 \dt 
\plt 0 209 205 \dt 
\plt 0 209 207 \dt 
\plt 1 209 209 \dt 
\plt 0 209 212 \dt 
\plt 0 209 214 \dt 
\plt 0 209 216 \dt 
\plt 27 210 1 \dt 
\plt 26 210 30 \dt 
\plt 26 210 59 \dt 
\plt 25 210 88 \dt 
\plt 25 210 117 \dt 
\plt 24 210 146 \dt 
\plt 24 210 175 \dt 
\plt 16 210 204 \dt 
\plt 219 211 1 \dt 
\plt 23 212 1 \dt 
\plt 23 212 26 \dt 
\plt 22 212 52 \dt 
\plt 22 212 77 \dt 
\plt 21 212 103 \dt 
\plt 21 212 128 \dt 
\plt 20 212 154 \dt 
\plt 19 212 180 \dt 
\plt 15 212 205 \dt 
\plt 219 213 1 \dt 
\plt 43 214 1 \dt 
\plt 20 214 46 \dt 
\plt 20 214 69 \dt 
\plt 19 214 92 \dt 
\plt 18 214 115 \dt 
\plt 17 214 138 \dt 
\plt 17 214 161 \dt 
\plt 16 214 184 \dt 
\plt 13 214 207 \dt 
\plt 219 215 1 \dt 
\plt 39 216 1 \dt 
\plt 18 216 42 \dt 
\plt 17 216 63 \dt 
\plt 17 216 83 \dt 
\plt 16 216 104 \dt 
\plt 15 216 125 \dt 
\plt 14 216 146 \dt 
\plt 14 216 166 \dt 
\plt 13 216 187 \dt 
\plt 12 216 208 \dt 
\plt 219 217 1 \dt 
\plt 35 218 1 \dt 
\plt 16 218 38 \dt 
\plt 16 218 57 \dt 
\plt 15 218 76 \dt 
\plt 14 218 95 \dt 
\plt 14 218 114 \dt 
\plt 13 218 133 \dt 
\plt 12 218 152 \dt 
\plt 12 218 171 \dt 
\plt 10 218 191 \dt 
\plt 10 218 209 \dt 
\plt 219 219 1 \dt 
\plt 15 220 1 \dt 
\plt 15 220 18 \dt 
\plt 15 220 35 \dt 
\plt 14 220 53 \dt 
\plt 14 220 70 \dt 
\plt 12 220 88 \dt 
\plt 12 220 105 \dt 
\plt 11 220 123 \dt 
\plt 11 220 140 \dt 
\plt 10 220 158 \dt 
\plt 9 220 176 \dt 
\plt 8 220 193 \dt 
\plt 7 220 211 \dt 
}}

\vfil\eject

\References

\bibitem{\refAH} J. Alexander and A. Hirschowitz. {\it An asymptotic
vanishing theorem for generic unions of multiple points},
Invent. Math. 140 (2000), no. 2, 303--325.

\bibitem{\refB} E. Ballico. {\it Curves of minimal degree
with prescribed singularities}, Illinois J.\ Math.\ 45 (1999), 672--676.

\bibitem{\refBZ} A. Buckley and M. Zompatori. {\it Linear systems of plane curves with
a composite number of base points of equal multiplicity}, preprint (2001).

\bibitem{\refCat} M. V. Catalisano. {\it Linear Systems of Plane Curves
through Fixed ``Fat'' Points of \pr2},
J.\ Alg.\  142 (1991), 81-100.

\bibitem{\refCatb} M. V. Catalisano. {\it ``Fat'' points on a conic},
Comm.\ Alg.\  19(8) (1991), 2153--2168.

\bibitem{\refCCMO} C. Ciliberto, F. Cioffi, R. Miranda and
F. Orecchia. {\it Bivariate Hermite interpolation via computer algebra
and algebraic geometry techniques.},
Preprint (2001).

\bibitem{\refCM} C. Ciliberto and R. Miranda. {\it Degenerations
of planar linear systems},
Journ.\ Reine Angew.\ Math. 501 (1998), 191--220.

\bibitem{\refCMb} C. Ciliberto and R. Miranda. {\it Linear systems
of plane curves with base points of equal multiplicity},
Trans.\ Amer. Math. Soc. 352 (2000), 4037--4050.

\bibitem{\refCMc}  C. Ciliberto and R. Miranda. {\it
The Segre and Harbourne-Hirschowitz conjectures}, in:  
Applications of algebraic geometry to coding theory, 
physics and computation (Eilat, 2001),   37--51,  
NATO Sci. Ser. II Math. Phys. Chem., 36,  
Kluwer Acad. Publ., Dordrecht,  2001. 

\bibitem{\refE} L. \'Evain. {\it Une minoration du degre des courbes
planes \`a singularit\'es impos\'ees}, Bull.\ Soc.\ Math.\
France 126 (1998), no. 4, 525--543.

\bibitem{\refEb}
L. \'Evain. {\it La fonction de {H}ilbert de la r{\'e}union de $4^h$
gros points g{\'e}n{\'e}riques de \pr2 de m{\^e}me multiplicit{\'e}},
J.\ Algebraic Geometry (1999), 787--796.

\bibitem{\refFHH} S. Fitchett, B. Harbourne and S. Holay.
{\it Resolutions of Fat Point Ideals involving Eight General Points of \pr2},
J. Algebra 244 (2001), 684--705.

\bibitem{\refGGR} A. V. Geramita, D. Gregory, L. Roberts.
{\it Monomial ideals and points in projective space},
J.\ Pure Appl.\ Algebra 40 (1986), 33--62.

\bibitem{\refGi} A. Gimigliano.
{\it Regularity of Linear Systems of Plane Curves},
J.\ Alg. 124 (1989), 447--460.

\bibitem{\refGIda} A. Gimigliano and M. Id\`a.
{\it The ideal resolution for generic 3-fat points in \pr2},
preprint (2002).

\bibitem{\refvanc} B. Harbourne. {\it The geometry of rational surfaces
and Hilbert functions of points in the plane},
Can.\ Math.\ Soc.\ Conf.\ Proc.\ 6
(1986), 95--111.

\bibitem{\refigp} B. Harbourne. {\it  The Ideal Generation
Problem for Fat Points},
J.\ Pure Appl.\ Algebra 145 (2000), 165--182.

\bibitem{\refnagconj} B. Harbourne.
{\it On Nagata's Conjecture}, J. Alg. 236 (2001), 692--702.

\bibitem{\refsurv} B. Harbourne.
{\it Problems and Progress: A survey on fat points
in \pr2}, v. 123, 2002, Queen's
papers in pure and applied mathematics,
The curves seminar at Queen's.

\bibitem{\refSC} B. Harbourne. {\it  Seshadri constants and very ample
divisors on algebraic surfaces},
Journ.\ Reine Angew.\ Math. (to appear).

\bibitem{\refHHF} B. Harbourne, S. Holay and S. Fitchett.
{\it Resolutions of ideals of quasiuniform fat point subschemes of 
\pr2}, Trans. Amer. Math. Soc. 355 (2003), no. 2, 593--608.

\bibitem{\refHi} A. Hirschowitz.
{\it Une conjecture pour la cohomologie
des diviseurs sur les surfaces rationelles g\'en\'e\-ri\-ques},
Journ.\ Reine Angew.\ Math. 397
(1989), 208--213.

\bibitem{\refHib} A. Hirschowitz.
{\it La m\'ethode d'Horace pour l'interpolation \`a plusieurs
variables}, Manus. Math. 50 (1985), 337--388.

\bibitem{\refHo} M. Homma.
{\it A souped up version of Pardini's theorem and its aplication
to funny curves}, Comp.\ Math. 71 (1989), 295--302.

\bibitem{\refI} M. Id\`a.
{\it The minimal free resolution for the first infinitesimal
neighborhoods of $n$ general points in the plane},
J.\ Alg. 216 (1999), 741--753.


\bibitem{\refMig}
T. Mignon. {\it Syst\`emes de courbes planes \`a singularit\'es
impos\'ees: le cas des multiplicit\'es inf\'erieures ou \'egales
\`a quatre}, J.\ Pure Appl.\ Algebra 151 (2000), no. 2, 173--195.

\bibitem{\refMirb}
R. Miranda. {\it Algebraic Curves and Riemann Surfaces},
Graduate Studies in Mathematics 5, Amer. Math. Soc.,
(1995), xxi $+$ 390 pp.


\bibitem{\refN} M. Nagata. {\it On the 14-th problem of Hilbert},
Amer.\ J.\ Math.\ 81 (1959), 766--772.

\bibitem{\refNb} M. Nagata. {\it On rational surfaces, II}, Mem.\ Coll.\
Sci.\ Univ.\ Kyoto, Ser.\ A Math.\ 33 (1960), 271--293.

\bibitem{\refR} J. Ro\'e. {\it On the existence of plane curves
with imposed multiple points}, J. Pure Appl. Alg. 156 (2001), 115--126.

\bibitem{\refRb} J. Ro\'e. {\it Linear systems of plane curves
with imposed multiple points}, Illinois J. Math. 45 (2001), 895--906.

\bibitem{\refS} B. Segre. {\it Alcune questioni su insiemi finiti
di punti in Geometria Algebrica},
Atti del Convegno Internaz. di Geom. Alg., Torino (1961).

\bibitem{\refXb}
G. Xu. {\it Ample line bundles on smooth surfaces},
Journ.\ Reine Angew.\ Math.
{469} (1995), 199--209.


\bye